%% file: arxivthese.tex
\documentclass[fleqn,12pt,final]{book}

\usepackage{amsmath} 
\usepackage{amssymb,latexsym} 

\newcommand{\emptyevenpage}{}

\newcommand{\Ccal}{\mathcal{C}}

\let\ol=\overline
\let\fr=\partial
\let\Om=\Omega
\let\ck=\check

\newcommand{\Y}{\mathcal{Y}}
\newcommand{\K}{\mathcal{K}}
\let\d=\delta
\let\e=\varepsilon
\let\phi=\varphi
\let\sub=\subseteq
\let\w=\omega
\let\lra=\longrightarrow

\newcommand{\R}{\ensuremath{\mathbb{R}}}
\newcommand{\N}{\ensuremath{\mathbb{N}}}
\newcommand{\C}{\ensuremath{\mathbb{C}}}
\newcommand{\Q}{\ensuremath{\mathbb{Q}}}

\let\oldS=\S

\renewcommand{\S} {\ensuremath{\mathcal{S}}} 
\renewcommand{\P} {\ensuremath{\mathcal{P}}}
\newcommand{\U} {\ensuremath{\mathcal{U}}}
\newcommand{\X} {\ensuremath{\mathcal{X}}}

\newcommand{\Frac}[2]{\frac{\displaystyle{#1}}{\displaystyle{#2}}}

\newcommand{\DP}[2]{\frac{\partial #1}{\partial #2}} 

\newcommand{\rk}{\mathrm{rank\,}}
\let\gr=\nabla

\newtheorem{theorem}{Theorem}
\newtheorem{lemma}[theorem]{Lemma}
\newtheorem{corollary}[theorem]{Corollary}
\newtheorem{definition}[theorem]{Definition}
\newtheorem{proposition}[theorem]{Proposition}
\newtheorem{example}[theorem]{Example}
\newtheorem{rem}[theorem]{Remark}
\newtheorem{notation}[theorem]{Notation}

{ \medskip \noindent \rmfamily
 {\it Remark. } }{ \medskip}

 \newenvironment{proof}%
 { \medskip
   \rmfamily \noindent
   {\bf Proof:\/} } {\hfill $\Box$ \medskip}

\newcommand {\iso}         {\cong}

\newcommand{\V}{\mathcal{V}}

\newcommand{\Sig} {\mathfrak{S}}
\newcommand{\s} {\sigma}
\let\fr=\partial
\let\w=\omega

\let\ra=\rightarrow
\newcommand{\Z}{\mathbb{Z}}

\newcommand{\dist}{\mathrm{dist}}
\newcommand{\T}{\mathcal{T}}

\newcommand{\Ecal}{\mathcal{E}}

\let\bs=\backslash
\newcommand{\f}{\boldsymbol{f}}
\newcommand{\n}{\boldsymbol{n}}

\newcommand{\Bcal}{\mathcal{B}}

\let\inc=\hookrightarrow

\newcommand{\xx}{\mathbf{x}}
\newcommand{\y}{\mathbf{y}}

\newcommand{\Zcal}{\mathcal{Z}}

\let\mpty=\varnothing
\newcommand{\W}{\mathcal{W}}
\newcommand{\NQ}{\widetilde{Q}}
\newcommand{\sib}{\bar{\sigma}}
\newcommand{\A}{\mathcal{A}}
\newcommand{\I}{\mathcal{I}}
\newcommand{\HB}{\bar{H}}
\newcommand{\HT}{\widetilde{H}}
\newcommand{\HBT}{\widetilde{\bar{H}}{}}
\newcommand{\bbm}{b^{\mbox{\tiny{BM}}}}
\newcommand{\Hbm}{H^{\mbox{\tiny{BM}}}}
\newcommand{\limind}{\mathop{\lim}_{\lra}}
\newcommand{\sign}{\mathrm{sign}}
\newcommand{\J}{\mathcal{J}}

\usepackage{xypic}

\begin{document}

\begin{titlepage}
\vfill

\begin{center}
{\huge Quantitative study of semi-Pfaffian sets}\\ 

\vfill

{\large \'Etude quantitative des ensembles semi-pfaffiens}

\vfill
{\large Thierry Zell \\ December 2003}

\bigskip

03C98 (Secondary).

\vfill
School of Mathematics, Georgia Institute of Technology \\
686 Cherry Street, Atlanta GA 30332-0160\\
{\tt zell}@{\tt math.gatech.edu}
\end{center}

\end{titlepage}


\input{abstract}

\eject
\thispagestyle{empty}

\tableofcontents

\cleardoublepage

\input{index}

\input{ack}

\pagestyle{myheadings}
\markboth{ \hfill INTRODUCTION}{INTRODUCTION \hfill}

\chapter*{Introduction}
\addcontentsline{toc}{chapter}{Introduction}
\input{intro}

\pagestyle{headings}
\chapter[Preliminaries]{Preliminaries}

\input{aa1}

\chapter[Betti numbers of semi-Pfaffian sets]{Betti numbers of semi-Pfaffian sets}
\setcounter{theorem}{0}
\input{aa2}

\chapter[Betti numbers of sub-Pfaffian sets]{Betti numbers of sub-Pfaffian sets}
\setcounter{theorem}{0}
\input{aa3}

\chapter[Connected components of limit sets]{Connected components of limit sets}
\setcounter{theorem}{0}
\input{aa4}

\chapter[Topology of Hausdorff limits]{Topology of Hausdorff limits}
\setcounter{theorem}{0}
\input{aa5}

\appendix
\setcounter{theorem}{0}
\input{spectral}

\input{bibliographie}

\eject 

\end{document}

%% file: abstract.tex

\hbox{ }
\thispagestyle{empty}

\newpage
\thispagestyle{empty}
\begin{center}
{\large \textsc{\'Etude quantitative des ensembles semi-pfaffiens} 

\medskip

{\bf R\'esum\'e}}
\end{center}

  Dans la pr\'esente th\`ese, on \'etablit des bornes sup\'erieures
  sur les nombres de Betti des ensembles d\'efinis \`a l'aide de
  fonctions pfaffiennes, en fonction de la complexit\'e pfaffienne
  (ou format) de ces ensembles.

  \smallskip

  Les fonctions pfaffiennes ont \'et\'e d\'efinies par Khovanskii,
  comme solutions au comportement quasi-polynomial de certains
  syst\`emes polynomiaux d'\'equations diff\'erentielles. Les
  ensembles semi-pfaffiens satisfont une condition de signe bool\'eene
  sur des fonctions pfaffiennes, et les ensembles sous-pfaffiens sont
  projections de semi-pfaffiens. Wilkie a d\'emontr\'e que
  les fonctions pfaffiennes engendrent une structure o-minimale, et
  Gabrielov a montr\'e que cette structure pouvait \^etre efficacement
  d\'ecrite par des ensembles pfaffiens limites.

  \smallskip

  \`A l'aide de la th\'eorie de Morse, de d\'eformations, de
  recurrences sur le niveau combinatoire et de suites spectrales, on
  donne dans cette th\`ese des bornes effectives pour  toutes les
  cat\'egories d'ensembles pr\'e-cit\'ees.

\vfill

\begin{center}
 {\large \textsc{Quantitative study of semi-Pfaffian sets}

\medskip

{\bf Abstract}}
\end{center}

  In the present thesis, we establish upper-bounds on the Betti
  numbers of sets defined using Pfaffian functions, in terms of the
  natural Pfaffian complexity (or format) of those sets.

  \smallskip

  Pfaffian functions were introduced by Khovanskii, as solutions of
  certain polynomial differential systems that have polynomial-like
  behaviour over the real domain. Semi-Pfaffian sets are sets that
  satisfy a quantifier-free sign condition on such functions, and
  sub-Pfaffian sets are linear projection of semi-Pfaffian sets.
  Wilkie showed that Pfaffian functions generate an o-minimal
  structure, and Gabrielov showed that this structure could be
  effectively described by Pfaffian limit sets.

  \smallskip

  Using Morse theory, deformations, recursion on combinatorial levels
  and a spectral sequence associated to continuous surjections,
  we give in this thesis effective estimates for sets belonging to all
  of the above classes.

\vfill

%% file: index.tex

\eject

\thispagestyle{plain}

\begin{center}
{\large {\bf Index of notations}}
\end{center}
\addcontentsline{toc}{chapter}{Index of notations}

\bigskip

\begin{tabular}{ll}
   \hspace{3cm}   &  \cr

   $\N$  \dotfill &    the natural numbers. $\N=\{0,1,2,\cdots\}$ \cr

   $\R$  \dotfill &    the field of real numbers\cr

   $\f$  \dotfill &    a Pfaffian chain $\f=(f_1, \ldots, f_\ell)$ \cr

   $\ell$ \dotfill &   the length of $\f$\cr

   $\alpha$ \dotfill & the degree of $\f$\cr

   $\U$  \dotfill &  a domain where $\f$ is defined \cr

   $\gamma$ \dotfill & the degree of $\U,$ a measure of its topological complexity \cr

   $q$ \dotfill & a Pfaffian function in the chain $\f;$ 
$q(x)=Q(x,f_1(x),\ldots,f_\ell(x))$\cr

   $\deg_{\f} q$ \dotfill & the degree of $q$ in the chain $\f$\cr

   $\P$ \dotfill &     a finite collection of Pfaffian functions\cr

   $ \beta$ \dotfill & a bound on $\deg_{\f} q$ for $q \in \P$\cr

   $|x|$\dotfill &    the Euclidean norm of $x$ \cr

   $\ol{A}$ \dotfill & the closure of $A$ for the Euclidean topology \cr
  
   $ \fr A$ \dotfill & the frontier of $A,$ $\fr A=\ol{A} \bs A$\cr

   $\Gamma(f)$ \dotfill & the graph of $f:$  $\Gamma(f)=\{(x,y) \mid y=f(x)\}$ \cr

   $\lambda\ll 1$\dotfill &    abbreviates $\exists \lambda_0, \, 
\forall \lambda \in (0, \lambda_0)\ldots$ \cr

   $H_i(X)$ \dotfill & the $i$-th singular homology group of $X,$ with 
coefficients in $\Z$  \cr  

   $b_i(X)$ \dotfill & the $i$-th Betti number of $X,$ $b_i(X)=\rk H_i(X)$\cr 

   $b(X)$ \dotfill & the sum of the Betti numbers of $X$\cr 

   $X_\lambda$\dotfill &  the fiber of $X$ at $\lambda:$ 
$X_\lambda=\{x \mid (x,\lambda) \in X\}$  \cr
\end{tabular}

%% file: ack.tex

\chapter*{Acknowledgments}
\addcontentsline{toc}{chapter}{Acknowledgments}

This document is the final version of the my Ph.D. thesis for
archival on {\tt arXiv.org}, under the reference math.AG/0401079.

\bigskip 

This work was undertaken under the joint direction of Andrei Gabrielov
in Purdue University and Marie-Fran\c{c}oise Roy in the Universit\'e
de Rennes~1. It goes without saying that this document owes a lot to
both of them, and I am very grateful to both of them for their
guidance and their patience.

\bigskip

My committee in Purdue was composed of Alex Eremenko, Leonard
Lipshitz, Andrei Gabrielov and Kenji Matsuki. Krysztof Kurdyka
(Universit\'e de Savoie) and Jean-Philippe Rolin (Universit\'e de
Bourgogne) refereed the thesis for the defense in Rennes and sat on
the defense jury, along with Michel Coste, Jean-Marie Lion and
Marie-Fran\c{c}oise Roy from the Universit\'e de Rennes~1.
Many thanks go to all of them for their help.

\bigskip

In addition to the people mentioned above, I'd also like to thank
Ricky Pollack, Saugata Basu, Dan Richardson, Mohab Safey El Din,
\'Eric Schost and Patrick Speissegger for their encouragements and
opportunities for fruitful discussions.

\bigskip

\begin{flushright}
Atlanta, Georgia\\
January 8, 2004.
\end{flushright}

\cleardoublepage

%% file: intro.tex


\begin{quote}
{\it 
Je n'ai jamais \'et\'e assez loin pour bien sentir l'application de
l'alg\`ebre à la g\'eom\'etrie. Je n'aimais point cette mani\`ere
d'op\'erer sans voir ce qu'on fait; et il me semblait que r\'esoudre
un probl\`eme de g\'eom\'etrie par les \'equations, c'\'etait jouer un
air en tournant une manivelle.}
\begin{flushright}
J.-J. \textsc{Rousseau} ({\em Les Confessions,} Livre \textsc{VI})
\end{flushright}
\end{quote}

\section*{Origins of Pfaffian functions}
\addcontentsline{toc}{section}{Origins of Pfaffian functions}

Pfaffian functions were introduced by
Khovanskii~\cite{kh:class,kh:mani,kh:few} in the late seventies. They
are analytic functions which have 
{\bf polynomial-like finiteness properties} 
over the reals. If $\f=(f_1, \ldots, f_\ell)$ is a
sequence of analytic functions defined on a domain $\U \sub \R^n,$ we
will say that $\f$ forms a {\bf Pfaffian chain} on $\U$ if there
exists real polynomials $P_{i,j}$ such that the following {\bf
triangular} differential system holds
\begin{equation*}
df_i(x)=\sum_{j=1}^n P_{i,j}(x, f_1(x), \ldots, f_i(x)) \, dx_j;
\qquad \hbox{for all }1 \leq i \leq \ell.
\end{equation*}
(The system is triangular because for all $i,$ $df_i$ depends only on
the functions $f_1, \ldots, f_i.$)

\medskip

More generally, a {\bf Pfaffian function} is a real analytic function
$q$ that can be written in the form $q(x)=Q(x, f_1(x), \ldots,
f_\ell(x))$ for some polynomial $Q$ and some Pfaffian chain
$(f_1,\ldots, f_\ell).$ They form a large class of functions that
contains, in particular, all Liouvillian functions and all elementary
functions provided they involve sines and cosines that have been
restricted to bounded intervals. In~\cite{kh:few}, a more general
setting involving nested integral manifolds of polynomial 1-forms is
also introduced. This wider setting gives locally the same kind of
functions, so no generality is lost by restricting ourselves to the
first definition.

\medskip

The central result of the theory is the following: any system of $n$
Pfaffian equations in $n$ variables $q_1(x)=\cdots=q_n(x)=0$ has only
finitely many {\bf non-degenerate} real solutions, {\em i.e.}
solutions for which the Jacobian determinant $|\fr q_i / \fr x_j|$
does not vanish.  Moreover, this number can be explicitly bounded from
above in terms of the discrete parameters that define the system. This
result is known as {\bf Khovanskii's theorem}, and it gives Pfaffian
functions their polynomial-like behaviour.  When $q_1, \ldots, q_n$
are polynomials of degree at most $d,$ the number of non-degenerate
solutions of the system $q_1(x)=\cdots=q_n(x)=0$ is bounded by $d^n$
according to the {\bf B\'ezout inequality}.  The idea behind Khovanskii's
theorem is to replace the functions in the Pfaffian chain by variables
to reduce to the polynomial case, using {\bf Rolle's theorem} to build
additional constraints to keep a system with non-degenerate
solutions. 

\medskip

The fact that the Pfaffian chain from which the functions $q_1,
\ldots, q_n$ are derived satisfies a {\em triangular} system is a
crucial point for Khovanskii's theorem to hold. If that restriction is
lifted (giving rise to what is known as a chain of {\bf Noetherian
functions}), one cannot hope to have global finiteness anymore%
\footnote{However, finiteness is preserved on the local level, see~\cite{gakh}.}
, as can be seen from the
example $f_1(x)=\cos x,$ $f_2(x)=\sin x.$ The chain $(f_1,f_2)$ is a
Noetherian chain on the whole real line, since $f_1'=-f_2$ and
$f_2'=f_1,$ but the equation $f_1(x)=0$ has infinitely many
non-degenerate solutions on $\R.$

\bigskip

Khovanskii introduced Pfaffian functions to use in investigations
related to the second part of {\bf Hilbert's sixteenth problem.}  In this
problem, one considers a polynomial vector field in the real plane
given by
\begin{equation*}
\frac{dy}{dx}=\frac{P(x,y)}{Q(x,y)};
\end{equation*}
and Hilbert's original question was to know how many limit cycles
(isolated periodic solutions) could this equation have, and where they
were situated, as a function of the degrees of $P$ and $Q.$ This
problem is still unsolved and has given rise to many related problems,
and Pfaffian functions play an important role in this theory. In
particular, they had a key part in the solution of the so-called
{\bf local Hilbert-Arnold problem} for elementary polycycles: Ilyashenko and
Yakovenko~\cite{iy} proved that the number of limit cycles generated
by an elementary polycycle in a generic smooth family of planar vector
fields with $k$ parameters is finite, and Kaloshin~\cite{kaloshin}
established an effective upper-bound in terms of $k.$ We refer the
reader to~\cite{ilya} for an historical overview of the many
developments surrounding the second part of Hilbert's sixteenth
problem.

\bigskip

Another important application of Pfaffian functions is to {\bf
fewnomials}, or sparse polynomials. Recall that for a real univariate
polynomial $p(x)=\sum_{i=1}^{\, r} a_i x^{m_i}$ (where $a_i\neq 0$ for all
$i$), Descartes's rule says that the number of positive real roots of
$p(x)$ (counted with multiplicity) is bounded by the number of indices
$i$ such that $a_ia_{i+1}<0.$ As a corollary, the number of zeros of
$p(x)$ on the whole real line (this time counted without multiplicity)
must be bounded by $2r-1,$ independently of the degree of $p(x).$

\medskip

Using Pfaffian functions, one can generalize the above result to
polynomial systems in several variables, and produce an effective
bound which is exponential in $r^2$ and polynomial in the number of
variables. Similar results can be established for polynomials with a
low {\bf additive complexity.} Informally, the additive complexity of
a polynomial is the minimum number of additions necessary for a
program to evaluate this polynomial. Thus, the polynomial
$f(x)=(1+x^p+x^q)^r$
has an additive complexity of~2 for any values of $p,q$ and $r,$ and
the number of real roots of $f(x)$ can be bounded independently of
$p,q$ and $r.$\footnote{Here, the fact that we are dealing with polynomials is
irrelevant: the Pfaffian function technique applies also to Laurent
polynomials, or even positive real exponents.}

\section*{The tame topology of Pfaffian sets}
\addcontentsline{toc}{section}{Tame topology of Pfaffian sets}

The polynomial-like behaviour of Pfaffian functions means that the
sets crafted from such functions have very nice geometrical
properties, avoiding the most pathological examples of topology.
These sets have a well-defined notion of dimension (which is an
integer) and many of their topological and geometrical characteristics
are finite, 
and in the same fashion that the
B\'ezout inequality can be used to derive effective bounds for the
complexity of semi-algebraic sets, we can use Khovanskii's theorem to
give these finiteness results a quantitative aspect.

\bigskip

\begin{quote}
\noindent
{\bf Object of the dissertation. --}
\itshape
The object of this thesis is to illustrate the above remark
concretely, by giving
effective bounds for the Betti numbers of Pfaffian sets.
\upshape
\end{quote}

\bigskip

Such an endeavor involves more, however, than simply translating
well-known results about semi-algebraic sets into some result about
Pfaffian functions. This is because Pfaffian sets come in different
flavors, more so than semi-algebraic sets, each of them being relevant
by appearing naturally in some contexts, and each type needing a
different treatment. We will see shortly what those different flavors
of Pfaffian sets are.

\medskip

Fortunately, this tame topology can be still studied in a unifying
framework, the theory of {\bf o-minimal structures.}  As the reader
will soon see, Pfaffian sets offer a confusing picture where some
questions are still left unanswered. But o-minimality will offer us
powerful tools that will let us treat Pfaffian sets without having to
solve those questions.




\subsection*{O-minimality and Pfaffian functions}
\addcontentsline{toc}{subsection}{O-minimality and Pfaffian functions}

Before we can explain what is the o-minimal structure associated to
Pfaffian functions and why it is relevant to us, let us go back to
semi-algebraic sets.

\medskip

The semi-algebraic subsets of $\R^n$ are, by definition, sets
belonging to the Boolean algebra $SA_n$ generated by the sets of the
form $\{q>0\}$ for all $n$-variate real polynomials $q.$ Thus, $SA_n$
is {\bf stable under finite unions, finite intersections and
complementation.}  When considering different dimensions, we have
additional stability properties: a fairly obvious {\bf stability by
Cartesian product} ({\em i.e.} if $A\in SA_m$ and $B\in SA_n,$ we must have
$A\times B\in SA_{m+n}$), and by a result known as the
Tarski-Seidenberg theorem, {\bf stability under linear projection}: if
$\pi$ is the canonical projection $\R^{m+n} \to \R^n$ and $A\in
SA_{m+n},$ then $\pi(A)\in SA_n.$ These properties make the collection
of all semi-algebraic sets a {\bf structure.}  In practice, this means
that starting from semi-algebraic data to be studied, we should expect
to encounter only semi-algebraic sets.%
\footnote{This is to be contrasted with a situation of the following
type: $V \sub \R^n$ is a real algebraic variety. The set of
non-singular points $V^* \sub V$ is always semi-algebraic, but is not
generally algebraic.}

\medskip

Similarly, we can define {\bf semi-Pfaffian sets} in $\R^n$ as the
sets that belong to the Boolean algebra generated by the sets of the
form $\{q>0\},$ where $q$ is a Pfaffian function in $n$ variables.
These sets are not stable under linear projection, as can be seen from
the classical counter-example of Osgood~\cite{osgood}. Worse: if $X \sub \R^n$
is a {\bf sub-Pfaffian set,} {\em i.e.} the projection of a
semi-Pfaffian set, it is not known in general if the complement
$\R^n \bs X$ is again sub-Pfaffian.%
\footnote{If this was known to be always true, sub-Pfaffian sets would
form a structure. See Remark~\ref{rem:sub}.}

\bigskip

In general, a structure is called {\bf o-minimal} if all sets in the
structure have finitely many connected components. An example of such
a structure is given by the structure of semi-algebraic sets.
O-minimality is a
notion that originated in mathematical logic, more precisely in model
theory~\cite{vdd:pfaff,vdd:book,om2,om1}.  Sets in o-minimal
structures are notoriously well-behaved from the topological
viewpoint: the structure stability property plus the finiteness of
connected components is enough to ensure that o-minimal structures
admit an analogue of the cylindrical algebraic decomposition of
semi-algebraic sets.  Thus, all such structures are very similar
geometrically and topologically.  In particular, the sum of the Betti
numbers of any set that belongs to an o-minimal structure is always
finite.

\bigskip


Since all the structure operations may come naturally in our study of
sets defined by Pfaffian functions, another way to look at things is
to define the {\bf Pfaffian structure} as the smallest collection of
sets containing all semi-Pfaffian sets and stable under all the
structure operations. Since Pfaffian functions are so well-behaved, it
was natural to hope that the Pfaffian structure would turn out to be
o-minimal. This fact was finally proved by Wilkie~\cite{wilkie99},
building on ideas of Charbonnel~\cite{charb} (see
also~\cite{kamac:pf-cl,lr:vol,sp:pf-cl}). 
One of the difficulties is to construct the Pfaffian structure in a
more effective fashion. Since no theorem of the complement is
available for sub-Pfaffian sets, it was necessary to build the
structure using an additional operation, the {\bf closure at infinity}.
The good news is that we can now work without too many worries: any
set we construct in the Pfaffian structure will be well-behaved. The
problem we set out to solve makes sense for any set in the Pfaffian
structure, since any such set has a finite sum of Betti numbers.
Moreover, tameness arguments will be very helpful when manipulating
these sets. 

\bigskip

A problem remains: Wilkie's construction does not yield a manageable
notion of complexity for Pfaffian sets.  As an alternative to
Wilkie's construction, Gabrielov~\cite{ga:rc} defined {\bf limit sets}
to describe the Pfaffian structure in a way that would allow to answer
quantitative questions.  If $X$ is a {\bf 1-parameter semi-Pfaffian family},
defined for a parameter $\lambda>0,$ we let $\ck{X}$ be the {\bf Hausdorff
limit} of the sets $\ol{X_\lambda}$ when $\lambda$ goes to zero. If
$(X,Y)$ is a {\bf semi-Pfaffian couple}, {\em i.e.} a couple of such
families verifying some additional conditions, the {\bf relative
closure} of $(X,Y)$ is defined as $(X,Y)_0=\ck{X}\bs \ck{Y}.$
Gabrielov proved in~\cite{ga:rc} that any set in the Pfaffian
structure is a {\em limit set:} a finite union of such relative
closures.

\subsection*{Distinguished types of definable sets}
\addcontentsline{toc}{subsection}{Distinguished types of definable sets}

Sets belonging to a structure are usually called {\bf definable} in
the structure. Most of the time in this thesis, {\em definable} will
refer to the Pfaffian structure. From the point of view of complexity
of their description we will distinguish in our treatment the
following types of set.
\begin{itemize}
\item A {\bf Pfaffian  variety} $V$ is a set defined by a condition of the
form $q_1(x)=\cdots=q_r(x)=0.$ We will write $V=\Zcal(q_1,
\ldots,q_r).$ 
\item A {\bf semi-Pfaffian set} is given by a Boolean combination of
sign conditions on Pfaffian functions.
\item A {\bf sub-Pfaffian set} is the linear projection of a
semi-Pfaffian set.
\item The {\bf relative closure} of a semi-Pfaffian couple $(X,Y)$ is
the set $(X,Y)_0=\ck{X} \bs \ck{Y}.$
\item A {\bf Pfaffian limit set} is a finite union of relative
closures. 
\end{itemize}

Note that as we go through this list, the definitions get more and
more complicated.  For any set in one of the above type, we can
associate a {\bf format} which is a tuple of integers that measures
the combinatorial (number of Pfaffian functions) and algebraic
(degrees involved, length of the chain, number of variables)
complexity of the description of the set. All upper-bounds are
functions of the format.

\bigskip

The reader will find a much more detailed treatment of all this
material in Chapter~1, including Pfaffian functions and Khovanskii's
theorem, of course, but also a detailed account of o-minimal
structures, from their basic properties to more specific details about
their relevance in the present work.

\section*{Overview of the results}
\addcontentsline{toc}{section}{Overview of the results}

As mentioned previously, our interest is with establishing
upper-bounds for the various Betti numbers of sets of any of the above
type. 
The fundamental remark for our purpose is the following: Khovanskii's
theorem allows us to bound the sum of Betti numbers of any Pfaffian
variety $V$ in terms of the degrees of the equations defining $V.$
Indeed, it is well-known from classical Morse theory that if $V$ is
smooth and compact, the sum of its Betti numbers is bounded by the
number of critical points of a Morse function on $V.$ For a
well-chosen Morse function, those critical points will be solution of
a Pfaffian system and Khovanskii's estimate will apply. Tameness of
the Pfaffian structure lets us reduce the case of an arbitrary variety
to this favorable situation.

\medskip

The bound will be denoted by $\V(\cdots)$ in the text, %
\footnote{Where $\cdots$ is the format of the variety, see
Definition~\ref{df:format} for more details.}  and all the results in
this thesis can be expressed in terms of this bound $\V$ (and most of
them are). The ingredients we will need to reduce all of our problems
to questions about varieties range from algebraic topology techniques
to general position arguments and tame topology arguments coming from
o-minimality.

\subsection*{Topology of semi-Pfaffian sets (Chapter~2)}
\addcontentsline{toc}{subsection}{Topology of semi-Pfaffian sets}

Chapter~2 is devoted to semi-Pfaffian sets, and as such, the results
in the chapter are closely related to analogous results for
semi-algebraic sets. 

\medskip

First, the bound $\V$ for Pfaffian varieties that was
mentioned above is established (Theorem~\ref{TH:varieties}).  The
main ideas were explained above, and the method follows the classical
ideas of Oleinik, Petrovsky Thom and
Milnor~\cite{oleinik,op,thom,milnor:betti} in the algebraic case.

\bigskip

From Theorem~\ref{TH:varieties}, we then derive a bound on the Betti
numbers of a compact semi-Pfaffian set defined by a $\P$-{\bf closed formula}
\footnote{If $\P=\{p_1, \ldots, p_s\}$ is a set of Pfaffian functions,
a $\P$-closed sign condition is a sign condition obtained by
conjunctions and disjunctions of atoms, each atom being of the form
$\{p_i\geq 0\},$ $\{p_i \leq 0\}$ or $\{p_i=0\}$ for any $p_i \in \P.$
In particular, {\bf no negations} are allowed.}
(Theorem~\ref{th:Pclosed}).  Then, in Theorem~\ref{th:cells}, we
establish a bound on the {\bf number of connected sign cells}
$\Ccal(\P)$of a family of Pfaffian functions $\P$ ($\Ccal(\P)$ is a
quantity that acts both as a bound on the number of connected
components of any semi-Pfaffian set defined from $\P$ and as a bound
on the number of consistent sign assignment for $\P$). We deduce from
these results a bound on the rank of the {\bf Borel-Moore homology} of
locally-closed semi-Pfaffian sets (Theorem~\ref{th:bm}).  This last
inequality can be used to bound the ordinary Betti numbers of a
compact semi-Pfaffian set, when that set is not defined by a
$\P$-closed formula (but the bound is worse).

\bigskip

Similar estimates had already been derived in the algebraic case from
the Oleinik, Petrovsky Thom and Milnor bound, see for
instance~\cite{ba:euler,bpr:cells,bpr:book,burgisser,momopa,yao}.

\bigskip

The common thread throughout Chapter~2 is the use of induction on the
{\bf combinatorial level} of a family $\P$ of functions. This notion
was introduced in the papers of Basu, Pollack and Roy; see for
instance the book~\cite{bpr:book}. The combinatorial level of $\P$ is
defined as the largest number of functions in $\P$ that can
simultaneously vanish.  By reducing to problems in general position
(at a small combinatorial cost), we can limit ourselves to cases where
the combinatorial level is bounded by the dimension of the ambient
space. Since a combinatorial level of zero means the set is a variety,
induction on this level lets us express all our bounds on
semi-Pfaffian sets in terms of $\V.$ In particular, this allowed to
give some {\em exact} bounds for the Betti numbers, as well as
asymptotic estimates.%
\footnote{Exact in the sense that they do not depend on unknown
constants. Chapter~4 is the only other instance where giving exact
bounds will be practical. Chapters~3 and~5 contain only asymptotic
estimates.}

\bigskip

We finish Chapter~2 with a recent application of
Theorem~\ref{th:Pclosed} due to Gabrielov and Vorobjov. The result,
stated without proof, is an upper-bound on the sum of the Betti
numbers of {\em any} semi-Pfaffian set (no topological assumption or
assumption on the defining formula necessary). Such single-exponential
upper-bounds were not known previously even for semi-algebraic sets.

\subsection*{Betti numbers of sub-Pfaffian sets (Chapter~3)}
\addcontentsline{toc}{subsection}{Betti numbers of sub-Pfaffian sets}

In Chapter~3, we introduce a {\bf spectral sequence} $E_{p,q}^r$ associated
to any continuous surjection $f:X \to Y,$ and we prove that it converges to the
homology of $Y$ whenever $f$ is {\bf compact-covering} %
\footnote{$f:X \to Y$ is compact covering if and only if for any compact
$L\sub Y,$ there exists a compact $K \sub X$ with $f(K)=L.$}
(Theorem~\ref{th:ccss}).
An analogue of this spectral sequence is known as {\em cohomological
descent} in sheaf cohomology (see for instance~\cite{deligne}).
The term $E^1_{p,q}$ is isomorphic to the $q$-th homology group of the
$(p+1)$-fold fibered product $X\times_Y\cdots\times_YX,$ which is the
set of $(p+1)$-tuples of points of $X$ having the same image by $f.$
Thus, such a sequence allows to estimate the Betti numbers of $Y$ when
the Betti numbers of those fibered products are known (Theorem~\ref{th:ss}).

\bigskip

We apply this in the case where $f$ is a projection of a compact or
open semi-Pfaffian subset $X$ of a cube. The spectral sequence converges,
and since the fibered products $X\times_Y\cdots\times_YX$ are
semi-Pfaffian, we can use the results from Chapter~2 and the spectral
sequence estimate to derive a bound on the Betti numbers of the
{\bf sub-Pfaffian set} $Y=f(X)$ (Theorem~\ref{th:exist}). This technique can
also be used to estimate Betti numbers of sets defined by a universal
quantifier (Corollary~\ref{cor:univ}).

\bigskip

This technique can be applied inductively to derive an estimate on the
Betti number of any sub-Pfaffian subset of the cube that is defined
from a compact or open semi-Pfaffian set by $\nu$ {\bf quantifier
alternations}.  We obtain this way a recurrence relation
(Theorem~\ref{th:Erec}) that can be solved to establish a bound in the
Pfaffian case (Corollary~\ref{cor:EPfaff}) and in the algebraic case
(Corollary~\ref{cor:Ealg}). These bounds tend to improve on previously
known bounds when $\nu$ is small (but are worse when $\nu$ goes to
infinity).

\subsection*{Topological bounds for limit sets (Chapter~4 and~5)}
\addcontentsline{toc}{subsection}{Topological bounds for limit sets}

Chapters~4 and~5 are devoted to the study of the topology of limit
sets.

\bigskip

In Chapter~4, we establish a single-exponential bound on the number of
{\bf connected components} of the relative closure $(X,Y)_0$ of a Pfaffian
couple. Since in general limit sets are finite unions of such relative
closures, this gives explicit finiteness results for all limit sets.

\bigskip

This estimate is obtained by bounding the number of local extrema of a
generic distance function $\Phi$ restricted to a fiber $X_\lambda
\times Y_\lambda$ for $\lambda \ll 1.$ In the smooth case, the problem
reduces to studying critical points of a similar distance on an open
subset of $X_\lambda \times (Y_\lambda)^p,$ for $1\leq p \leq \dim(X)+1.$ 
Finally, we are reduced to estimating the number of connected
component of such a set, and using the results of Chapter~2, we obtain
the bound given in Theorem~\ref{th:smooth_bound}. The singular case
reduces to the smooth case by an appropriate deformation
(Theorem~\ref{th:sing_bound}). 

\bigskip

At last, Chapter~5 addresses the problem of the {\bf higher Betti numbers}
of relative closures. We start with the case of the closure of a Pfaffian
family $X_\lambda,$ where no $Y$ is subtracted. In that case, the
closure $X_0$ is simply the {\bf Hausdorff limit} of the compact family
$X_\lambda.$ We establish a general result about such Hausdorff limits 
in general o-minimal structures (Theorem~\ref{th:1param}): the Betti
numbers of the limits can be estimated in terms of Betti numbers of
simple deformations of the diagonal in Cartesian powers of the
fibers. Using the bounds from Chapter~2, this allows to estimate the
Betti numbers of $X_0$ in the relative closure case (Theorem~\ref{th:rcbnd}). 

\bigskip

In the presence of a {\bf non-empty} $Y,$ we establish that the Betti
numbers of $(X,Y)_0$ can be bounded from above in terms of the formats
of the fibers $X_\lambda$ and $Y_\lambda$ (Theorem~\ref{th:rc-betti}).
Thus, we confirm that the format of the fibers (rather than the format
of the whole family) is the right measure of the complexity relative
closures \footnote{at least from the point of view of Betti
numbers.}. The dependence in the parameter does not affect the
topology of the limits.

\bigskip

The proof of these results relies heavily on the spectral sequence
introduced in Chapter~3. We show that if $X$ is a definable
1-parameter family of compact sets, we can construct (in a
non-effective way) a {\bf family of continuous surjections} $f^\lambda$ from
generic fibers $X_\lambda$ to the Hausdorff limit $X_0.$ The core of
the proof is then to show that the tame behaviour of the family
$f^\lambda$ allows us to prove that for $\lambda$ small enough, the
Betti numbers of the fibered products coming from the spectral
sequence coincide with the Betti numbers of certain deformations of
the diagonal.

\section*{Some general remarks about the results}
\addcontentsline{toc}{section}{Some general remarks about the results}

It is widely believed that the bound from Khovanskii's theorem is not
sharp when the number $\ell$ of functions in the Pfaffian chain is
large. Thus, any tightening of that bound would also improve the bound
$\V$ for Pfaffian varieties, and thus all of the results presented
here.  Khovanskii's method relies essentially on bounding the number
of real isolated solutions of a system of $n$ Pfaffian functions by
the number of real solutions of a system of $n+\ell$ polynomial
equations in $n+\ell$ variables derived inductively from the original
Pfaffian system. The number of solutions of the polynomial system is
then estimated by the B\'ezout bound, but there is no reason to
believe that this bound would be sharp for those specific systems.

\bigskip

The results presented here have been published or are about to be for
the most part: sections~2.1 and~2.2 of Chapter~2 appeared
in~\cite{z99}, Chapter~3 appeared in~\cite{gvz} and Chapter~4
in~\cite{gz:cc}. A paper containing the results of Chapter~5 (with a
stronger emphasis on the aspect of Hausdorff limits in general
o-minimal structures) has been submitted~\cite{z03}, and explicit
estimates for a general relative closure are work in
progress~\cite{z03:rc}.

\bigskip

Lastly, it is worth mentioning that the recent result of Gabrielov and
Vorobjov~\cite{gv:qf} giving good ({\em i.e.} close to optimal)
estimates for sums of Betti numbers of sets defined by any quantifier-free
formulas 
came out too late to be completely incorporated into the text,
although it was used when it offered new possibilities or removed
substantial difficulties in the original version.


\emptyevenpage

%% file: aa1.tex
This chapter presents all the necessary background material about
Pfaffian functions and o-minimal structures. The material is organized
as follows.
The first section introduces {\em Pfaffian functions} along with the
bounds of Khovanskii about the number of solutions of a Pfaffian
system. Section~2 deals with {\em semi and sub-Pfaffian sets,} and
their formats; section~3 is about {\em o-minimal structures} on the
real field and their basic geometric properties.  At last, Pfaffian
{\em limit sets} are introduced in section~4. This section finishes
with some corollaries of the o-minimality of the structure of Pfaffian
functions that will be widely used in the other chapters.

\section{Pfaffian functions}

In this section, we define Pfaffian functions following
Khovanskii; we define the notion of complexity of Pfaffian
functions and state the fundamental result in the theory: any system
of Pfaffian equations has a finite number of {\em isolated}%
\footnote{Real soltions isolated over $\C,$ that is.}
solutions, that
can be effectively estimated from above by an expression involving
only the discrete parameters of the Pfaffian system (degrees, number
of variables, and chain length). These parameters are often referred
to as the {\em format} or {\em Pfaffian complexity} of the functions.

\subsection{Definition and examples}

\begin{definition}[Pfaffian chain]\label{df:chain}
Let $\f=(f_1, \ldots, f_\ell)$ be a sequence of real analytic
functions defined on a domain $\U \sub \R^n.$ We say that they
constitute a {\em Pfaffian chain} if there exists polynomials
$P_{i,j},$ each in $n+i$ indeterminates, such that the following equations
\begin{equation}\label{eq:chain}
\DP{f_i}{x_j}(x)=P_{i,j}(x, f_1(x), \ldots, f_i(x)), \quad 1 \leq i
\leq \ell, \; 1\leq j\leq n,
\end{equation}
hold for all $x\in \U.$ 
\end{definition}

This definition is sufficient when considering functions that are all
simultaneously defined.  However, in all generality, one
should use the following definition.

\begin{definition}[Pfaffian chain 2]
\label{def:1}
A sequence $\f=(f_1, \ldots, f_\ell)$ of analytic functions in $\U$ is
called a {\em Pfaffian chain} if it satisfies on $\U$ a differential
system of the form:
\begin{equation}
\label{eq:chain2}
df_i=\sum_{j=1}^n P_{i,j}(x, f_1(x), \ldots, f_i(x)) dx_j,
\end{equation}
where each $P_{i,j}$ is a polynomial in $n+i$ indeterminates, and the
following holds.
\begin{itemize}
\item[{\rm (P1)}] The graph $\Gamma_i=\{t=f_i(x)\}$ of $f_i$ is
contained in a domain $\Omega_i$ defined by polynomial inequalities in
$(x,f_1(x), \ldots, f_{i-1}(x),t),$ and such that $\fr \Gamma_i \sub
\fr \Omega_i.$
\item[{\rm (P2)}] $\Gamma_i$ is a {\em separating submanifold} in
$\Omega_i,$ {\em i.e.}  $\Om_i \setminus \Gamma_i$ is a disjoint union
of the two sets $\Om_i^+=\{f_i>0\}$ and $\Om_i^-=\{f_i<0\}.$
(\/See~\cite[p.~38]{kh:few}. This is also called the {\em Rolle leaf}
condition in the terminology of~\cite{lr:hs,lr:vol}.\/)
\end{itemize}
\end{definition}

\begin{definition}[Pfaffian function]\label{df:pf}
Let $\f=(f_1, \ldots, f_\ell)$ be a fixed Pfaffian chain, and $q(x)$ be
an analytic function on the domain of that chain. We say that $q$ is
{\em a Pfaffian function in the chain} $\f$ %
if there exists a polynomial $Q$ such that $q=Q(x,\f),$ {\em i.e.}
\begin{equation}\label{eq:pf}
q(x)=Q(x, f_1(x), \ldots, f_\ell(x)) \quad \forall x\in\U.
\end{equation}
\end{definition}

\begin{definition}[Format]\label{df:format}
Let $\f=(f_1, \ldots, f_\ell)$ be a Pfaffian chain. We call $\ell$ the
{\em length} (also called {\em depth} or {\em order}) of $\f.$ We say $\f$
is of {\em degree} $\alpha$ if all the polynomials $P_{i,j}$ appearing
in~\eqref{eq:chain} are of degree at most $\alpha.$ If $Q$ is a
polynomial of degree $\beta$ in $n+\ell$ variables and $q=Q(x,\f),$ we
say that $\beta$ is the {\em degree of $q$ in $\f,$} and we will write
$\beta=\deg_{\f}(q).$ 
\end{definition}

\bigskip

\noindent
{\bf Examples of Pfaffian functions}
\begin{enumerate}
\item The polynomials are the Pfaffian functions such that $\ell=0.$

\item The exponential function $f_1(x)=e^x$ is Pfaffian, with $\ell=1$
  and $\alpha=1,$ because of the equation $f_1'=f_1.$ More generally,
  we can define the iterated exponential functions by the induction
  $f_r(x)=\exp(f_{r-1}(x))$ for all $x.$ Then, $(f_1, \ldots, f_r)$ is
  a Pfaffian chain of length $r$ and degree $r$ for all $r,$ since
  $f'_r=f'_{r-1}f_r=f_1\cdots f_r$ (by induction).

\item Let $\U=\R\bs\{0\},$ and let $f(x)=x^{-1}$ and  $g(x)=\ln |x|.$
Then, $(f,g)$ is a Pfaffian chain of degree $\alpha=2$ on $\U,$ since
we have $f'=-f^2$ and $g'=f.$

\item Let $f(x)=(x^2+1)^{-1}$ and $g(x)=\arctan x.$ Then, $(f,g)$ is
a Pfaffian chain of degree $\alpha=3$ on $\R$ since we have $f'=-2xf$
and $g'=f.$

\item Let $f(x)=\tan x$ and $g(x)=\cos^2 x.$ We have $f'=1+f^2$ and
$g'=-fg,$ so $(f,g)$ is a Pfaffian chain of degree $\alpha=2$ on the 
domain $\{x \not\equiv \frac{\pi}{2} \, [\mathrm{mod}\,  \pi] \}.$ The
function $h(x)=\cos(2x)$ is Pfaffian in this chain, since  we have
$h(x)=2\,g(x)-1.$

\item Let $m\geq 2$ be an integer, and $f$ and $g$ be as above. Then,
the previous example shows that $(f(x/2m), g(x/2m))$ is a Pfaffian
chain on the domain $\{x \not\equiv m\pi\, [\mathrm{mod}\, 2m\pi] \}.$
Then $\cos x$ is a Pfaffian function of degree $m$ in that chain,
since $\cos x$ is a polynomial of degree $m$ in
$\cos(x/m)=2\,g(x/2m)-1.$ 

\item $f(x)=\cos x$ is not Pfaffian on the whole real line, since
$f(x)=0$ has infinitely many isolated solutions (see
Theorem~\ref{th:kh}).

\end{enumerate}

More generally, if we consider the following functions (in any finite
number of variables): polynomials, exponentials, trigonometric
functions and their composition inverses wherever applicable. Then,
the {\em real elementary functions} is the class obtained from these
by taking the closure under arithmetical operations and composition.
If $f$ is in this class and the functions $\sin$ and $\cos$ appear in
$f$ only through their restriction to {\em bounded} intervals, the $f$
is Pfaffian on its domain of definition (See~\cite[\oldS1]{kh:few}). 

\bigskip

Still, one of the most important applications of Pfaffian functions is
to polynomials themselves, and more specifically to the so-called {\em
fewnomials.}

\begin{definition}[Fewnomials]\label{df:fewnomials}
Fix $\K=\{m_1, \ldots, m_r\} \in \N^n$ a set of exponents. 
The polynomial $q$ is a {\em $\K$-fewnomial} if it is of the form:
\begin{equation*}
q(x)=Q(x^{m_1}, \ldots, x^{m_r}),
\end{equation*}
where $Q$ is a polynomial in $r$ variables. If $\beta=\deg(Q),$ we say
that $q$ has {\em pseudo-degree} $\beta$ in $\K.$
\end{definition}

Let $\ell=n+r,$ and $\f=(f_1, \ldots, f_\ell)$ be the chain
defined by:
\begin{equation}
\label{eq:few_chain}
f_i(x)=
\begin{cases}
&x_i^{-1} \hbox{ if $1 \leq i \leq n,$}\\
& x^{m_{i-n}} \hbox{ if $i>n.$}\\
\end{cases}
\end{equation}
It is easy to see that $\f$ is a Pfaffian chain of length $\ell$ and
degree $\alpha=2$ in the domain $\U=\{x_1\cdots x_n \neq 0\}$, since
we have:
\begin{equation*}
\DP{f_i}{x_j}=
\begin{cases}
-f_i^2 \hbox{ if $i=j \leq n,$} \\
f_j f_i \hbox{ if $i >n.$}\\
\end{cases}
\end{equation*}
Then, a $\K$-fewnomial $q$ can be seen as a Pfaffian function in $\f,$
with $\deg_{\f} q$ equal to the pseudo-degree of $q,$ but its format
is completely independent of the usual {\em degree} of $q.$ This fact
will enable us to generalize the well-known consequence of Descartes's
rule: a univariate polynomial with $m$ non-zero monomials has at most
$m-1$ positive roots.

\begin{rem}\label{rem:log}
This is not necessarily the best way to see a $\K$-fewnomial as a
Pfaffian function. On the first quadrant $(\R_+)^n,$ one can make the
change of variables $t_i=\log x_i,$ and questions about
$\K$-fewnomials can thus be reduced to questions about Pfaffian
functions in the chain $(e^{m_1\cdot t}, \ldots , e^{m_r \cdot t}),$
which is of length only $r$ compared to the chain~\eqref{eq:few_chain}
which has length $n+r.$
\end{rem}

These considerations about fewnomials can be generalized in many ways:
we do not need the exponents to be integers, and we can consider
functions in these chains of degree larger than one. Though the number
of monomials of such a function may depend on the values of $m_1,
\ldots, m_r$ they are still well-behaved. More generally, one can
consider the additive complexity of polynomials.

\begin{definition}[Additive complexity~\cite{br,kh:few}]\label{df:addcplx}
Let $m\in \N^n$ and $c \in \R\bs\{0\}.$ Then, the polynomial $c+x^m$ is said
to have {\em additive complexity } 1. If $q$ is a polynomial, we say
its {\em additive complexity is bounded} by $k+1$ if $q(x)=c+x^{m_0}
p_1(x)^{m_1}\cdots p_k(x)^{m_k},$ where $m_0 \in \N^n,$ and for all $1
\leq i \leq k,$ $m_i \in \N$ and $p_i$ is a polynomial of additive
complexity bounded by $i.$
\end{definition}

In particular, if $p$ has an additive complexity bounded by $k,$ it
means that it can be evaluated using at most $k$ additions. Since a
function of the form $p(x)^m$ is Pfaffian with a complexity
independent of $m$ on the domain $\{p(x) \neq 0\},$ so this notion can
be approached from the point of view of Pfaffian functions. Such an
approach yields explicit bounds on the number of roots of such
polynomials (see Theorem~\ref{th:addcplx}).

\begin{proposition}\label{prop:algebra}
Let $\f=(f_1, \ldots, f_\ell)$ be a Pfaffian chain on a domain $\U
\sub \R^n.$ Then, the algebra generated by $\f$ is stable under
differentiation. Moreover, the degree in $\f$ of the sum, product, and
partial derivatives of functions from this algebra can be estimated in
terms of the format of the original functions.
\end{proposition}

\begin{proof}
Let $g=G(x,\f)$ and $h=H(x,\f)$ be two functions from the algebra
generated by $\f,$ with $\deg(G)=\beta_1$ and $\deg(H)=\beta_2.$ We
have
\begin{equation*}
(g+h)(x)=G(x,f_1(x), \ldots, f_\ell(x))+H(x,f_1(x), \ldots,
f_\ell(x)),
\end{equation*}
so $g+h$ is in the algebra generated by $\f,$ and we have
$\deg_{\f}(g+h)\leq\max(\beta_1,\beta_2).$

\medskip

Similarly, we have
\begin{equation*}
(gh)(x)=G(x,f_1(x), \ldots, f_\ell(x))\,H(x,f_1(x), \ldots,
f_\ell(x)),
\end{equation*}
so $gh$ is in the algebra generated by $\f$ 
with $\deg_{\f}(gh)=\beta_1+\beta_2.$

\medskip

At last, we have by the chain rule,
\begin{equation*}
\DP{g}{x_j}(x)=\DP{G}{X_j}(x, \f(x))
+\sum_{k=1}^\ell \DP{G}{Y_k}(x, \f(x)) \, P_{k,j}(x, \f(x)).
\end{equation*}
The stability under derivation of the algebra generated by $\f$
follows. If the degree of the chain $\f$ is $\alpha,$ the degree of
any first-order derivative of $g$ is bounded by $\alpha+\beta_1-1.$
\end{proof}

\begin{rem}
If $\f_1$ and $\f_2$ are two Pfaffian chains defined on the same
domain $\U,$ of length respectively $\ell_1$ and $\ell_2$ and degree
$\alpha_1$ and $\alpha_2,$ the concatenation of $\f_1$ and $\f_2$
gives a new Pfaffian chain $\f$ of length at most $\ell_1+\ell_2$ and
degree $\max(\alpha_1,\alpha_2).$ Thus, we can always work in the
algebra generated by a fixed chain $\f.$
\end{rem}

\subsection{Khovanskii's theorem}

The fundamental result about Pfaffian functions is the
following theorem.

\begin{theorem}[Khovanskii]\label{th:kh}
Let $\f$ be a Pfaffian chain of length $\ell$ and degree $\alpha,$
with domain $\R^n.$ Let $Q_1, \ldots ,Q_n$ be polynomials in $n+\ell$
variables of degrees respectively $\beta_1, \ldots, \beta_n,$ and let
for all $1 \leq i \leq n,$ $q_i(x)=Q_i(x,\f).$
Then the number of %
solutions of the system
\begin{equation}\label{eq:thkh}
q_1(x)=\cdots=q_n(x)=0,
\end{equation}
that are isolated in $\C^n$ is bounded from above by
\begin{equation}\label{eq:kh}
2^{\ell(\ell-1)/2} \, \beta_1\cdots\beta_n\; (\beta_1+\cdots+\beta_n
-n+\min(n,\ell)\alpha+1)^\ell.
\end{equation}
\end{theorem}

The above bound can be found in~\cite[\oldS3.12, Corollary~5]{kh:few}.
It also holds when the domain of the functions is the quadrant
$(\R_+)^n.$ Over $\C^n,$ the result is of course not true,
since $e^x$ is a Pfaffian function. The complex analogue of the above
result is a bound on the multiplicity of the root of a system of
complex Pfaffian functions~\cite{ga:multi-pf} (see also~\cite{gakh}).

\medskip

Roughly, the method of proof is the following: one has to replace
inductively the functions $f_i(x)$ by variables $y_i,$ starting from
$f_\ell(x).$ At each step, a Rolle-type argument allows to produce an
extra polynomial $Q_{n+i}$ so that the system
\begin{equation*}
Q_j(x, f_1(x), \ldots, f_{i-1}(x), y_i, \ldots, y_\ell)=0, 
\quad 1 \leq j \leq n+i; 
\end{equation*}
has at least as many isolated solutions as the original system.
Thus, after replacing $f_1(x),$ one obtains a system of $n+\ell$
polynomial equations in $n+\ell$ unknowns. The degrees of the
polynomials $Q_{n+1}, \ldots, Q_{n+\ell}$ can be effectively
estimated, and by B\'ezout's theorem, the number of isolated solutions
of the final system can be bounded.

\begin{rem}
In~\cite{kh:few}, Theorem~\ref{th:kh} is formulated as a bound on the
number of {\em non-de\-ge\-ne\-rate} roots of the
system~\eqref{eq:thkh}.  If $q$ is the map $(q_1, \ldots, q_n),$ the
number of non-de\-ge\-ne\-rate roots of the system is simply the
number of points $x$ in the preimage $q^{-1}(0)$ for which the rank of
$dq(x)$ is maximal. The two formulations are clearly equivalent.
\end{rem}

Considering systems defined by sparse polynomials involving $r$
non-zero monomials in the positive quadrant, one can use the change of
variables $t_i=\log x_i,$ -- as explained in Remark~\ref{rem:log}, -- to
reduce the problem to a problem about systems involving $r$
exponential functions. One can then bound the number of non-degenerate
solutions independently of the degrees of the polynomials, to obtain
the following estimate.

\begin{corollary}[Fewnomial systems]\label{cor:kh_few}
Let $q_1, \ldots, q_n$ be polynomials in $n$ variables such that $r$
monomials appear with a non-zero coefficient in at least one of these
polynomials. Then, the number of non-degenerate solutions of the system 
\begin{equation*}
q_1(x)=\cdots=q_n(x)=0,
\end{equation*}
in the quadrant $(\R_+)^n$ is bounded by
\begin{equation}\label{eq:kh_few}
2^{r(r-1)/2} \, (n+1)^r.
\end{equation}
\end{corollary}

For systems defined by polynomials of additive complexity bounded by
$k,$ (see Definition~\ref{df:addcplx}), there is a detailed proof
in~\cite[Chapter~4]{br} that the number of non-degenerate roots admits
a computable upper-bound in terms of $k.$ In the case $n=1,$ the
following explicit bound is given~\cite[Theorem~4.2.4]{br}.

\begin{theorem}[Bounded additive complexity]\label{th:addcplx}
Let $p(x)$ be a univariate polynomial of additive complexity bounded
by $k.$ The number of real roots of $p$ is at most
\begin{equation*}
(k+2)^{2k+1} \, 2^{2k^2+2k+1};
\end{equation*}
which is less than $5^{k^2}$ for $k$ large enough.
\end{theorem}

\subsection{Domains of bounded complexity}

We will now define a class of domains $\U$ over which Khovanskii's
result can be easily generalized. Note that in order to have the nice
topological and geometrical properties we hope for, one cannot
generalize these results to domains that would be too pathological.

\begin{definition}[Domain of bounded complexity]\label{df:U}
We say that $\U$ is a domain of bounded complexity $\gamma$ for the
Pfaffian chain $\f=(f_1, \ldots, f_\ell)$ if there exists a
function $g$ of degree $\gamma$ in the chain $\f$ such that the sets
$\{g \geq \e\}$ form an exhausting family of compact subsets of $\U$
for $\e\ll 1.$ We call $g$ an {\em exhausting function } for $\U.$
\end{definition}

\begin{example}
Let $\f=(f_1,\ldots,f_\ell)$ be a Pfaffian chain defined on a {\bf
bounded} domain $\U$ of the form
\begin{equation}\label{eq:U}
\U=\{x \in \R^n \mid g_1(x)>0, \ldots, g_r(x)>0\},
\end{equation}
where $(g_1, \ldots, g_r)$ are Pfaffian functions in the chain $\f.$
Then, $\U$ is a domain of bounded complexity, since $g=g_1\cdots g_r$
is clearly an exhausting function for $\U.$
\end{example}

Note that the assumption of boundedness of $\U$ can be dropped: let
\begin{equation}\label{eq:rho}
\rho(x)= \frac{1}{1+|x|^2}.
\end{equation}
The function $\rho$ is a Pfaffian function defined on $\R^n,$ with a
degree $\alpha=3,$ since we have
\begin{equation*}
d\rho(x)=-2\rho^2(x)\,(x_1dx_1+\cdots+x_ndx_n).
\end{equation*}
Moreover, $\rho(x)>0$ on $\R^n$ and the sets $\{\rho \geq \e\}$ are
compact for $0<\e<1.$ So even an {\bf unbounded} domain $\U$ of the
form~\eqref{eq:U} is a domain of bounded complexity for any Pfaffian
chain of the form $(\rho, f_1, \ldots, f_\ell),$ with exhausting
function $g=g_1\cdots g_r + \rho.$

\bigskip

Over a domain of bounded complexity, Khovanskii's estimates becomes
the following.

\begin{theorem}[Khovanskii's theorem for a domain of bounded
complexity]\label{th:khbc} Let $\f$ be a Pfaffian length of degree
$\alpha$ and length $\ell$ defined on a domain $\U$ of bounded
complexity $\gamma$ for $\f.$ Let $Q_1, \ldots, Q_n$ be polynomials in
$n+\ell$ variables of degree respectively $\beta_1,
\ldots, \beta_n$ and let $q_i=Q(x,\f)$ for all $i.$
Then, the number of solutions in $\U$ of the system
\begin{equation}\label{eq:systS}
q_1(x)=\cdots=q_n(x)=0;
\end{equation}
which are isolated in $\C^n$ is bounded by 
\begin{equation}\label{eq:kh_eff}
2^{\ell(\ell-1)/2} \beta_1 \cdots \beta_n
\frac{\gamma}{2}
[\beta_1+\cdots+\beta_n+\gamma-n+\min(n+1,\ell)\alpha]^\ell
\end{equation}
\end{theorem}

\begin{proof}
Introduce a new variable $t$ and consider the system given by
\begin{equation}\label{eq:systSS}
q_1(x)=0, \ldots, q_n(x)=0, \, g(x)-t^2=\e;
\end{equation}
where $\e$ is a fixed positive real number. Then, for any values of
$\e,$ every isolated solution of the system~\eqref{eq:systS} that is
contained in the domain $\Omega_\e=\U \cap \{g(x)>\e\}$ gives rise to
exactly two isolated solutions for~\eqref{eq:systSS}. So it is enough to bound
the number of isolated solutions of~\eqref{eq:systSS} for a value of
$\e$ such that all the isolated solutions of~\eqref{eq:systS} are contained in
$\Omega_\e.$ The choice of the parameter $\e$ does not affect the
complexity of the new system, and the bound~\eqref{eq:kh_eff} can then
be established following the results appearing in~\cite{kh:few}.
\end{proof}

\bigskip

\section{Semi and sub-Pfaffian sets}

Semi and sub-Pfaffian sets occur naturally in the study of Pfaffian
functions: semi-Pfaffian sets are sets that can be defined by a
quantifier-free sign condition on Pfaffian functions, and sub-Pfaffian
sets are linear projections of semi-Pfaffian sets, or equivalently,
defined by existential sign conditions on Pfaffian functions.

\begin{example}\label{ex:semi-sub}
Let $q$ be a Pfaffian function defined on a domain $\U\sub\R^n.$ Then, the
set of critical points of $q$ is semi-Pfaffian and the set of its
critical values is sub-Pfaffian.
\end{example}

\begin{proof}
This is straightforward. The critical locus of $q$ is defined by 
\begin{equation*}
X=\left\{x \in \U \mid \DP{q}{x_1}(x)=\cdots=\DP{q}{x_n}(x)=0\right\};
\end{equation*}
and the set of critical values is 
\begin{equation*}
Y=\{y \in \R \mid \exists\, x\in X, \, q(x)=y\}.
\end{equation*}
Since the partial derivatives of $q$ are again Pfaffian functions, it
is clear that $X$ is semi-Pfaffian and $Y$ is sub-Pfaffian.
\end{proof}

\bigskip

Semi-Pfaffian and sub-Pfaffian sets have a lot of finiteness
properties. The present section contains mainly definitions and
relevant examples, and we refer the reader to the
bibliography~\cite{gv:strat, ga:sub, ga:fr-and-cl,ga:multi-pf,
gv:cyldec,pv:cyldec} for more details. A comprehensive
survey~\cite{gv:survey} will be available soon.

\bigskip

From now on, $\f=(f_1, \ldots, f_\ell)$ will be a fixed Pfaffian chain
of degree $\alpha$ defined on a domain of bounded complexity $\U \sub
\R^n,$ and we will consider only functions fro, the algebra  generated by $\f.$

\subsection{Semi-Pfaffian sets}

As mentioned in the beginning of this section, semi-Pfaffian sets are
given by quantifier-free sign conditions on Pfaffian functions. We
start by recalling the definition of quantifier-free formulas, and  we
define a notion of {\em format} for such formulas. This format will be
all the data we need to establish bounds on the topological complexity
of semi-Pfaffian sets.

\begin{definition}[QF formula]\label{df:qf}
Let $\P=\{p_1, \ldots, p_s\}$ be a set of Pfaffian functions. A {\em
quantifier-free} (QF) formula with atoms in $\P$ is constructed as follows:
\begin{enumerate}
\item An atom is of the form $p_i \star 0,$ where $1 \leq i \leq s$
  and $\star \in \{=,\leq,\geq\}.$ It is a QF formula;
\item If $\Phi$ is a QF formula, its {\em negation} $\neg \Phi$ is a
  formula;
\item If $\Phi$ and $\Psi$ are QF formulas, then their {\em
conjunction} $\Phi \wedge \Psi$
  and their {\em disjunction} $\Phi \vee \Psi$ are QF formulas.
\end{enumerate}
\end{definition}

\begin{definition}[Format of a formula]\label{df:qf_format}
Let $\Phi$ be a QF formula as above.  If 
the number of variables is $n,$ the length of $\f$ is $\ell,$
the degrees of the polynomials $P_{i,j}$ in~\eqref{eq:chain} 
is bounded by $\alpha,$ 
$s=|\P|$ and $\beta$ is the
maximum of the degrees in the chain of the functions in $\P,$ we call
$(n,\ell, \alpha, \beta, s)$ the {\em format} of $\Phi.$
\end{definition}

\begin{definition}[$\P$-closed formulas]\label{df:Pclosed}
We will say that the formula $\Phi$ is $\P$-closed if it was derived
without negations, {\em i.e.} using rules 1 and 3 only.
\end{definition}

\begin{definition}[Semi-Pfaffian set]
A set $X \sub \U$ is called {\em semi-Pfaffian} if there exists a
finite set $\P$ of Pfaffian functions and a QF formula $\Phi$ with
atoms in $\P$ such that
\begin{equation*}
X=\{ x \in \U \mid \Phi(x)\}.
\end{equation*}
We call {\em format} of $X$ the format of the defining formula $\Phi.$
\end{definition}

The notion of format is important to establish the kind of
quantitative bounds we want on the topology of semi-Pfaffian sets.
But the above definition can be improved.

\medskip

Indeed, taking example on the algebraic case, one expects equalities
and inequalities to affect very differently bounds on the topology of
the sets they define. If $V=\{x \in \R^n \mid
p_1(x)=\cdots=p_s(x)=0\}$ is a real algebraic variety defined by
polynomials of degree at most $d,$ we know by a result of
Oleinik-Petrovsky Thom and Milnor~\cite{op,oleinik,milnor:betti,thom}
that the sum of its Betti numbers is bounded by $d(2d-1)^{n-1},$ so
does not depend on $s.$

\medskip

On the other hand, when dealing with inequalities, the number $s$ of
functions does matter, as the following example shows: take
$p_i(x)=(x-i)^2$ and let $S=\{x \in \R
\mid p_1(x)>0, \ldots, p_s(x)>0\}.$ Then $S=\R\bs\{1,2,\ldots,s\},$ so
it has $s+1$ connected components.

\bigskip

To make full use of that difference between equalities and
inequalities in our formulas, we will introduce the following
definitions. 

\begin{definition}[Variety]
The semi-Pfaffian set $V \sub \U$ is called a {\em variety} if it is
defined using only equations. We will use the notation
\begin{equation*}
\Zcal(q_1, \ldots, q_r)=\{x \in \U \mid q_1(x)=\cdots=q_r(x)=0\}.
\end{equation*}
\end{definition}

\begin{definition}[Semi-Pfaffian subsets of a variety]
If $V=\Zcal(q_1, \ldots, q_r)$ is a Pfaffian variety and $\Phi$ a QF
formula, one can consider the semi-Pfaffian set $X=\{x \in V \mid
\Phi(x)\}.$ Then, the {\em format} of $X$ is defined as $(n,\ell,
\alpha, \max(\beta_1, \beta_2), s)$ where $\beta_1$ is a bound on the
degrees of $q_1, \ldots, q_r$ and $(n,\ell, \alpha,\beta_2, s)$ is
the format of $\Phi.$
\end{definition}

\begin{rem}
Such a cumbersome definition and notion of format may seem strange,
but we will see in Chapter~2 that this will allow us to
establish more precise bounds on the topology of $X,$ for which $r$ is
irrelevant and the parameter $d=\dim(V)$ plays a part.
\end{rem}

A more usual definition for semi-Pfaffian sets is to define them as
finite unions of {\em basic} semi-Pfaffian sets, where a basic set is
of the form
\begin{equation}
\label{eq:basic_sp}
B=\{x \in \U \mid \phi_1(x)=\cdots=\phi_I(x)=0, \,\psi_1(x)>0, \ldots,
\, \psi_J(x)>0 \};
\end{equation}
for some Pfaffian functions $\phi_1, \ldots, \phi_I $ and $\psi_1,
\ldots, \psi_J.$ (Writing as semi-Pfaffian set as a union of basic
ones is just putting the defining formula $\Phi$ in disjunctive normal form,
so the two definitions are clearly equivalent.)

\begin{rem}\label{rem:basic-format}
Semi-Pfaffian sets presented as union of basic sets occur frequently
in the literature. The reader should be aware that the definition of
their {\em format} in that case is different. The format of a basic
set of the form~\eqref{eq:basic_sp} is then defined as
$(I,J,n,\ell,\alpha,\beta),$ and if $X$ is the union of $N$ basic sets
$B_1, \ldots, B_N,$ of respective formats $(I_i,J_i,n,\ell,\alpha,
\beta),$ the format of $X$ is $(N, I,J,n,\ell,\alpha, \beta),$ where
$I=\max\{I_1, \ldots, I_N\}$ and $J=\max\{J_1, \ldots, J_N\}.$ The two
notions of formats are comparable.
\end{rem}

\begin{definition}[Effectively non-singular set]\label{def:eff-ns}
If $X$ is a basic semi-Pfaffian set, we'll say that $X$ is {\em
effectively} non-singular if the functions $\phi_1, \ldots, \phi_I$
appearing in~\eqref{eq:basic_sp} verify
\begin{equation*}
\forall x \in X, \quad d\phi_1(x) \wedge \cdots \wedge d\phi_I(x) \neq 0.
\end{equation*}
If $X$ is effectively non singular, it is a smooth submanifold of
$\R^n$ of dimension  $n-I.$
\end{definition}

Basic sets appear rather naturally because they are easier to handle
algorithmically. In particular, effectively non-singular basic sets is
what is used in~\cite{gv:strat} to produce a weak stratification
algorithm (using an oracle) for semi-Pfaffian sets.

\begin{definition}[Restricted set]\label{df:restricted}
We say that a semi-Pfaffian set $X$ is {\em restricted} if it is
relatively compact in $\U.$
\end{definition}

\medskip

Let us introduce now the notations we will use for the topological
invariants we want to bound.

\begin{notation}\label{not:betti}
Throughout this thesis, if $X$ is a topological space, $H_i(X)$ will
denote its $i$-th homology group with integer coefficients, $b_i(X)$
will be the $i$-th Betti number of $X,$ which is the rank of $H_i(X),$ 
and $b(X)$ will denote the sum $\sum_i b_i(X).$
\end{notation}

\bigskip

If $\P=\{p_1, \ldots, p_s\}$ a family of Pfaffian functions,
we denote by $\Sig$ be the set of strict sign conditions on $\P.$
If $\s \in \Sig,$ we have 
\begin{equation}\label{eq:Sig}
\s(x)=p_1(x)\s_1 0 \wedge \cdots \wedge p_s(x) \s_s 0;  \quad \s_i \in 
\{<,>,=\} \hbox{ for } 1 \leq i \leq s. 
\end{equation}

Then, for any fixed Pfaffian variety $V,$ we can consider the following.

\begin{definition}[Connected sign cells]\label{df:cells}
A {\em cell} of the family $\P$ on the variety $V$ is a connected
component of the basic semi-Pfaffian set $S(V;\s)=\{x \in V \mid
\s(x)\}$ for some $\s \in \Sig.$
\end{definition}

Then, we define the {\em number of connected sign cells} of $\P$ over
$V$ simply as the sum
\begin{equation}\label{eq:nbcells}
\Ccal(V;\P)=\sum_{\s \in \Sig} b_0(S(V;\s)).
\end{equation}

\begin{rem}\label{rem:cc}
In particular, for any semi-Pfaffian set $X=\{x \in V \mid \Phi(x)\},$
if $\Phi$ as atoms in $\P,$ the number of connected components of $X$
is bounded by $\Ccal(V;\P).$ 
\end{rem}

\begin{proof}
Indeed, we can assume without loss of generality that $\Phi$ is in
disjunctive normal form, in which case it is of the form $\s_1 \vee
\cdots \vee \s_N$ for some $\s_1, \ldots,
\s_N \in \Sig.$ Then, $X=S(V;\s_1) \cup \cdots \cup S(V;\s_N),$ so we have
\begin{equation*}
b_0(X) \leq b_0(S(V;\s_1))+\cdots+b_0(S(V;\s_N)) \leq \Ccal(V;\P).
\end{equation*}
\end{proof}

Of course, higher Betti numbers are not sub-additive, so a similar
procedure cannot be followed in general.

\medskip

\begin{definition}[Consistent sign assignment]
Let $V$ be a Pfaffian variety, $\P$ a family of Pfaffian functions,
and $\Sig$ the set of strict sign conditions on $\P.$ Then $\s \in
\Sig$ is a {\em consistent sign assignment} of $\P$ on $V$ if the
basic set $S(V;\s)$ is not empty.
\end{definition}

Then, $\Ccal(V;\P)$ bounds also the number of consistent sign
assignments $\s \in \Sig.$ Theorem~\ref{th:cells} shows that for a
fixed variety $V,$ $\Ccal(V;\P)$ is a polynomial in the number $s$ of
functions in $\P,$ and thus, the number of consistent sign assignments
is asymptotically much less than the trivial bound of $3^s.$

\subsection{Sub-Pfaffian sets}

\begin{definition}
The set $Y\sub \R^n$ is a {\em sub-Pfaffian set} if there exists a
semi-Pfaffian set $X \sub \U \sub \R^{n+p}$ such that 
$Y$ is the image of $X$ by the canonical projection $\pi: \R^{n+p} \to
\R^n.$ Equivalently, this can be formulated by using an existential formula;
\begin{equation}\label{eq:subpf}
Y=\{y \in \R^n \mid \exists x \in \R^p, \, (x,y) \in X\}.
\end{equation}
\end{definition}

Unlike semi-algebraic sets, semi-Pfaffian sets are not stable under
projections.  

\begin{example}[Osgood~\cite{osgood}]
The following sub-Pfaffian set is not semi-Pfaffian.
\begin{equation*}
X=\{(x,y,z) \in \R^3 \mid \exists u \in [0,1], y=xu, z=xe^u \}.
\end{equation*}
\end{example}

\begin{proof}
The set $X$ is clearly a strict subset of $\R^3$ that contains $0.$ If
$X$ is semi-Pfaffian, there exists a non-zero analytic function $F$
that vanishes on $X$ in a neighbourhood of $0.$ We can write $F$ as a
convergent series of homogeneous polynomials $F_d,$ were
$\deg(F_d)=d.$ Then, we must have for all $u \in [0,1],$
\begin{equation*}
F(x,xu,xe^u)=\sum_{d\geq 0} x^d F_d(1,u,e^u)=0.
\end{equation*}
Thus, we must have $F_d(1,u,e^u)=0$ for all $d\geq 0$ and all $u \in
[0,1],$ which implies $F_d \equiv 0$ for all $d\geq 0.$ Thus, $F \equiv
0$ is the only analytic function that vanishes on $X$ in a
neighbourhood of $0.$ Since $X$ is a strict subset of $\R^3,$ it
cannot be semi-Pfaffian.
\end{proof}

\begin{rem}[Sub-fewnomial sets]
Let $X\sub \R^{n+p}$ be a semi-algebraic set and $Y\sub \R^n$ be its
projection. By the Tarski-Seidenberg theorem, $Y$ is certainly
semi-algebraic too. But even though $X$ may be a fewnomial set, $Y$ is only
sub-fewnomial: describing $Y$ with a fewnomial quantifier-free formula
may not always be possible.
\end{rem}

\begin{example}[Gabrielov~\cite{ga:cex}]\label{ex:sub-few}
Consider for all $m\in \N$ the set
\begin{equation}\label{eq:ga-cex}
Y_m=\{(x,y)\in \R^2 \mid \exists t \in \R, \, t^m-xt=1, \,
(y-t)^m-x(y-t)=1\}.
\end{equation}
Then, there is no quantifier-free fewnomial formula describing $Y_m$
having a format independent of $m.$
\end{example}

\begin{rem}[Open problem]\label{rem:open}
If $Y$ is a sub-Pfaffian set and $Y$ is not subanalytic, it is not
known whether its complement is also sub-Pfaffian or not.  This is one
of the reasons that motivates the introduction of Pfaffian limit sets
in Section~\ref{sec:limitsets}.
\end{rem}

\section{Basic properties of o-minimal structures}

In this section, we describe the main definitions and results
concerning o-minimal structures. O-minimal structures appear in model
theory, and provide a framework for the ideas of {\em tame
topology}~\cite{esquisse}. Many surveys are available to the reader
for more details, for
instance~\cite{coste:pisa,vdd:book,vddm:geom}. 
For more details about model-theoretic aspects, see
also~\cite{vdd:96,vdd:keele}.

\subsection{O-minimal expansions of the real field}

\begin{definition}[o-minimal structure] For all $n\in\N$ let $\S_n$ be
  a collection of subsets of $\R^n,$ and let $\S=(\S_n)_{n\in\N}.$ We
  say that $\S$ is an {\em o-minimal structure on the field} $\R$ if
  the following axioms hold.
\begin{itemize}
\item[{\bf(O1)}] For all $n,$ $\S_n$ is a Boolean algebra.
\item[{\bf(O2)}] If $A\in\S_m$ and $B\in\S_n,$ then $A\times B\in\S_{m+n}.$
\item[{\bf(O3)}] If $A\in\S_{n+1},$ and $\pi$ is the canonical
  projection $\R^{n+1}\to\R^n,$ then $\pi(A)\in\S_n.$
\item[{\bf(O4)}] $\S_n$ contains all the semi-algebraic subsets of $\R^n.$
\item[{\bf(O5)}] All sets in $\S_1$ have a finite number of connected
  components. 
\end{itemize}
\end{definition}

Recall that {\bf (O1)} means that the collections $\S_n$ are stable by
finite intersection, finite unions and taking complements. The axioms
{\bf (O1)} through {\bf (O4)} mean that $\S$ is a {\em structure.}
Axiom {\bf (O5)} is called the {\em o-minimality axiom.}

\begin{definition}[Definability]
Let $\S$ be a structure. If $A\in\S_n,$ we say that $A$
is $\S$-{\em definable}. 
A map $f:A\sub\R^m \to B\sub\R^n$ is called $\S$-{\em definable} if
and only if its graph belongs to $\S_{m+n}.$
\end{definition}

\begin{example}
Let $\S_n$ be the set of semi-algebraic subsets of $\R^n.$ Then,
$\S=(\S_n)_{n \in \N}$ is an o-minimal structure.
\end{example}

\begin{proof}
Recall that a subset of $\R^n$ is {\em semi-algebraic} if it can be
defined by a quantifier-free sign condition on polynomials. Then, $\S$
is clearly a Boolean algebra and is stable by Cartesian
products. Elements of $\S_1$ can only have finitely many connected
components since polynomials in one variable have only finitely many
zeros. Thus, the only non-trivial axiom is~{\bf (O3)}: stability by
projection, which is the result of the classical Tarski-Seidenberg
theorem~\cite{bcr}.
\end{proof}

\begin{example}
Let $\S_n$ be the set of {\em globally subanalytic sets:} subsets of
$\R^n$ that are subanalytic in $\R P^n.$ Then, $\S=(\S_n)_{n \in \N}$
is an o-minimal structure.
\end{example}

\begin{proof}
Here, $\S$ is stable by projections by definition, and the fact that
it is a Boolean algebra follows from Gabrielov's theorem of the
complement~\cite{ga:68}. The axioms 2 and 4 are clear, and the
finiteness of the number of connected components follows from the
local properties of semi-analytic sets~\cite{lo:semi}.
\end{proof}

\bigskip

\begin{definition}[Generated structure]
Let $\S$ be a structure and $\A=(\A_n)_{n \in \N}$ a collection of
subsets of $\R^n$ for all $n\in \N.$ If the closure of $\A$ under 
the Boolean operations, Cartesian product and linear projections is
$\S,$ we say that $\S$ is {\em generated} by $\A.$
\end{definition}

For example, the structure of semi-algebraic sets is generated by the
sets $\{f=0\}$ for all polynomials $f,$ and the structure of globally
subanalytic sets is generated by all the restrictions  $f|_{[-1,1]^n}$ of
all the graphs of functions $f$ that are analytic in a neighbourhood of
$[-1,1]^n.$ 

\bigskip

After the general setting of o-minimal structures was introduced, a
lot of effort was put into constructing new examples. Our main
interest here is the fact that o-minimal functions do generate an
o-minimal structure. This fact, proved first by Wilkie
in~\cite{wilkie99}, is the object of the next section. For now, let us
mention two other cases that seem of particular interest.

\medskip

The structure $\R_{\exp}$ generated
by the exponential function is
o-minimal~\cite{wilkie96}. This is of special interest since it
relates to Tarski's problem about the decidability of real
exponentiation~\cite{Macintyre-wilkie}, which was one of the problems
which first motivated the introduction of o-minimal structures.

\medskip

More recently, Rolin, Speissegger and Wilkie~\cite{rsw:denjoy}
constructed new o-minimal structures using certain quasi-analytic
Denjoy-Carleman classes. This construction allowed to settle two open
problems: (1) if $\A_1$ and $\A_2$ generate o-minimal structures, the
structure generated by $\A_1 \cup \A_2$ is not necessarily o-minimal,
(and thus, there is no 'largest' o-minimal structure), and (2) there
are some o-minimal structures that do not admit analytic cell
decomposition (see Theorem~\ref{th:celldec} and
Remark~\ref{rem:ckdec}).

\bigskip

We will now describe the main properties of o-minimal structures.
Essentially, such a structure has a geometrical and topological
behaviour which is very similar to what is observed in semi-algebraic
sets.  For the remaining of the chapter, $\S$ will be a fixed
o-minimal structure and we will write {\em definable} for
$\S$-definable. In the next chapters, definable will always mean
definable in the o-minimal structure generated by Pfaffian functions.

\subsection{The cell decomposition theorem}

The cell decomposition theorem is an o-minimal analogue of the
cylindrical algebraic decomposition used in real algebraic
geometry~\cite{br,bcr}. Since most features of semi-algebraic sets
follow from that decomposition, they will have an equivalent for
definable sets in o-minimal structures. We fix $\S$ an o-minimal
structure.

\begin{definition}[Cylindrical cell]
Cylindrical cells are defined by induction on the dimension of the
ambient space $n.$ A subset $C$ of $\R$ is a cell if and only if it is
an open interval or a point. A set $C\sub\R^n$ is a cell if and only
if there exists a cell $D\sub \R^{n-1}$ such that one of the two
following conditions is true.
\begin{enumerate}
\item There exists a {\em continuous} definable function $f:D \to \R$
  such that $C$ is one of the following sets,
\begin{align*}
C_0(f)&=\{(x',x_n) \mid x_n=f(x')\},\\
C_+(f)&=\{(x',x_n) \mid x_n>f(x')\},\\
C_-(f)&=\{(x',x_n) \mid x_n<f(x')\}.
\end{align*}
 \item There exists {\em continuous} definable functions $f$ and $g$
 from $D$ into $\R$ such that $f<g$ on $D$ and 
\begin{equation*}
C(f,g)=\{(x',x_n) \mid f(x')<x_n<g(x')\}.
\end{equation*}
\end{enumerate}
\end{definition}

\begin{definition}[Cell decomposition] A cell decomposition of $\R$ is
  a finite partition of $\R$ into open intervals and points. For
  $n>1,$ we say that a finite set $\Ccal$ of cylindrical cells of $\R^n$
  is a {\em cell decomposition} of $\R^n$ if $\Ccal$ is a partition of
  $\R^n$ such that the collection $\{\pi(C)\mid C\in\Ccal\}$ is a cell
  decomposition of $\R^{n-1}.$ (Here again, $\pi$ is the canonical
  projection $\R^n \to \R^{n-1}.$)

If $A_1, \ldots, A_k$ are definable subsets of $\R^n$ and $\Ccal$ is a
cell decomposition of $\R^n,$ the partition $\Ccal$ is said to be {\em
  compatible} with $A_1, \ldots, A_k$ if each $A_i$ is a finite union
of cells in $\Ccal.$
\end{definition}

\begin{theorem}[Cell decomposition theorem]\label{th:celldec}
\begin{enumerate}
\item
Let $A_1, \ldots, A_k$ be definable subsets of $\R^n.$ Then there
exists a cell decomposition of $\R^n$ compatible with $A_1, \ldots,
A_k.$
\item
If $f: A \to \R$ is definable, $A \sub \R^n,$ there exists a cell
decomposition of $\R^n$ compatible with $A$ such that for each cell $C
\sub A,$ the restriction $f|_C$ is continuous.
\end{enumerate}
\end{theorem}

\begin{rem}[$C^k$ cell decompositions]\label{rem:ckdec}
The notions of cells and cell decomposition can be generalized to the
$C^k$ setting for any $k\in\N\cup\{\infty,\omega\}:$ we can define
$C^k$-cells by requiring that the graphs that appear in the definition
of the cells are graph of $C^k$ functions. Then, a $C^k$ cell
decomposition of $\R^n$ would be of course a cell decomposition where
all the cells considered are $C^k.$ For any fixed $k\in \N,$ a $C^k$
analogue of Theorem~\ref{th:celldec} holds for any o-minimal
structure.  Although analytic cell decomposition holds in many known
cases, including the Pfaffian
case~\cite{ls:analytic,vddm:exp,loi:anal}, it does not hold in
general~\cite{rsw:denjoy}. 
\end{rem}

As mentioned previously, Theorem~\ref{th:celldec} can be interpreted as a
generalization of the cylindrical algebraic decomposition. However, it
is worth noting that the proof of Theorem~\ref{th:celldec} is much
more technical (it takes about ten pages in~\cite{vdd:book}). 

\medskip

In the course of proving Theorem~\ref{th:celldec}, the following
theorem is necessary.

\begin{theorem}[Monotonicity theorem]\label{th:monotone}
Let $-\infty \leq a < b \leq \infty$ and let $f: (a,b) \to \R$ be
definable. Then, there exists $a_0=a<a_1<\cdots<a_k=b$ such that on
each interval $(a_i,a_{i+1}),$ the function $f$ is either constant or
strictly monotonous and continuous.
\end{theorem}

\bigskip

The following results are immediate corollaries of the existence of
cell decomposition.

\begin{corollary}
Any definable set has a finite number of connected components. 
\end{corollary}

\begin{proof}
Let $A \sub \R^n$ be definable and $\Ccal$ be a definable cell
decomposition of $\R^n$ compatible with $A.$ Each cell $C\in \Ccal$ is
connected, so the number of connected components of $A$ is at most the
number of cells $C \in \Ccal$ such that $C \sub A.$
\end{proof}

By construction, all cells in a cell decomposition are definably
homeomorphic to a cube $(0,1)^d$ for some $d.$ If $C$ is a cell
homeomorphic to $(0,1)^d,$ we let $d=\dim(C).$ Then, it is natural to
define the dimension of a definable set $A \neq \mpty$ as the maximum of
$\dim(C)$ taken over all cells $C$ contained in $A,$ for a given cell
decomposition $\Ccal$ compatible with $A.$ Then, the following holds.

\begin{proposition}[Dimension is well behaved]
Let $A$ be a definable set, $A \neq \mpty,$ and $f:A \to \R^m$ a
definable map. The dimension of $A$ is well-defined (independent of
the choice of the cell decomposition $\Ccal$) and dimensions verifies
the following properties.
\begin{enumerate}
\item $\dim (\fr A) < \dim(A);$
\item $\dim (f(A)) \leq \dim(A).$
\end{enumerate}
\end{proposition}

\begin{corollary}[Uniform bound on fibers]
Let $A\sub\R^m\times\R^n$ be definable, and define for all $x\in\R^m$
the fiber 
\begin{equation*}
A_x=\{y \in \R^n \mid (x,y)\in A\}.
\end{equation*}
Then, there exists $N$ such that for all $x\in \R^m$ the fiber $A_x$
has at most $N$ connected components.
\end{corollary}

\begin{proof}
Let $\Ccal=\{C_i\}$ be a cell decomposition of $\R^m \times \R^n$ compatible
with $A.$ If $A=C_1 \cup \cdots \cup C_N,$ we have $A_x=(C_1)_x \cup
\cdots \cup (C_N)_x.$ Thus, it is enough to check that each $(C_i)_x$
is connected, which is an easy induction on $n.$ 
\end{proof}

\begin{proposition}[Definable choice]
Let $A \sub \R^m \times \R^n$ be definable. Denote by $\pi$ the
canonical projection $\R^m \times \R^n \to \R^m,$ and let $B=\pi(A).$
Then, there exists a definable $s: B \to \R^n$ such that its graph
$\Gamma(s)$ is contained in $A.$
\end{proposition}

\begin{proof}
Suppose $n=1,$ and let $\Ccal$ be a cell decomposition of $\R^m \times\R$
compatible with $A.$ Fix $x \in B;$ there exists a cell $C \in \Ccal$
such that $C \sub A$ and $x\in\pi(C).$ According to the type of the
cell $C,$ we can define $s(x)$ so that $(x,s(x))\in C.$ 
\begin{itemize}
\item If $C=C_0(f),$ we let $s(x)=f(x);$
\item if $C=C_+(f),$ we let $s(x)=f(x)+1;$
\item if $C=C_-(f),$ we let $s(x)=f(x)-1;$
\item and if $C=C(f,g),$ we let $s(x)=(f(x)+g(x))/2.$
\end{itemize}
This solves the case $n=1,$ since $s$ is certainly a definable
function  $\pi(C) \to \R.$

\medskip

We finish the proof by induction on $n.$ Assume the result holds up to
$n-1,$ and let $\pi'$ be the canonical projection $\R^m \times \R^n
\to \R^m \times \R^{n-1}$ and $A'=\pi'(A).$ By induction, there exists
$s':B \to \R^{n-1}$ definable such that $\Gamma(s') \sub A'.$ We can
use the case $n=1$ for $\Gamma(s')$ now, so there exists $s'':
\Gamma(s') \to \R$ such that $\Gamma(s'') \sub A.$ The projection
$\pi$ restricted to $\Gamma(s'')$ must be a bijection onto $B,$ so
$\Gamma(s'')$ is the graph of a definable function $s:B \to \R^n.$
\end{proof}

\begin{corollary}[Curve lemma]
Let $A \sub \R^n$ be definable and $a \in \fr A.$ Then, there exists
a definable arc $\gamma: (0,1) \to A$ such that $\lim_{t\to 0}
\gamma(t)=a.$ 
\end{corollary}

\begin{proof}
Let $a \in \fr A,$ and let $B=\{|x-a|, x \in A\}.$ The set $B$ is
definable and since $0 \in \ol{B},$ there must be an interval $(0,\e)$
contained in $B.$ By the definable choice theorem above, there exists
a function $\gamma: t \in B \mapsto \gamma(t)\in A$ such that $|a
\gamma(t)|=t.$ By Theorem~\ref{th:monotone}, $\gamma$ is continuous on
an interval $(0,\delta)$ for some $\d\leq \e,$ and by rescaling the
variable $t,$ we can always assume that $\d=1.$
\end{proof}

\subsection{Geometry of definable sets}\label{sec:geom}

We will now list some of the deeper consequences of the o-minimality
axiom. First, asymptotic behaviour of definable functions is very
controlled, and there is a dichotomy between {\em polynomially
bounded} o-minimal structures and structures where $e^x$ is
definable~\cite{mi:exp}.  Although the usual, polynomial
{\L}ojasiewicz inequality does not hold in o-minimal structures that
are not polynomially bounded, the following version does hold.

\begin{theorem}[Generalized {\L}ojasiewicz inequality]
Let $f,g: A \sub \R^n \to \R$ be definable functions such that
$\{f=0\} \sub \{g=0\}$ and $A$ is compact. Then, there exists a
definable $C^p$ function $\phi$ such that $|\phi\circ g(x)| \leq
|f(x)|$ for all $x \in A.$
\end{theorem}

\begin{rem}
For Pfaffian functions, a more explicit inequality, the {\em
exponential} {\L}ojasiewicz inequality, holds. See
Proposition~\ref{prop:exp-loj}.
\end{rem}

\begin{definition}[Stratification]
Let $p\geq 1$ be an integer. A $C^p$ stratification of a set $A$ is a
partition of $A$ into strata such that each stratum is a $C^p$-smooth
submanifold and if $X$ and $Y$ are two strata such that $X \cap \ol{Y}
\neq \mpty,$ then we have $X \sub \ol{Y}.$
\end{definition}

\begin{theorem}[Existence of stratifications]
Let $p\geq 1$ be a fixed integer.
\begin{enumerate}
\item Let $A$ be definable. There exists a definable $C^p$
stratification of $A.$
\item Let $A$ be a closed definable set and $f:A \to \R$ be a
continuous definable function. Then, there exists a definable $C^p$
stratification of $A$ such that for each stratum $X,$ the restriction
$f|_X$ is $C^p$ and of constant rank.
\end{enumerate}
\end{theorem}

\begin{rem} More precise results about stratification with specific
regularity conditions exist: e.g. Whitney, Thom~\cite{loi:thom},
Verdier~\cite{loi:verdier}, etc\ldots
\end{rem}

The next result is about the local triviality of continuous definable
maps. It originated with Hardt in the semi-algebraic case~\cite{hardt}.

\begin{definition}[Trivial map]
Let $f: A \to C$ be a definable map. The map $f$ is called {\em
(definably) trivial} if there exists a definable set $F$ and a
definable homeomorphism $h: A \to C\times F$ such that the following
diagram commutes.
\begin{equation*}
\xymatrix{A\ar[rr]^h \ar[dr]_f & & C\ar[dl]^{\pi_1} \times F \\
& C & 
}
\end{equation*}
where $\pi_1: C \times F \to C$ is the canonical projection.
\end{definition}

\begin{theorem}[Generic triviality]\label{th:generic}
Let $f: A \to C$ be a continuous definable map.
there exists a finite definable
partition $C=C_1 \cup \cdots \cup C_r$ such that $f$ is definably
trivial over each $C_i.$ 
\end{theorem}

\bigskip

We finish our discussion of the general properties of o-minimal
structures with questions about triangulations of sets and maps, which
will play a large role in the proofs of the results of Chapter~5.
For simplicial complexes, we use
the terminology of~\cite{coste:pisa} rather than~\cite{vdd:book}.

\begin{definition}[Simplex]\label{df:simplex}
Let $a_0, \ldots, a_d$ be affine-independent points in $\R^n,$ (not
contained in any $(d-1)$-dimensional affine subspace). We define the
{\em closed} simplex $\sib=[a_0, \ldots, a_d]$ as the subset of
$\R^n$ defined by
\begin{equation}\label{eq:closed-simplex}
\sib=\left\{\sum_{i=0}^d w_i\,a_i \mid \sum_{i=0}^d w_i=1, w_1\geq 0,
\ldots, w_d\geq 0\right\}.
\end{equation}
The {\em open} simplex $\s=(a_0, \ldots, a_d)$ is defined as above,
with the additional condition that all weights $w_i$ are positive.
The points $a_0, \ldots, a_d$ are called {\em vertices} of the (open
or closed) simplex. The {\em dimension} of the simplex is $d.$
\end{definition}

Note that the condition $w_0+\cdots+w_d=1$
in~\eqref{eq:closed-simplex} implies that the weights $w_i$ are
uniquely determined. 

\begin{definition}[Faces]
If $\sib=[a_0, \ldots, a_d]$ is a closed simplex, its faces are all
the closed simplexes of the form $[a_i, i \in I]$ where $I$
is any non-empty subset of $\{0, \ldots, d\}.$
\end{definition}

\begin{definition}[Simplicial complex]\label{df:scomplex}
A {\em (finite\/) simplicial complex} $K$ of $\R^n$ is a finite collection
$\{\sib_1, \ldots, \sib_k\}$ of closed simplices of $\R^n$ such that
the following two conditions hold.
\begin{itemize}
\item For any $i,j\in \{1, \ldots, k\},$ the intersection $\sib_i \cap
  \sib_j$ is a common face of $\sib_i$ and $\sib_j;$
\item $K$ is closed under taking faces.
\end{itemize}
We denote by $|K|$ the subset $\sib_1 \cup \cdots \cup \sib_k$ of $\R^n.$
\end{definition}

\begin{theorem}[Triangulation of compact definable sets]
Let $A \sub \R^n$ be a compact definable set, and $B_1, \ldots, B_k$
be definable subsets of $A.$ There exists a finite simplicial complex
$K$ with vertices in $\Q^n,$ sets $S_1, \ldots ,S_k$ of open simplices
of $K$ and a definable homeomorphism $\Phi: |K| \to A$ such that for
each $i,$ we have $B_i=\cup_{\sigma\in S_i}\Phi(\sigma).$ 
\end{theorem}

Note that the result still holds when $A$ is not compact, provided we
take the weaker notion of simplicial complex where we do not require
$K$ to be closed under taking faces. See~\cite{vdd:book} for more
details.

\bigskip

It is well-known that definable maps are not always triangulable: for
example, the blow-up map $f(x,y)=(x,xy)$ is not.  The result below
says that definable continuous maps from a compact into $\R$ always
are.  Recall that the function $f: A \to \R$ is {\em triangulable} if
there exists a finite simplicial complex $K$ and a homeomorphism
$\Phi: |K| \to A$ such that $f \circ \Phi$ is affine. By affine map,
we mean the following.

\begin{definition}[Affine map]
Let $g: \sib=[a_0, \ldots, a_d] \to \R$ be a function defined on a
simplex. It is {\em affine} if it satisfies the equality
\begin{equation}\label{eq:affine}
g\left( \sum_{i=0}^d w_i\, a_i \right)= \sum_{i=0}^d w_i \,g(a_i);
\end{equation}
for all non-negative $(w_0, \ldots, w_d)$ such that
$w_0+\cdots+w_d=1.$
\end{definition}

The proof of the following theorem
in the o-minimal setting can be found in~\cite{coste:few,coste:pisa}.

\begin{theorem}[Triangulation of functions]\label{th:triang-map}
Let $A \sub \R^n$ be a compact definable subset in an o-minimal
structure $\S$ and $f: A \to \R$ be a definable continuous function.
Then, there exists a finite simplicial complex $K$ in $\R^{n+1}$ and a
definable homeomorphism $\Phi: |K| \to A$ such that $f \circ \Phi$ is
affine on each simplex of $K.$

Moreover, given finitely many definable subsets $B_1, \ldots, B_k$ of
$A,$ we can choose the triangulation $\Phi: |K| \to A$ so that each
$B_i$ is the union of images of open simplices of $K.$
\end{theorem}

\section{Pfaffian functions and o-minimality}
\label{sec:limitsets}

As mentioned earlier, works of Charbonnel~\cite{charb} and
Wilkie~\cite{wilkie99} led to the following result, which we will use
extensively in the present work.
This theorem was then generalized extensively in~\cite{kamac:pf-cl,sp:pf-cl,lr:vol}. 

\begin{theorem}[Wilkie]\label{th:wilkie}
The structure generated by Pfaffian functions is o-minimal.
\end{theorem}

\bigskip

The main result in~\cite{wilkie99} is a theorem of the complement:
Wilkie shows that the Pfaffian structure can be obtained by starting
from semi-Pfaffian sets and iterating the operations of closure under
finite unions, projections and {\em closure at infinity,} where the
last operation consists in considering all sets of the form $A_0 \cap
\ol{A_1} \cap \cdots \cap\ol{A_p}$ for sets $A_0, \ldots, A_p$ already
constructed. The end result is called the {\em Charbonnel closure},
and \cite[Theorem~1.8]{wilkie99} says that the Charbonnel closure
obtained from semi-Pfaffian sets is closed under complementation (and
thus a {\em bona fide} structure) and o-minimal.

\bigskip

This construction of the Pfaffian structure, however, is not very
convenient for quantitative purposes, especially since if $\T_m$
denotes the pre-structure obtained after $m$ iteration and if $X\sub
\R^n$ denotes a definable set that can be constructed within $\T_m,$
there doesn't seem to be any way to derive from Wilkie's work an
upper-bound a the number $p$ such that $\R^n \bs X$ can be constructed
within $\T_{m+p}.$ This is what made it desirable to find an
alternative construction for the Pfaffian structure.

\bigskip

We will now describe in some details the construction of the Pfaffian
structure via limit sets that was suggested by Gabrielov
in~\cite{ga:rc}. Limit sets will provide a notion of format for
arbitrary definable sets, and we will show this format can effectively
be used to derive upper-bounds (Chapter~4 and~5).

\begin{rem}[About sub-Pfaffian sets]\label{rem:sub}
We come back to the open problem evoked in Remark~\ref{rem:open}: if
$(C)$ is the statement: {\em the complement of any sub-Pfaffian set is
again sub-Pfaffian}, we do not know whether $(C)$ holds or not. If we
knew that $(C)$ was true, we could deduce easily that the Pfaffian
structure is o-minimal, since then all definable sets would be
sub-Pfaffian and semi-Pfaffian sets (and thus sub-Pfaffian sets)
always have finitely many connected components.
\end{rem}

\begin{proof}
Let us show that $(C)$ would imply that sub-Pfaffian sets form a
structure.  Since sub-Pfaffian sets are clearly stable under
projection and Cartesian products. Also, sub-Pfaffian sets are always
closed under finite unions since if $X$ and $Y$ are sub-Pfaffian, we
can assume that $X=\pi(X_1)$ and $Y=\pi(Y_1)$ for some semi-Pfaffian
subsets $X_1$ and $Y_1$ of $\R^{n+p},$ with $\pi$ the canonical
projection $\R^{n+p}\to \R^n,$ and thus $X\cup Y=\pi(X_1 \cup Y_1),$ so
it is clearly sub-Pfaffian.

\medskip

Hence, all we have to show is that if $X$ and $Y$ are sub-Pfaffian and
$(C)$ holds, then $X \cap Y$ is sub-Pfaffian too.  Since we're
assuming $(C),$ it is enough to show that the complement $\R^n \bs
(X\cap Y)$ is sub-Pfaffian. But this is obvious, since $\R^n \bs
(X\cap Y)=(\R^n \bs X)\cup (\R^n \bs Y):$ by $(C)$, both $(\R^n \bs X)$
and $(\R^n \bs Y)$ are again sub-Pfaffian, and we have just showed
that sub-Pfaffian sets were stable under finite unions.
\end{proof}

\subsection{Relative closure and limit sets}

From now on, we consider semi-Pfaffian subsets of $\R^n \times \R_+$
with a fixed Pfaffian chain $\f=(f_1, \ldots, f_\ell)$ in a domain $\U$
of bounded complexity.
We write $(x_1, \ldots ,x_n)$ for the coordinates in $\R^n$ and $\lambda$
for the last coordinate (which we think of as a parameter.) If $X$ is
such a subset and $\lambda>0,$ $X_\lambda$ is its fiber
\begin{equation*}
X_\lambda=\{x \mid (x,\lambda) \in X \} \sub \R^n;
\end{equation*}
and we consider $X$ as the family of its fibers $X_\lambda.$
We let 
\begin{equation*}
X_+=X \cap \{\lambda>0\}, \hbox{ and } \ck{X}=\{x \in \R^n \mid (x,0)
\in \ol{X_+}\}.
\end{equation*}
(Thus, $\ck{X}$ is the Hausdorff limit of the family $\ol{X_\lambda}$ 
when $\lambda$ goes to zero.)
The following definitions appear in~\cite{ga:rc}.

\begin{definition}[Semi-Pfaffian family]
\label{df:family}
Let $X$ be a relatively compact semi-Pfaffian subset of $\R^n \times
\R_+.$ The family $X_\lambda$ is said to be a {\em semi-Pfaffian family} if
for any $\e>0,$ the set $X \cap \{\lambda>\e\}$ is restricted. (See
Definition~\ref{df:restricted}.)  The {\em format} $(n,\ell, \alpha,
\beta, s)$ of the family $X$ is the format of the fiber $X_\lambda$ for a
small $\lambda>0.$
\end{definition}

\begin{rem}[Format]
Note that the format of $X$ as a semi-Pfaffian set is different from
its format as a semi-Pfaffian family. We will sometimes refer to it as
the {\em fiber-wise format} to emphasize that fact. Note also
that~\cite{ga:rc} uses the format discussed in
Remark~\ref{rem:basic-format} rather than the formula-based format,
both being of course valid measures of the descriptive complexity of
limit sets.
\end{rem}

\begin{definition}[Semi-Pfaffian couple]
\label{df:couple}
Let $X$ and $Y$ be semi-Pfaffian families in $\U$ with  a common chain
$(f_1, \ldots, f_\ell).$ They form a {\em semi-Pfaffian couple} if the
following properties are verified:
\begin{itemize}
\item $(\bar{Y})_+=Y_+;$
\item $(\partial X)_+ \sub Y.$
\end{itemize}
Then, the {\em format} of the couple $(X,Y)$ is the component-wise
maximum of the format of the families $X$ and $Y.$
\end{definition}

\begin{definition}[Relative closure]
Let $(X,Y)$ be a semi-Pfaffian couple in $\U.$ We define the {\em
relative closure} of $(X,Y)$ at $\lambda=0$ by 
\begin{equation}
\label{def:rc}
(X,Y)_0=\ck{X} \setminus \ck{Y} \sub \ck{\U}.
\end{equation}
We will use the notation $X_0=(X,\fr X)_0.$
\end{definition}

\begin{definition}[Limit set]
Let $\Om \sub \R^n$ be an open domain.  A {\em limit set} in $\Om$ is
a set of the form $(X_1,Y_1)_0 \cup \cdots \cup (X_k,Y_k)_0,$ where
$(X_i,Y_i)$ are semi-Pfaffian couples respectively defined in domains
$\U_i \sub \R^n \times \R_+,$ such that $\ck{\U_i}=\Om$ for $1 \leq i
\leq k.$ If the formats of the couples $(X_i,Y_i)$ is bounded
component-wise by  $(n,\ell, \alpha, \beta, s)$ we say that the
{\em format} of the limit set is $(n,\ell, \alpha, \beta, s,k)$
\end{definition}

\begin{rem}\label{rem:relcompact}
We assumed that the semi-Pfaffian families $X$ are all relatively
compact. This restriction allows us to avoid a separate treatment of
infinity: we can see $\R^n$ as embedded in $\R{}P^n,$ in which case
any set we consider can be subdivided into pieces that are relatively
compact in their own charts.
\end{rem}

\begin{example}
Any (not necessarily restricted\/) semi-Pfaffian set $X$ is a limit
set. 
\end{example}

\begin{proof}
It is enough to prove the result for a basic set $X \sub \U,$ 
\begin{equation*}
X=\{x \in \U \mid \phi_1(x)=\cdots=\phi_I(x)=0, \,\psi_1(x)>0, \ldots,
\, \psi_J(x)>0 \};
\end{equation*}
Let $\psi=\psi_1\cdots \psi_J$ and 
let $g$ be an exhausting function for $\U.$
Define the sets
\begin{align*}
W&=  \left\{(x,\lambda)\in X \times \Lambda \mid g(x)>\lambda\right\};\\
Y_1&=\left\{(x,\lambda)\in \U \times \Lambda \mid \phi_1(x)=\cdots=\phi_I(x)=0,\:
\psi(x)=0,\: g(x) \geq \lambda\right\};\\
Y_2&=\left\{(x,\lambda)\in \U \times \Lambda \mid \phi_1(x)=\cdots=\phi_I(x)=0,\:
g(x) = \lambda\right\}.
\end{align*}
where $\Lambda=(0,1].$ If $Y=Y_1 \cup Y_2,$ it is clear that $(W,Y)$
satisfies the requirements of a semi-Pfaffian couple in
Definition~\ref{df:couple}; its relative closure is clearly $X.$
\end{proof}

For all $n \in \N$ we let $\S_n$ be the collection of limit sets in
$\R^n,$ and $\S=(\S_n)_{n \in \N}.$ The following theorem sums up
the results in~\cite[Theorems~2.9~and~5.1]{ga:rc}.

\begin{theorem}
The collection $\S$ is a structure, and it is o-minimal. Moreover, if
$X$ is a definable set obtained by a combination of Boolean operations
and projections of the limit sets $L_1, \ldots, L_N,$ the set $X$ can
be presented as a limit set whose format is bounded by an effective
function of the formats of $L_1, \ldots, L_N.$
\end{theorem}

Moreover, it is clear that $\S$ coincides with the structure
constructed by Wilkie.
A key result to work with limit sets is the following inequality.

\begin{proposition}[Exponential {\L}ojasiewicz inequality]
\label{prop:exp-loj}
Let $\f$ be a Pfaffian chain of length $\ell$ defined on a domain of
bounded complexity $\U$ for $\f.$ Let $q(x)=Q(x,\f),$ and suppose that
$0 \in \mathrm{cl}(X \cap \{q>0\}).$ Then, there exists $N \in \N$ such that
\begin{equation*}
0 \in \mathrm{cl}(\{x \in X \mid q(x) \geq 1/\exp_\ell(|x|^{-N})\});
\end{equation*}
where $\exp_\ell$ is the $\ell$-th iterated exponential.
\end{proposition}

The proof relies on proving that the rank of the Hardy field generated
by $\f$ at $0$ is bounded by $\ell+1$ (see~\cite{ros:87}).  A detailed
proof can be found in~\cite{ga:rc}, see
also~\cite{gri:sigmo,lion:loj,lmp:diff-eq}.

\subsection{Special consequences of o-minimality}

When giving bounds on the topology of sets defined using Pfaffian
functions, one invokes constantly the o-minimality of the structure
generated by those functions. In this section are gathered a few minor
results that will be often used in the next chapters.

\begin{lemma}[Existence of limits]\label{lem:limit}
Let $f: (0,\e) \to \R$ be definable. %
Then, the function $f$ has a well-defined limit in $\R \cup \{\pm \infty\}.$ 
\end{lemma}

\begin{proof}
This is a simple consequence from the monotonicity theorem
(Theorem~\ref{th:monotone}). There exists $\d>0$ such that the
restriction of $f$ to $(0,\delta)$ is continuous, and one of strictly
increasing, constant, or strictly decreasing. The case where $f$ is
constant on $(0,\delta)$ is trivial. If $f$ is strictly increasing on
that interval, then either it is bounded from above, and then $f$ must
have a finite limit at $0,$ or it is not bounded and the limit of $f$
is $+\infty.$ The case where $f$ is decreasing is similar.
\end{proof}

\begin{lemma}[Critical values]\label{lem:sard}
Let $q: \U \sub \R^n \to \R$ be a $C^\infty$ definable function.
Then, $q$ has finitely many critical values.
\end{lemma}

\begin{proof}
The set of critical values of $q$ is a definable subset of $\R.$ (If
$q$ is a Pfaffian function, this set is actually sub-Pfaffian, see
Example~\ref{ex:semi-sub}.)  By Sard's lemma, it must be of measure
zero, and a definable subset of $\R$ of measure zero can only be
finite.
\end{proof}

Recall that we denote by $b_i(X)$ the $i$-th Betti number of $X$ (see
Notation~\ref{not:betti}) and $b(X)=\sum_i b_i(X).$

\begin{lemma}[Deformation of basic sets]\label{lem:deform}
Let $\U$ be a domain of bounded complexity, $g$ an exhausting function
for $\U$ and let $X=\{x \in \U \mid q_1(x)=\cdots=q_r(x)=0, p_1(x)>0,
\ldots, p_s(x)>0\}$ be a basic semi-Pfaffian set, and for $\e>0$ and
$t\in (\R_+)^s,$ define $X_\e=\{x \in \U \mid q_1(x)=\cdots=q_r(x)=0,
p_1(x)\geq \e t_1,
\ldots, p_s(x)\geq \e t_s, g(x) \geq \e\}.$
Then for all $\e\ll 1,$ $b(X_\e)=b(X).$
\end{lemma}

\begin{proof}
The groups $H_*(X)$ are the direct limit of the singular homology
groups of the compact subsets of $X$ (\cite[Theorem~4.4.5]{spanier}.) 
Thus, $b(X)=\lim_{\e \to 0} b(X_\e),$ and by the generic triviality theorem 
(Theorem~\ref{th:generic}), this sequence is eventually stationary.
\end{proof}

\begin{lemma}[topology of compact limits]\label{lem:proj}
Let $K_\e$ be a decreasing sequence of compact definable sets defined
for $\e>0,$ and let $K$ be their intersection. Then, for all $\e\ll
1,$ and all $0\leq i\leq n,$ we have
\begin{equation*}
b_i(K_\e)= b_i(K) 
\end{equation*}
\end{lemma}

\begin{proof}
Since all the sets considered are triangulable, their homological type
is that of a polyhedron, and the \v{C}ech homology $\check{H}_*$ and
the singular homology $H_*$ are isomorphic. Since the sequence $K_n$
is compact and decreasing and the limit is compact too, we
have~\cite{es}
\begin{equation}\label{eq:inv}
H_*(K)=\lim_{\longleftarrow} H_*(K_\e).
\end{equation}
But by generic triviality, the sets $K_\e$ are homeomorphic for $\e\ll
1,$ hence the limit in~\eqref{eq:inv} becomes eventually stationary.
\end{proof}

\emptyevenpage

%% file: aa2.tex
This chapter is devoted to the study of the possible complexity bounds
that can be proved on the Betti numbers of semi-Pfaffian sets defined
on a domain of bounded complexity. These results include bounds for
the sum of Betti numbers of compact and non-compact Pfaffian varieties
(Theorem~\ref{th:var} and~Theorem~\ref{th:var-nc}), bounds for the sum
of Betti numbers of basic semi-Pfaffian sets (Lemma~\ref{lem:basic})
and semi-Pfaffian sets given by $\P$-closed formulas (the main result
of this chapter, Theorem~\ref{th:Pclosed}). Theorem~\ref{th:cells}
gives a bound on $\Ccal(V;\P),$ the number of connected sign cells of
the family $\P$ on $V$ that was introduced in
Definition~\ref{df:cells}, and this allows to establish in
Theorem~\ref{th:bm} a bound on the Borel-Moore homology of arbitrary
(locally closed) semi-Pfaffian sets.  In particular, this last result
provides an upper-bound for the sum of Betti numbers of any compact
semi-Pfaffian set, without requiring the defining formula to be
$\P$-closed.

\medskip

Recently, Gabrielov and Vorobjov~\cite{gv:qf} generalized the results
of the present chapter: they established a general, single-exponential
bound for the sum of the Betti numbers of {\em any} semi-Pfaffian set,
without any assumption on its topology or defining formula. Such a
result was not known even in the algebraic case, and the precise
statement was added at the end of this chapter (Theorem~\ref{th:qf})
for reference purposes.

\medskip

The setting for the present chapter will be the following: we will
consider a fixed Pfaffian chain $\f$ of length $\ell$ and degree
$\alpha$ in a domain $\U \sub \R^n$ of bounded complexity for the
chain $\f.$ We will let $g$ be an exhausting function for $\U,$ and
$\gamma=\deg_{\f} g.$

\medskip

Throughout this chapter, $p_1, \ldots, p_s,$ and $q_1, \ldots, q_r,$
will be Pfaffian functions in the fixed chain $\f,$ and we'll write
$\P$ for $\{p_1, \ldots, p_s\}$ and $\deg_{\f} \P$ for $\max
\{\deg_{\f} p_i \mid 1\leq i \leq s\}.$ 
The number $\beta$ will be a common upper-bound for $\deg_{\f} p_i$
and $\deg_{\f} q_j.$
We'll let $q=q_1^2+\cdots+q_r^2$ and $V=\Zcal(q_1, \ldots,
q_r)=\Zcal(q).$ The dimension of $V$ will be denoted by $d.$ We'll let
\begin{equation}\label{eq:defV}
\V(n,\ell, \alpha, \beta, \gamma)=
2^{\ell(\ell-1)/2}
\beta(\alpha+\beta-1)^{n-1}\frac{\gamma}{2}
[n(\alpha+\beta-1)+\gamma+\min(n,\ell)\alpha]^\ell.
\end{equation}
The chapter is organized as follows.
\begin{itemize}
\item In the first section, we show that $b(V)$ can be bounded
in terms of $\V(n,\ell,\alpha,\beta,\gamma).$ 
\item In section~\ref{sec:semi-betti}, we show that if $X$ is a
semi-Pfaffian subset of a compact variety $V$ given by a $\P$-closed
formula, $b(X)\leq (5s)^d\, \V(n,\ell, \alpha, 2\beta, \gamma).$
\item In section~2.3, we show that $\Ccal(V;\P) \leq \Sigma(s,d) \,
\V(n,\ell,\alpha, \beta^\star,\gamma),$ where
$\beta^\star=\max(\beta,\gamma)$ and 
\begin{equation*}
\Sigma(s,d)=\sum_{0 \leq i \leq d}\binom{4s+1}{i}.
\end{equation*}
\item The last section is devoted to proving that the rank of the
Borel-Moore homology groups of a locally closed $X$ with the above
format is bounded by an expression of the form
\begin{equation*}
\bbm(X) \leq s^{2d} \, 2^{\ell(\ell-1)} O(n \beta+
\min(n,\ell) \alpha)^{2(n+\ell)};
\end{equation*}
for some constant depending on $\U.$
\end{itemize}

The main results of this chapters are inspired by similar results in
the semi-algebraic case by Basu, Pollack and
Roy~\cite{ba:euler,bpr:cells,bpr:pt}. Note that more analogues could
also be formulated for more recent results in the same
vein~\cite{ba:single,ba:diff,ba:arr}.  Indeed, the o-minimality of the
structure generated by Pfaffian functions ensures that most arguments
can still be used. The use of infinitesimals in those papers can be
avoided most of the time by placing oneself in a compact setting and
replacing the infinitesimals in small real numbers. (The proof of
Theorem~\ref{th:Pclosed} is an example of how one can compute a real
number $r$ such that the condition that $\e$ is an infinitesimal can
be replaced by $\e<r.$)

\bigskip

As in the work of Basu, Pollack and Roy, one of the ideas behind the
bounds is the notion of {\em combinatorial level} of a family of
functions $\P.$%

\begin{definition}[Combinatorial level]\label{df:comblvl}
Let $X \sub \U$ be a semi-Pfaffian set and $\P$ a family of functions
on $\U.$ The {\em combinatorial level} of the couple $(X,\P)$ is the
largest integer $m$ such that there exists $x$ in $X$ and $m$
functions in $\P$ vanishing at $x.$
\end{definition}

This leads to a combinatorial definition of the idea of general
position.

\begin{definition}[General position]
Let $V$ be a Pfaffian variety. The set $\P$ is said to be in {\em
general position} with $V$ if the combinatorial level of $(V,\P)$ is
bounded by $\dim(V).$
\end{definition}

\section{Betti numbers of Pfaffian varieties}

This section is devoted to proving the following analogue for Pfaffian
varieties of the Oleinik-Petrovskii-Thom-Milnor upper
bound~\cite{oleinik,op,thom,milnor:betti} on the Betti numbers of real
algebraic sets. As explained above, $\f$ is a fixed Pfaffian chain of
degree $\alpha$ and length $\ell,$ and $\U$ is a domain of bounded
complexity for $\f$ with an exhausting function $g$ 
such that $\deg_{\f} g=\gamma.$
We fix $V=\Zcal(q_1, \ldots, q_r)$ a Pfaffian variety and let
$q=q_1^2+\cdots+q_r^2.$ The result we will prove is the following.

\begin{theorem}\label{TH:varieties}
Let $V=\Zcal(q)$ a Pfaffian variety with $q$ as above, 
$\deg_{\f} q= 2\beta.$ If $V$ is compact,
its Betti numbers verify
\begin{equation*}%
b(V) \leq \V(n,\ell,\alpha,\beta,\gamma);
\end{equation*}
where $\V(n,\ell,\alpha,\beta,\gamma)$ is defined in~\eqref{eq:defV}.

If $V$ is not compact, we have
$b(V)\leq\V(n,\ell,\alpha,\beta^*,\gamma),$ where
$\beta^*=\max(\beta,\frac{\gamma}{2}).$
\end{theorem}

In practice, it makes sense to assume that the chain $\f$ and the
domain $\U$ are fixed, and to let the degree $\beta$ go to
infinity. We then obtain a more manageable estimate.

\begin{corollary}\label{cor:V-O}
Let $\U$ be a fixed domain of bounded complexity for a Pfaffian chain
$\f.$ If $V=\Zcal(q_1, \ldots, q_r)$ is a Pfaffian variety and
$\deg_{\f} q_i \leq \beta$ for all $i,$ the following asymptotic
estimate holds.
\begin{equation*}
b(V) \leq 2^{\ell(\ell-1)/2}\;O( n \beta + \min(n,\ell) \alpha)^{n+\ell}.
\end{equation*}
(Here, the constant in the O term depends on $\gamma.$)
\end{corollary}

\subsection{Bound for compact Pfaffian varieties}
\label{varieties}

In this section, we assume that the variety $V\sub \U$ is compact.
Recall that if $V=\Zcal(q_1, \ldots, q_r),$ we let $q=q_1^2+\cdots+q_r^2.$
Let $K_r=\{x \in \U \mid q(x) \leq r \}.$ According to
Lemma~\ref{lem:sard}, $q$ has only a finite number of critical values,
and so the $K_r$ are smooth manifolds with boundaries for all but
finitely many values of $r.$ Let $K^*_r \sub K_r$ be the union of the
connected components of $K_r$ that intersect $V.$ We want to show that
$b(V)$ is equal to $b(K^*_r)$ for small values of $r.$
We shall start by proving that $K^*_r$ is compact if $r$ is small
enough. 

\begin{lemma}\label{lem:tube}
  Let $d_V(x)$ be the distance of $x$ to $V,$ and for all $ \delta>0,$ let
  $T( \delta )=\{x\in\U \mid d_V(x)\leq \delta\}.$
There exists $ \delta_1 >0$ such that 
$K^*_r=K_r \cap T(\delta_1 )$ for $r \ll 1.$
\end{lemma}

\begin{proof}
Define, for any set $C,$ $\dist(C,V)=\min\{d_V(x) \mid x\in C\}.$
Let $\delta_0=\dist(\fr\U,V).$ Since $V$ is compact, we have
$\delta_0>0.$ Fix $\d_1\in(0,\d_0).$

\medskip

For all $r>0,$ let $C_r=K_r \bs K^*_r.$ This set is closed for all
$r.$ We will show that $\dist(C_r,V)>\d_1$ when $r \ll 1,$ by
contradiction. If $C_r \cap T(\delta_1)\neq \mpty$ for all $r>0,$
their intersection $\cap_{r>0} (C_r \cap T(\delta_1))$ must be
non-empty too, since those sets are compact. But a point in this
intersection cannot be in $V.$ Since $V=\cap_{r>0} K_r,$ we have a
contradiction. 
\end{proof}

\begin{rem}
  It is important to consider $K^*_r,$ since $K_r$ itself is not
  necessarily compact. The following example comes from~\cite{br}. 
\end{rem}

Let $P:\,\R^2 \rightarrow \R$ be the map:
  \begin{equation*}
    P(x,y)=(x^2+y^2)((y(x^2+1)-1)^2+y^2).
  \end{equation*}
$P^{-1}(0)=\{0\}$ is compact, but as $P(x,(1+x^2)^{-1})$ goes to $0$
as $x$ goes to infinity, the sets $\{P \leq r\}$ are not bounded for
$r>0.$

\bigskip

Since the set $T(\delta_1)$ is compact, Lemma~\ref{lem:tube} implies
that $K^*_r$ is compact for $r \ll 1.$ 
Since $V=\cap_{r>0} K^*_r$ is compact too, we can apply
Lemma~\ref{lem:proj} to conclude that $b(V)=b(K^*_r)$ for $r \ll 1.$
To obtain a bound on $b(V),$ we need to establish a relation between
the topology of $K^*_r$ and the topology of its boundary.

\begin{lemma}
\label{betti:MwB}
 Let $K=K^*_r.$ Then $b(\partial K)=2 b(K).$ 
\end{lemma}

\begin{proof}
  Let $K^c=\R^n \backslash K.$ The
  Mayer-Vietoris sequence in reduced homology of $(K,\overline{K^c})$ is:
\begin{equation} 
\cdots \longrightarrow \tilde{H}_{i+1}(\R^n) \longrightarrow
\tilde{H}_i(\partial K) \longrightarrow \tilde{H}_i(K) \oplus
\tilde{H}_i(\overline{K^c}) \longrightarrow \tilde{H}_i(\R^n)
\longrightarrow \cdots
\end{equation}
As $\tilde{H}_*(\R^n)=0,$ this yields $ \tilde{H}_i(\partial K) \iso
\tilde{H}_i(K) \oplus \tilde{H}_i(\overline{K^c}),$ and as $\partial
K$ has a collar in $\overline{K^c},$ we have
$\tilde{H}_i(\overline{K^c}) \iso \tilde{H}_i(K^c).$

Alexander duality gives $\tilde{H}_i(K^c) \iso \check{H}^{n-i-1}(K).$ 
This yields the relations
\begin{equation}
\label{betti1}
 b_i(\partial K) = b_i(K) + b_{n-i-1}(K), \qquad 0 \leq  i \leq n-1.
\end{equation}
We have $b_n(K)=b_n(\fr K)=0,$ so summing all the equalities in
\eqref{betti1} gives the result $b(\fr K)=2 b(K).$
\end{proof}

\begin{theorem}[Compact varieties]
\label{th:var}
Let $\f$ be a Pfaffian chain of length $\ell$ and degree $\alpha$
defined in a domain $\U$ of bounded complexity $\gamma.$ Let $q_1,
\ldots, q_r$ be Pfaffian functions such that $\deg_{\f} q_i \leq \beta,$
and let $V=\Zcal(q_1,\ldots, q_r).$ If $V$ is {\em compact,} we have
  \begin{equation}
    \label{betti:var}
    b(V) \leq \V(n,\ell, \alpha, \beta, \gamma);
  \end{equation}
where $\V(n,\ell, \alpha, \beta, \gamma)$ is defined in~\eqref{eq:defV}.
\end{theorem}

\begin{proof}
For $r \ll 1,$ we know from Lemma~\ref{lem:proj} that $b(K^*_r)=b(V).$
According to Lemma~\ref{betti:MwB}, it is enough to estimate the Betti
numbers of $W=\partial K^*_r,$ which is a smooth compact manifold for $r
\ll 1.$

\medskip

Up to a rotation of the coordinate system, -- which does not alter the
complexity of $V,$ -- we can assume that the projection map $\pi$ of
$W$ on the $x_1$ axis is a Morse function, {\em i.e.} has only
non-degenerate critical points with distinct critical values.  By
standard Morse theory~\cite{milnor:morse}, $b(W)$ is bounded by the
number of critical points of $\pi,$ which in turn is bounded by the number
of (non-degenerate) solutions of the system;
\begin{equation*}
q(x)-r=\DP{q}{x_2}(x)=\cdots=\DP{q}{x_n}(x)=0.  
\end{equation*} 
The first equation has degree $2\beta$ in the chain $\f,$ and the
others have degree $\alpha+2\beta-1.$ From Theorem~\ref{th:khbc}, such
a system has at most $2\,\V(n\ell,\alpha,\beta,\gamma)$ solutions, and
the bound on $b(V)$ follows.
\end{proof}

\subsection{The case of non-compact varieties}
\label{cells:section:Betti}

Assume now that $V \sub \U$ is not compact. Let $g$ be an exhausting
function for $\U,$ and define for all $\e>0,$ 
\begin{equation}\label{eq:Ve}
V_\e =\{x \in V \mid g(x) \geq \e\}.
\end{equation}
The set $V_\e$ is compact for all $\e>0.$

\begin{proposition}
\label{cells:proposition:compact}
For all $\e \ll 1,$ we have $b(V_\e)=b(V).$
\end{proposition}

\begin{proof}
By generic triviality, there exists $\e_0>0$ such that $g$ restricted
to $V$ is a trivial fibration over $(0,\e_0).$ In particular, this
implies that for $0<\e'<\e<\e_0,$ the inclusion $V_\e \inc V_{\e'}$ is
a homotopy equivalence, and thus $b(V_\e)$ is constant for $\e \in(0,
\e_0).$ By \cite[Theorem~4.4.6]{spanier}, $H_*(V)$ is the inductive
limit of the groups $H_*(V_\e),$ and the result follows.
\end{proof}

\begin{theorem}[Non-compact varieties]
\label{th:var-nc}
Let $\f$ be a Pfaffian chain of length $\ell$ and degree $\alpha$
defined in a domain $\U$ of bounded complexity $\gamma.$ Let $q_1,
\ldots, q_r$ be Pfaffian functions in $\f$ of degree at most $\beta$
and $V=\Zcal(q_1,\ldots, q_r).$ If $V$ is {\em not} compact, we have
  \begin{equation}
    \label{betti:var-nc}
    b(V) \leq \V(n,\ell,\alpha,\beta^*,\gamma);
  \end{equation}
where $\V$ is defined in~\eqref{eq:defV} and $\beta^*=\max(\beta,
\frac{\gamma}{2}).$ 
\end{theorem}

\begin{proof}
Choose $\e>0$ such that $b(V_\e)=b(V),$ and define, for $\w>0$ and
$\eta>0,$ the set $K=\{\w q+g \geq \e-\eta\}.$ Note that $K$ is a
compact subset of $\U.$ 

\medskip

We can choose sequences $\w_\nu$ and $\eta_\nu$ such that the
corresponding sets $K_\nu$ are smooth manifolds with boundary, but
also such that the sequence is decreasing, and that $V_\e=\cap_\nu
K_\nu.$ In order to do that, it is enough to take 
a sequence $\eta_\nu$ that decreases to $0,$ and, if $M_\nu=\max_{K_\nu}
q,$ to choose $\w_\nu \ra_\nu \infty$ such that
$(\w_{\nu+1}-\w_\nu)M_\nu \leq \eta_\nu -\eta_{\nu+1}.$ 

\medskip

Since the decreasing sequence of compacts $K_\nu$ has the compact set
$V_\e$ as a limit, Lemma~\ref{lem:proj} gives that $b(V_\e)=\lim_{\nu
\ra \infty} b(K_\nu),$ and by the same arguments as in
Lemma~\ref{betti:MwB}, we have $2b(K_\nu)=b(\fr K_\nu).$ As in the
proof of Theorem~\ref{th:var}, we reduced our the problem to the one
of estimating the Betti numbers of a compact smooth hypersurface given
by a single Pfaffian equation $\{\w q+g=\e-\eta\}.$ 
This estimate is established by the same method as
in the compact case, by counting critical points of a projection on a
coordinate axis. After a shift of coordinates, we must estimate the
number of non-degenerate solutions of the system
\begin{equation}\label{eq:sysnc}
h(x)-\e+\eta=\DP{h}{x_2}(x)=\cdots=\DP{h}{x_n}(x)=0;
\end{equation} 
where $h(x)=\w q(x)+g(x).$ Since $\deg_{\f} q=2\beta$ and $\deg_{\f}
g=\gamma,$ we must have $\deg_{\f} h \leq \max(2\beta, \gamma),$ and
Khovanskii's bound from Theorem~\ref{th:khbc} gives that the
system~\eqref{eq:sysnc} has at most $2\,
\V(n,\ell,\alpha,\beta^*,\gamma)$ non-degenerate solutions.
\end{proof}

\section{Betti numbers of semi-Pfaffian sets}
\label{sec:semi-betti}

In this section, $\f$ is still a Pfaffian chain defined on a domain of
bounded complexity $\U.$ Let $V$ be a compact Pfaffian variety of
dimension $d$ and $\Phi$ be a $\P$-closed QF formula, with atoms in
a finite set of Pfaffian functions $\P=\{p_1,\ldots, p_s\}.$

\subsection{Going to general position}
\label{general}

 Recall that $\P$ and $V$ are said to be in general position if the
 combinatorial level $(V,\P),$ introduced in
 Definition~\ref{df:comblvl}, is bounded by $d.$ 
The following proposition shows that one can reduce to this case at a
low complexity cost.

\begin{proposition}[General position]
\label{prop:general}
  Let $X=\{x \in V \mid \Phi(x)\},$ where $V$ and $\Phi$ are as
  above. Then there exists a set of $2s$ Pfaffian functions $\P^*$ and
  a $\P^*$-closed QF formula $\Phi^*$ such that the set $X^*=\{x \in V \mid
  \Phi^*(x)\}$ %
  verifies $b(X)=b(X^*).$
Moreover, we have $\deg_{\f}\P=\deg_{\f} \P^*.$
\end{proposition}

\begin{proof}
Of course, the result is non-trivial only if the combinatorial level
of $(V,\P)$ is at least $d+1,$ which implies in particular that $s\geq
d+1.$ 

\medskip

Fix $t \in (\R_+)^s$ and let $\P_\e=\{p_1 \pm \e t_1, \ldots, p_s\pm
\e t_s\}.$ For all $\e>0,$ we build from $\Phi$ a QF formula $\Phi_\e$
which is $\P_\e$-closed, replacing the atoms of $\Phi$ according to
the following procedure.
\begin{itemize}
\item An atom of the form $\{p_i \geq 0\}$ is replaced by 
$\{p_i \geq -\e t_i\};$
\item an atom of the form $\{p_i \leq 0\}$ is replaced by 
$\{p_i \leq \e t_i\};$
\item an atom of the form $\{p_i=0\}$ is replaced by the conjunction 
$\{p_i \leq \e t_i\} \wedge \{p_i \geq -\e t_i\}.$
\end{itemize}
Then, let $X_\e=\{x \in V \mid \Phi_\e(x)\}.$ The sets $X_\e$ are
compact and $X = \cap_{\e >0} X_\e.$ By Lemma~\ref{lem:proj}, there
exists $\e\ll 1$ such that $b(X_\e)=b(X).$

\medskip

Assume that $V$ is a $C^1$-smooth submanifold, and let $p=(p_1,
\ldots, p_s).$ By Sard's lemma, the set of critical values of $p|_V$
has measure zero. Hence, for a generic choice of $(t_1, \ldots, t_s),$
we can find $\e>0$ arbitrarily small such that any element of the form
$(\pm \e t_1 , \ldots, \pm \e t_s)$ is a regular value of $p|_V.$
For such a choice, $(V,\P_\e)$ is in general position and we can take
$X^*=X_\e.$ If $V$ is not a submanifold, it can be stratified as a
disjoint union of submanifolds, and we can choose a $t$ that will work
for every stratum.
\end{proof}

\begin{proposition}[Mayer Vietoris inequalities]\label{prop:MV}
  Let $X_1$ and $X_2$ be two compact semi-Pfaffian  sets. 
Then, for all $i,$ the following inequalities hold.
\begin{align}
    \label{MV1} 
& b_i(X_1)+b_i(X_2)  \leq b_i(X_1 \cup X_2) + b_i(X_1 \cap X_2); \\    
\label{MV2} 
& b_i(X_1 \cup X_2) \leq b_i(X_1)+b_i(X_2)+b_{i-1}(X_1 \cap X_2).
\end{align}
\end{proposition}

\begin{proof}
The Mayer-Vietoris sequence~\cite{bredon} for $X_1$ and $X_2$ is 
the following.
\begin{equation*}
\cdots \rightarrow H_{i+1}(X_1\cup X_2) \rightarrow H_i(X_1\cap X_2)
\rightarrow H_i(X_1)\oplus H_i(X_2) \rightarrow 
H_i(X_1\cup X_2)
\rightarrow \cdots%
\end{equation*}
This sequence is exact when $X_1$ and $X_2$ are compact, and the above
inequalities follow easily.
\end{proof}

\subsection{Betti numbers of a basic open set}

If $\P=\{p_1, \ldots, p_s\},$ the basic open set defined by $\P$ on
the variety $V$ is the set
\begin{equation}\label{eq:basic-open}
X(V;\P)=\{x \in V \mid p_1(x)>0, \ldots, p_s(x)>0 \}.
\end{equation}

\begin{definition}\label{df:B0}
Let $\Bcal_0(n,\ell,\alpha,\beta,\gamma,s,m)$ be the maximum of $b(X)$
where $X=X(V;\P)$ for some set of Pfaffian functions $\P$ on a
Pfaffian variety $V=\Zcal(q_1, \ldots, q_r),$ with the following
conditions. 
\begin{itemize}
\item All functions are polynomial in some Pfaffian chain
$\f$ of length $\ell$ and degree $\alpha,$ defined on a domain $\U$
of bounded complexity $\gamma$ for $\f;$
\item $|\P|=s;$ and the combinatorial level of $(V,\P)$ is $m;$
\item $\deg_{\f} p_i$ and $\deg_{\f} q_j$ are bounded by $\beta$
for $1\leq i \leq s$ and $1\leq j\leq r.$
\end{itemize}
\end{definition}

Then, $\Bcal_0$ admits the following upper-bound.

\begin{lemma}[Basic set bound]
\label{lem:basic}
Let $\Bcal_0(n,\ell,\alpha,\beta,\gamma,s,m)$ be as in
Definition~\ref{df:B0}.
Then, 
  \begin{equation}\label{eq:B0bd}
\Bcal_0(n,\ell,\alpha,\beta,\gamma,s,m)
\leq 2^m \,\binom{s}{m} \,\V(n,\ell, \alpha, \beta, \gamma);
  \end{equation}
where $\V(n,\ell, \alpha, \beta, \gamma)$ is defined
in~\eqref{eq:defV}.
In particular, if $\U$ is fixed and $m\leq n,$ we have
\begin{equation}\label{eq:B0-O}
\Bcal_0(n,\ell,\alpha,\beta,\gamma,s,m)
\leq \binom{s}{m} \, 2^{\ell(\ell-1)/2} O(n\beta+
\min(n,\ell)\alpha)^{n+\ell};
\end{equation}
for a constant depending on $\gamma.$
\end{lemma}

\begin{proof}
  Let $X_\e$ be the set:
  \begin{equation}\label{eq:X_e}
    X_\e = \{x \in V \mid  p_1(x) \geq \e ,
    \ldots , p_s(x) \geq \e \}.
  \end{equation}
Applying Lemma~\ref{lem:deform}, we have $b(X_\e)=b(X)$ for $\e \ll 1.$
Consider now the sets:
\begin{align*}
  T &= \{ x \in V \mid p_2 \geq \e, \ldots, p_s \geq \e \} \supset X_\e, \\
  X^-_\e &=  T \cap \{ p_1 \leq -\e \}, \\
  R &= X_\e \cup X^-_\e, \\
  S &= T \cap \{ -\e \leq p_1 \leq \e \}. \\
  W^+ &= T \cap  \{p_1 =\e \}, \\
  W^- &= T \cap  \{p_1 =-\e \}, \\
  W &=  W^+ \cup W^-.
\end{align*}

As $T=R \cup S$ and $W=R \cap S$ and $R$ and $S$ are compact, the
Mayer-Vietoris inequality~\eqref{MV1} gives: 
\begin{equation*}
  b(R)+b(S) \leq b(T)+b(W).
\end{equation*}
As the union $R = X_\e \cup X^-_\e$
is a disjoint union, the Mayer-Vietoris inequality~\eqref{MV2} gives
$b(R)=b(X_\e)+b(X^-_\e).$ This yields:
\begin{equation}
\label{rec1}
  b(X)=b(X_\e) \leq b(R) \leq b(T)+b(W).
\end{equation}

Let $\P_1=\{p_2, \ldots, p_s\}.$ For $\e \ll 1,$ the set $T$ has the
same Betti numbers as the basic set $X(V;\P_1),$ and
$b(W)=b(X(V^+_1;\P_1))+b(X(V^-_1;\P_1)),$ where $V^+_1=V \cap
\Zcal(p_1+\e)$ and $V^-_1=V \cap \Zcal(p_1-\e).$ The set $\P_1$ has
$s-1$ elements, and the corresponding combinatorial levels are
bounded by $m$ for
$(V;\P_1)$ and by $m-1$ for $(V^+_1;\P_1)$ and $(V^-_1;\P_1).$ Thus, the
relation~\eqref{rec1} gives the following inequality.
\begin{equation*}
\label{rec2}
  \Bcal_0(n,\ell, \alpha, \beta, \gamma,s,m)  \leq 
\Bcal_0(n,\ell, \alpha, \beta, \gamma,s-1,m)+2 \; 
\Bcal_0(n,\ell, \alpha, \beta, \gamma,s-1,m-1).
\end{equation*}

When $s=0,$ we have $X(V;\mpty)=V,$ and when $m=0,$ the functions
$p_i$ have constant sign on $V,$ so that $X(V;\P)=V$ or
$X(V;\P)=\mpty,$ depending on whether all functions of $\P$ are
positive on $V$ or not. Thus, we obtain the following initial
conditions for the induction.
\begin{equation}
\label{rec:init}
\begin{cases}
\Bcal_0(n,\ell, \alpha, \beta, \gamma,0,m) &\leq \V(n,\ell, \alpha,
\beta, \gamma); \\
\Bcal_0(n,\ell, \alpha, \beta, \gamma,s,0) & \leq 
\V(n,\ell, \alpha, \beta, \gamma).
\end{cases}
\end{equation}
We will prove~\eqref{eq:B0bd} by induction on $s,$ for all $m \leq s.$
Assume that
\begin{equation*}
\Bcal_0(n,\ell, \alpha, \beta, \gamma,s-1,m) \leq
2^m \,\binom{s-1}{m} \,\V(n,\ell, \alpha, \beta, \gamma);
\end{equation*}
holds for all integers $m\leq s-1.$ We have:
\begin{align*}
  \Bcal_0(n,\ell, \alpha, \beta, \gamma,s,m) & \leq \Bcal_0(n,\ell, \alpha, \beta, \gamma,s-1,m)+2 \Bcal_0(n,\ell, \alpha, \beta, \gamma,s-1, m-1) \\
          & \leq 2^m \,\binom{s-1}{m} \,\V(n,\ell, \alpha, \beta, \gamma)
+ 2 \cdot 2^{m-1} \, \binom{s-1}{m-1} \,\V(n,\ell, \alpha, \beta, \gamma)\\
          & \leq 2^m \,\binom{s}{m} \,\V(n,\ell, \alpha, \beta, \gamma);
\end{align*}
where the last line follows from the Newton identity. This proves the
estimate~\eqref{eq:B0bd}, and the asymptotic estimate follows easily from
this and Corollary~\ref{cor:V-O}.
\end{proof}

\subsection{Bound for a $\P$-closed formula}

\begin{definition}\label{df:B}
Let $\Bcal(n,\ell,\alpha,\beta,\gamma,s,m)$ be the maximum of $b(X)$
where $X=\{x\in V \mid \Phi(x)\}$ for some $\P$-closed formula having
atoms in a set of Pfaffian functions $\P,$ where  $V=\Zcal(q_1,
\ldots, q_r)$ and the following holds.
\begin{itemize}
\item All functions are polynomial in some Pfaffian chain
$\f$ of length $\ell$ and degree $\alpha,$ defined on a domain $\U$
of bounded complexity $\gamma$ for $\f;$
\item $|\P|=s;$ and the combinatorial level of $(V,\P)$ is $m;$
\item $\deg_{\f} p_i$ and $\deg_{\f} q_j$ are bounded by $\beta$
for $1\leq i \leq s$ and $1\leq j\leq r.$
\end{itemize}
\end{definition}

Recall that the notion of $\P$-closed formula was introduced in
Definition~\ref{df:Pclosed}. It is a quantifier-free formula with
atoms of the form $\{p=0\},$ $\{p\geq 0\}$ and $\{p \leq 0\}$ for
$p\in\P$ that is derived without using negations.

\begin{theorem}%
\label{th:Pclosed-induction}
Let $\Bcal(n,\ell,\alpha,\beta,\gamma,s,m)$ be as in
Definition~\ref{df:B}. Then, the following inequality holds
\begin{equation}\label{eq:Bbd}
\Bcal(n,\ell,\alpha,\beta,\gamma,s,m)\leq
\Bcal_0(n,\ell,\alpha,2\beta,\gamma,s,m)
+ 3s\Bcal(n,\ell,\alpha,\beta,\gamma,3s,m-1);
\end{equation}
where $\Bcal_0$ is as in Definition~\ref{df:B0}.
\end{theorem}

\begin{proof}
Let $\P=\{p_1, \ldots, p_s\}$ be a family of Pfaffian functions and
$X=\{x \in V \mid \Phi(x)\};$ where $V$ is a Pfaffian variety and
$\Phi$ is a $\P$-closed formula, and all the formats fit the
requirements of Definition~\ref{df:B}. We will decompose $X$ into sets
that do not involve the conditions $\{p_i=0\}$

\bigskip

We will start by bounding $b(X)$ with Betti numbers of sets where the
sign condition $\{p_1=0\}$ doesn't appear.  Assume that $m<s,$ and let
\begin{equation*}
\I=\{I=(i_1,\ldots,i_m)\mid 2 \leq i_1 < \cdots <i_s \leq s\};
\end{equation*}
and for all $I \in \I,$ define
\begin{equation*}
Z_I = \{x \in V \mid p_{i_1}(x)=\cdots=p_{i_m}(x)=0\}. 
\end{equation*}
Let $Z=\cup_{I \in \I} Z_I.$ The set $Z$ is compact, and the
restriction of $p_1$ to $Z$ is never zero, or it would contradict the
fact that the combinatorial level of $(V,\P)$ is bounded by $m.$ Thus,
$\e_1=\min_Z|p_1|>0.$ If $m=s,$ $\I$ is empty and we can take any
positive real for $\e_1.$

\medskip
Let $0<\eta_1 < \e_1 / 2.$ 
and consider the sets:
\begin{align*}
  R_1 &=\{x \in V \mid \Phi(x)\wedge
(p_1(x)\leq -\eta_1 \vee  p_1(x)\geq\eta_1)\},\\
  W^+_1 &= \{x \in V \mid \Phi(x) \wedge (p_1(x)=\eta_1)\}, \\
  W^-_1 &= \{x \in V \mid \Phi(x)\wedge (p_1(x)=-\eta_1)\}, \\
  S_1 &= \{x \in V \mid  \Phi(x)\wedge (-\eta_1 \leq p_1(x) \leq \eta_1) \}, \\
  S'_1 &= \{x \in V \mid \Phi(x) \wedge (p_1(x)=0)\}.
\end{align*}

As $X=R_1 \cup S_1$ the Mayer-Vietoris inequality~\eqref{MV2} yields
$b(X) \leq b(R_1)+b(S_1)+b(R_1 \cap S_1),$ and since $R_1 \cap
S_1=W_1^+ \cup W_1^-,$ we obtain $b(X)\leq b(R_1)+b(S_1)+b(W_1^+ \cup
W_1^-).$

\medskip

By Lemma~\ref{lem:proj}, we have for $\eta_1 \ll 1,$
$B(S_1)=b(S'_1),$ and since $W^+_1 \cap W_1^-=\mpty,$ we must have by the
Mayer-Vietoris inequality~\eqref{MV2} again that $b(W_1^+ \cup W_1^-)
\leq b(W_1^+)+ b(W_1^-).$ Thus, we have for $b(X)$ the following bound;
\begin{equation}
  \label{eq:rec-betti}
  b(X) \leq b(R_1)+b(S'_1)+b(W^+_1)+b(W_1^-).
\end{equation}

Introduce the following varieties:
\begin{equation*}
V_1=V\cap \Zcal(p_1), \quad V_1^+=V\cap \Zcal(p_1-\eta_1), \quad
V_1^-=V\cap \Zcal(p_1+\eta_1).
\end{equation*}

$S'_1$ is a semi-Pfaffian subset of the variety $V_1$ given
by sign conditions on $\P_1=\P \backslash \{p_1\}.$ If $V_1 \neq
\mpty,$ $(V_1,\P_1)$ has a combinatorial level that is at most $m-1:$
if there is $x \in V_1$ and $p_{i_1}, \ldots, p_{i_m},$ in
$\P_1=\{p_2,\ldots, p_s\}$ such that $p_{i_1}(x)= \cdots= p_{i_m}(x)=0,$
then $x$ is a point in $V$ be such that $p_1(x)=p_{i_1}(x)=\cdots=
p_{i_m}(x)=0.$ This  contradicts the hypothesis that the
combinatorial level of $(V,\P)$ is bounded by $m.$

\medskip

The sets $W^+_1,$ and $W^-_1$ are semi-Pfaffian subsets given by sign
conditions over the family $\P$ on the varieties $V^+_1$
and $V_1^-$ respectively.  According to the choice made
for $\eta_1,$ those varieties 
do not meet the set $Z=\cup_{I\in\I} Z_I,$ since
each $Z_I$ is a variety obtained by setting exactly $m$ of the
functions in $\{p_2,\ldots p_s\}$ to zero. Thus, 
the combinatorial level of $(V_1^+,\P),$ --  the system over
which $W^+_1$ is defined, -- is bounded by $m-1.$ The same holds for
$W_1^-.$  

\medskip

The above discussion for $W^+_1$ and $W^-_1$ works for $s>m$ only,
but when $m=s,$ the combinatorial levels of $(V^+_1,\P)$ and
$(V^-_1,\P)$ are still bounded by $m-1:$ if $m$ functions were to
vanish at a point $x\in V^+_1,$ $p_1$ would have to be one of them,
but it's impossible since $p_1\equiv\eta_1$ on $V^+_1.$

\bigskip

Thus, the relation~\eqref{eq:rec-betti} bounds $b(X)$ in terms of the
Betti numbers of three sets that have a lower combinatorial level and
one set, $R_1$ that can be defined by a sign condition that does not
involve the atom $\{p_1=0\}.$ 

\bigskip 

Now, the set $R_1$ is defined by sign conditions where the atom
$p_1=0$ doesn't appear any more. We can use a similar treatment
on this set to eliminate the atom $p_2=0$ by defining the following
sets:
\begin{align*}
  R_2 &=\{x \in R_1 \mid \Phi(x)\wedge
(p_2(x)\leq -\eta_2 \vee  p_2(x)\geq\eta_2)\},\\
  W^+_2 &= \{x \in R_1 \mid \Phi(x) \wedge (p_2(x)=\eta_2)\}, \\
  W^-_2 &= \{x \in R_1 \mid \Phi(x)\wedge (p_2(x)=-\eta_2)\}, \\
  S_2 &= \{x \in R_1 \mid  \Phi(x)\wedge 
(-\eta_2 \leq p_2(x) \leq \eta_2) \}, \\
  S'_2 &= \{x \in R_1 \mid \Phi(x) \wedge (p_2(x)=0)\}.
\end{align*}
Here, $\eta_2$ is any positive real number smaller than $\e_2/2,$ where
$\e_2$ is the minimum of $|p_2|$ over all the varieties
given by $m$ equations on $V$ chosen among $p_1=\eta_1,$
$p_1=-\eta_1,$ and $p_i=0$ for $i \geq 3.$ 
Repeating the previous arguments, we obtain
\begin{equation}
\label{induction}
  b(R_1) \leq b(R_2)+b(S'_2)+b(W^+_2)+b(W^-_2)
\end{equation}
when $\eta_2\ll 1.$ Note that again, the sets $S'_2,W^+_2$ and $W^-_2$
are defined with systems that have a combinatorial level at most
$m-1.$ Repeating this until all the functions $p_i$ have been
processed, we end up with a bound of the form:
\begin{equation}
\label{induction2}
  b(X) \leq b(R_s) + \sum\limits_{i=1}^s b(S'_i)+b(W^+_i)+b(W^-_i).
\end{equation}
In this relation, the sets $S'_i,$ $W^+_i$ and $W^-_i$ are all defined
by a system of combinatorial level at most $m-1.$ All that remains is
to estimate $b(R_s).$ We will show now that we can bound $b(R_s)$ by
the sum of Betti numbers of a certain basic semi-Pfaffian set.

\bigskip

Let $C$ be a connected component of $R_s.$ Then, $C$ is contained in
one of the basic closed sets of the form
\begin{equation*}
\{x \in V \mid p_1(x) \geq \pm \eta_1, \ldots, p_s(x) \geq \pm \eta_s
\}.
\end{equation*}
Indeed, if it wasn't the case, there would be points $y$ and $z$ of
$C$ and an index $i$ such that $p_i(y)\leq -\eta_i$ and $p_i(z)\geq
\eta_i.$ The points $x$ and $y$ would be joined by a curve contained
in $C$ that would have to intersect the variety $\Zcal(p_i).$ by
construction, $R_s$ does not meet any of the varieties $\Zcal(p_1),
\ldots, \Zcal(p_s),$ so $C$ is indeed contained in one of the sets of
the form above.

\medskip

Let's assume for simplicity that $C$ is contained in the subset of $V$
defined by $p_i(x)\geq \eta_i$ for all $1 \leq i \leq s.$
By Lemma~\ref{lem:deform}, the equality $b(C)=b(C')$ 
holds when $\eta_1, \ldots, \eta_s$ are small,
where $C'$ is the connected component of the basic set $\{x
\in V \mid p_1(x)>0, \ldots, p_s(x)>0\}$ such that $C \sub C'.$

\medskip

Define the basic set
\begin{equation*}
\Sigma=\{x \in V \mid  p_1^2(x)>0, \ldots, \; p_s^2(x)>0\};
\end{equation*}
and let $\Sigma^*$ be the union of all connected components $D$ of
$\Sigma$ such that $D \cap R_s \neq \mpty.$ Following the above
arguments, we have $b(R_s)=b(\Sigma^*)\leq b(\Sigma)$ for $\eta_1,
\ldots, \eta_s$ small enough. We can thus bound $b(R_s)$ using
Lemma~\ref{lem:basic}, and the inequality~\eqref{eq:Bbd} follows.
\end{proof}

\begin{theorem}\label{th:Pclosed}
For $\Bcal$ as in Definition~\ref{df:B}, we have the following
inequality.
\begin{equation}\label{eq:B-V}
\Bcal(n,\ell,\alpha,\beta,\gamma,s,m)\leq
(5s)^m \, \V(n,\ell,\alpha,2\beta,\gamma).
\end{equation}
In particular, if $X \sub V$ is a semi-Pfaffian subset of a compact
variety $V$ given by a $\P$-closed formula of format
$(n,\ell,\alpha,\beta,s),$ we have
\begin{equation}\label{eq:Betti-Pclosed}
b(X) \leq (10s)^d \, \V(n,\ell,\alpha,2\beta,\gamma).
\end{equation}
\end{theorem}

\begin{proof}
For $m=0,$ no function in $\P$ can change sign on $X,$ so any
connected component of $V$ is either not in $X$ or a connected
component of $X.$ For any space $X,$ its singular homology is the
direct sum of the homology of its connected components
\cite[Theorem~4.4.5]{spanier}. Thus, for $m=0,$ we have $b(X) \leq
b(V),$ and~\eqref{eq:B-V} holds by Theorem~\ref{th:var}, since we
certainly have $\V(n,\ell,\alpha,\beta,\gamma)\leq
\V(n,\ell,\alpha,2\beta,\gamma).$

\bigskip

Assume~\eqref{eq:B-V} holds at rank $m-1.$ 
Using the inductive relation proved in
Theorem~\ref{th:Pclosed-induction} and the bound on $\Bcal_0$ from
Lemma~\ref{lem:basic}, we obtain;
\begin{equation*}
\Bcal(n,\ell,\alpha,\beta,\gamma,s,m) \leq 
\left[ 2^m\,\binom{s}{m}+ 3 \cdot 5^{m-1} s^m \right]  \,
\V(n,\ell,\alpha,2\beta,\gamma).
\end{equation*}

We can bound the binomial coefficient with:
\begin{equation*}
  \binom{s}{m}=\frac{s!}{m!(s-m)!}=\frac{s (s-1)\cdots (s-m+1)}{m!}
\leq s^m;
\end{equation*}
which gives us:
\begin{equation*}
2^m \,\binom{s}{m}+ 3 \cdot 5^{m-1} s^m 
\leq \left(2^m+ 
 3 \cdot  5^{m-1} \right) s^m \leq (5s)^m.
\end{equation*}
This concludes the induction, proving~\eqref{eq:B-V}.

\medskip

The inequality~\eqref{eq:Betti-Pclosed} follows from this and from the
general position argument of Proposition~\ref{prop:general}.
\end{proof}

\begin{corollary}\label{cor:Pclosed}
Let $V\sub \U$ compact Pfaffian variety, $d=\dim(V),$ and let $X=\{x
\in V \mid \Phi(x) \}$ where $\Phi$ is a $\P$-closed Pfaffian
formula. If the format of $X$ is $(n, \ell , \alpha , \beta ,\gamma ,
s),$ the following bound holds.
\begin{equation*}
b(X) \leq s^d \, 2^{\ell(\ell-1)/2} O(n \beta+ \min(n,\ell) \alpha)^{n+\ell};
\end{equation*}
where the constant depends only on $\U.$
\end{corollary}

\begin{proof}
The result follows simply from~\eqref{eq:Betti-Pclosed} and the asymptotic
estimates for $\V$ appearing in Corollary~\ref{cor:V-O}.
\end{proof}

\subsection{Bounds for non-compact semi-Pfaffian sets}

The results of this section can be extended to the case where $V$ is
not compact, including the case $V=\U.$
Let $\Phi$ be a $\P$-closed formula and $V\sub \U$ be a non-compact Pfaffian
variety. Let $X=\{x \in V \mid \Phi(x)\}.$ 

\bigskip

First, define $V_\e=V \cap \{g(x) \geq \e\},$ where $g$ is an exhausting
function for $\U,$ and $X_\e=X \cap V_\e.$ For $\e \ll 1,$ we have
$b(X)=b(X_\e),$ so we are reduced to estimating $b(X_\e).$

\medskip

Proposition~\ref{prop:general} on general position can be repeated
verbatim for $X_\e$ instead of $X,$ so we can construct $X^*_\e$
compact defined by a $\P^*$ closed formula, where the combinatorial
level of $(V,\P^*)$ is bounded by $\dim(V),$ and $|\P^*| \leq 2s+2.$

\medskip

If the combinatorial level of $(V,\P^*)$ is zero, we have $b(X^*) \leq
b(V_\e).$ If $V \neq \U,$ this can be estimated using our bounds on
varieties. If $V=\U,$ Lemma~\ref{betti:MwB} indicates that $b(V_\e)$
can be estimated from $b(\Zcal(g-\e)).$ Since $\Zcal(g-\e)$ is a
compact variety, this last invariant can be estimated by
Theorem~\ref{th:var} without problem.

\medskip

Thus, the inductions can be initiated. The Mayer-Vietoris arguments
also hold in this case: for instance, we can restrict all the set to
$\{g(x) \geq \d\}$ for $\d \ll \e,$ so that we keep compact sets at
all times. 

\bigskip

Thus, analogues of Lemma~\ref{lem:basic} and Theorem~\ref{th:Pclosed}
hold for the case where $V$ is not compact. The precise bounds are
slightly different, but we obtain asymptotic bounds which are
identical to Corollary~\ref{cor:Pclosed}. We will finish this
discussion with the following result.

\begin{corollary}[Complements of $\P$-closed sets]\label{cor:ch2-Alex}
Let $\f$ be a Pfaffian chain defined on a domain $\U$ of bounded complexity.
Let $\Phi$ be a $\P$-closed Pfaffian formula of format $(n,\ell,
\alpha, \beta, s)$ 
and let $X=\{ x\in \U \mid \Phi(x)\}.$ We have
\begin{equation*}
b(\R^n \bs X) \leq s^d \, 2^{\ell(\ell-1)/2} O(n \beta+ \min(n,\ell)
\alpha)^{n+\ell}; 
\end{equation*}
where the constant depends on $\U.$
\end{corollary}

\begin{proof}
We can assume without loss of generality that $X$ is compact. By
Alexander duality~\cite{bredon}, the equality $b(\R^n\bs X)=b(X)+1$
holds, so the result follows from Corollary~\ref{cor:Pclosed}.
\end{proof}

\subsection{Applications to fewnomials}

Now, we apply the results of this section to semi-algebraic sets 
defined in the positive quadrant $(\R_+)^n.$ As explained in
Remark~\ref{rem:log}, we can reduce the problem by a change of
variables to a problem about Pfaffian functions in a chain of length
$r,$ where $r$ is the number of non-zero monomials appearing in the
polynomials defining the set. Thus, the following result follows from
Corollary~\ref{cor:Pclosed}. 

\begin{corollary} 
Let $X \sub V \sub (\R_+)^n$ be defined by a $\P$-closed formula,
where $V=\Zcal(q_1, \dots, q_r)$ and $p_i$ and $q_j$ are polynomials.
If $\dim(V)=d$ and the number of non-zero monomials appearing in the
polynomials $p_i$ and $q_j$ is $r,$ we have
\begin{equation*} 
b(X)\leq s^{d} \; 2^{O(n^2 r^4)}.  
\end{equation*}
\end{corollary}

\section{Counting the number of cells}
\label{sec:cells}

Let $\f$ be a Pfaffian chain defined on a domain of bounded complexity
$\U,$ and let $q_1, \ldots,q_r,$ and $\P=\{p_1, \ldots, p_s\}$ be
Pfaffian functions in $\f.$ We let $V =\Zcal(q_1, \ldots, q_r).$ 
This section is devoted to giving an upper-bound on the number of
cells $\Ccal(V;\P)$ introduced in Definition~\ref{df:cells}.

\bigskip

Recall that we denote by $\Sig$ the set of conjunctions strict sign
conditions $\s$ on $\P.$ For $\s \in \Sig,$  we let $S(V;\s)$ be the
corresponding basic set $\{x \in V \mid \s(x)\}.$ 

\bigskip

To count the number of connected components of $S(V;\s),$ we construct
a variety $V(\s)\sub S(V;\s)$ such that $b_0(V(\s)) \geq b_0(S(V;\s)).$

\subsection{Components deformation}

Fix positive numbers $a_1, \ldots, a_s, b_1, \ldots, b_s, \e$ and $\d,$
and let $\P'=\{p_1 \pm \d a_1, \ldots, p_s \pm \d a_s, p_1 \pm \eta b_1,
\ldots, p_s \pm \eta b_s, g-\e\}.$ We let $V_e=V \cap \{g(x)\geq
\e\},$ and we choose $\e \ll 1$ so that $V_\e$ meets every connected
component of every set $S(V;\s)$ for $\s \in \Sig.$

\bigskip
Fix $\d>0,$ and for any $\s \in \Sig,$ consider the set $C_1(\s)\sub
V_\e$ defined on $\P'$ by replacing any atom $\{p_i >0\}$ of $\s$ by
$\{p_i \geq \d a_i\}$ and any atom $\{p_i <0\}$ by $\{p_i \leq -\d
a_i\}.$

\begin{proposition}\label{prop:c1}
  There is $\d_0>0$ such that for all $\d \leq \d_0$ and for all
  strict sign condition $\s \in \Sig,$ we have $b_0(S(V;\s))\leq
b_0(C_1(\s)).$ 
\end{proposition}

\begin{proof}
It's enough to find a $\d_0$ for a fixed sign condition $\s.$ Clearly,
$C_1(\s) \sub S(V;\s),$ so all we need to do is prove that if $D$ is a
connected component of $S(V;\s),$ it meets $C_1(\s)$ when $\d$ is
small enough.  Fix $x^* \in D \cap V_\e.$ Then, $x^* \in C_1(\s)$ if
and only if for all $i$ such that $p_i(x^*) \neq 0,$ we have
$|p_i(x^*)| \geq \d a_i.$ Since $D \cap V_\e$ is compact, such a
condition will hold for $\d$ small enough.
\end{proof}

Fix $\eta>0$ and for a sign condition $\s \in\Sig,$ define $C_2(\s)$
to be the set defined on $V_\e$ by the following replacement rules: as
in the definition of $C_1(\s),$ any atom $\{p_i >0\}$ of $\s$ by
$\{p_i \geq \d a_i\}$ and any atom $\{p_i <0\}$ by $\{p_i \leq -\d
a_i\}.$ Moreover, the atoms of the type $\{p_i=0\}$ are replaced by
$\{-\eta b_i \leq p_i \leq \eta b_i\}.$

\begin{proposition}\label{prop:c2}
Let $\d \leq \d_0$ be fixed. Then there exists a $\eta_0>0$ such that for
all $\s \in \Sig,$ and for all $\eta \leq \eta_0,$ we have the equality 
$b_0(C_1(\s))=b_0(C_2(\s)).$
\end{proposition}

\begin{proof}
Again, it is enough to prove the result for a fixed $\s \in \Sig.$
The sets $C_2(\s)$ form a decreasing sequence of compacts converging
to $C_1(\s)$ when $\eta$ goes to zero, so the result follows readily. 
\end{proof}

\subsection{Varieties and cells}

The following result allows to extract from sets defined by weak
inequalities varieties that meet every connected components of those sets.

\begin{proposition}
\label{prop:cell2var}
Let $p_1, \ldots, p_s$ be Pfaffian functions, $V$ a Pfaffian
variety and $C$ be a connected component of the set $\{x \in V \mid
p_1(x) \geq 0, \ldots, p_s(x) \geq 0 \}.$ Then, there exists
$I \sub \{1, \ldots, s\},$ (\/possibly empty\/)
such that $C$ contains a connected component of the set $V_I=\{x \in V
\mid p_i(x)=0 \; \forall i \in I\}.$
\end{proposition}

\begin{proof}
Take $I$ a set such that $C \cap V_I \neq \mpty$ and $I$ is
maximal for inclusion. Let $x \in C \cap V_I$ and $D$ be the
connected component of $V_I$ containing $x.$
Assume $D \not \sub C.$
Let $y$ in $D \backslash C:$ there exists an index $j \not \in I$ 
such that $p_j(y)<0.$ Since $j \not \in I,$ and $x \in C,$ it implies
that $p_j(x)>0.$ Let $z(t)$ be a path connecting $x=z(0)$ to $y=z(1)$ in 
$D;$ by the intermediate value theorem, there exists $t_0$ such that
$p_j(z(t_0))=0.$ If $t_0$ is the smallest with this property, we must
have $z(t_0)\in C.$ But we have $p_i(z(t_0))=0$ for all $i \in I \cup
\{j\},$ and that contradicts the maximality of $I$ since $j \not\in I.$
\end{proof}

\begin{proposition}\label{prop:gen}
There exists $a_1, \ldots, a_s, \e$ positive real numbers such with
$\e<\e_0$ such that, for all $0<\d<1$ we can find positive real
numbers $b_1, \ldots, b_s$ for which for all $0<\eta<1,$ the family
\begin{equation*}
\P'=\{p_1 \pm \d a_1, \ldots, p_s \pm \d a_s, p_1 \pm \eta b_1,
\ldots, p_s \pm \eta b_s, g-\e\};
\end{equation*} 
is in general position over $V.$
\end{proposition}

\begin{proof}
This is essentially a repeat of the proof of
Proposition~\ref{prop:general}.
\end{proof}

We can now state the main result.

\begin{theorem}\label{th:cells}
Let $\f$ be a Pfaffian chain defined on a domain
$\U \sub \R^n,$  of bounded complexity $\gamma.$ 
Let $\P=\{p_1, \ldots, p_s\}$ be Pfaffian
functions defined in the chain $\f$ with degree
$\beta$ and $V$ be a Pfaffian variety of dimension $d$ given by
equations of degree at most $\beta$ in the same chain.  Then, 
\begin{equation}\label{eq:cells}
\Ccal(V;\P) \leq  \Sigma(s,d)\, \V(n,\ell,\alpha,\beta^\star,\gamma);
\end{equation}  
where $\beta^\star=\max(\beta,\gamma)$ and
\begin{equation*}
\Sigma(s,d)=\sum_{0 \leq i \leq d}\binom{4s+1}{i}.
\end{equation*}
\end{theorem}

\begin{proof}
According to the results proved in this section, it is enough to bound
the number of connected components of all the sets of the form
$C_2(\s)$ for a suitable choice of the real numbers $a_1, \ldots,
a_s, \e, \d, b_1, \ldots, b_s$ and $\eta.$

\medskip

For a fixed $\s \in \Sig,$ we can bound the number of connected
components of $C_2(\s),$ using Proposition~\ref{prop:cell2var}, by
counting the connected component of all the Pfaffian varieties defined
on $V$ by equations taken among the elements of $\P'.$

According to Proposition~\ref{prop:gen}, we can assume that the sets
of the form $C_2(\s)$ are given by functions which are in general
position over $V.$ Then, we have to count the connected components of
sets of the form
\begin{equation*}
\{ x \in V \mid p_{i_1}= \star_{i_1}, \ldots, p_{i_k}= \star_{i_k} \} 
\hbox{ or }
\{ x \in V \mid g=\e, p_{i_2}= \star_{i_2}, \ldots, p_{i_k}=
\star_{i_k} \}, 
\end{equation*}
where $\star_i \in \{-\d a_i, -\eta b_i, \d a_i, \eta b_i\},$ only for
$0 \leq k \leq d.$  

This gives $\Sigma(s,d)$ possible sets of equations over $V.$ We can then
apply Theorem~\ref{th:var} and the result follows.
\end{proof}

\begin{rem}[Combinatorial lemma]\label{rem:combi}
We have $\Sigma(s,d) \leq (4s+1)^d.$
\end{rem}

\begin{proof}
By definition, $\Sigma(s,d)$ is the number of subsets of cardinality
at most $d$ in a set with $4s+1$ elements. If $f$ is a function from
$A=\{1, \ldots, d\}$ to $B=\{1, \ldots, 4s+1\},$ we have $|f(A)| \leq
d,$ and thus $\Sigma(s,d)$ is bounded by the number of maps $f: A \to B$
which is $(4s+1)^d.$
\end{proof}

\begin{rem}
As explained in Chapter~1, the bound~\eqref{eq:cells} has two
corollaries: it bounds $b_0(X)$ for any semi-Pfaffian set $X,$
and bounds the cardinality of the set of consistent sign assignments:
$\{\s \in \Sig \mid S(V;\s) \neq \mpty\}.$ In particular, note that
for a fixed $d,$ the bound on $\Ccal(V;\P)$ is a polynomial in $s.$
\end{rem}

\begin{corollary}[Fewnomial case]
Let $\K$ be a set of $r$ exponents in $\N^n.$ If $V\sub \R^n$ is a
$d$-dimensional variety defined by $\K$-fewnomials and $\P$ is a set
of $s$ $\K$-fewnomials, the number of cells of $\P$ over $V$ is
bounded by
\begin{equation*}
\Ccal(V;\P) \leq \binom{s}{d}\, 2^{O(n^2 r^4)}.  
\end{equation*}
\end{corollary}

\begin{proof}
Divide $\R^n$ in $2^n$ quadrants and the $n$ coordinate
hyperplanes. By Theorem~\ref{th:cells}, a bound of this type holds for
each quadrant, and we can iterate this on the coordinate
hyperplanes. The number of cells is then bounded by the sum of the
number of cells of the restriction to each set in the partition.
Thus we get 
\begin{equation*}
\Ccal(V;\P) \leq 2^n \, \binom{s}{d}\, 2^{O(n^2 r^4)}
+ n\, \binom{s}{d}\, 2^{O((n-1)^2 r^4)};
\end{equation*}
and the result follows.
\end{proof}

\section{Borel-Moore homology of semi-Pfaffian sets}

We conclude this chapter by estimates on the Borel-Moore Betti numbers
of a locally closed semi-Pfaffian set. 
These estimates follow the techniques that appear
in~\cite{burgisser,momopa,yao} but yield a tighter bound even in the
semi-algebraic case because of our use of the improved bound on the
number of cells (Theorem~\ref{th:cells} in the previous section)
derived from~\cite{bpr:cells}.  In particular, this estimate lets us
bound $b(X)$ for $X$ a compact semi-Pfaffian set which is {\em not}
necessarily defined by a $\P$-closed formula.  However, this estimate
was recently outranked by recent work of Gabrielov and
Vorobjov~\cite{gv:qf}. Their result is stated in Theorem~\ref{th:qf}
for reference purposes.

\medskip

Throughout the rest of the section, we will assume without loss of
generality that all sets under consideration are {\em bounded.}

\bigskip

Recall that a {\em locally closed} subset of $\R^n$ is any set that
can be defined as the intersection of an open set and a closed set. In
particular, any {\em basic} semi-Pfaffian set is locally closed, but a
general semi-Pfaffian set is not necessarily so, since clearly the
subset of $\R^2$ defined by $\{x<0,y<0\} \cup \{x \geq 0, y \geq 0\}$ 
is not locally closed.

\begin{definition}[Borel-Moore homology]
Let $X$ be a locally closed semi-Pfaffian set. We then define its
Borel-Moore homology by
\begin{equation*}
\Hbm_*(X)=H_*(\ol{X},\fr X;\Z).
\end{equation*}
We will denote by $\bbm(X)$ the rank of $\Hbm_*(X).$ 
\end{definition}

Note that when $X$ is {\em compact,} we have of course $b(X)=\bbm(X).$
The key property of Borel-Moore homology is the following result.

\begin{lemma}\label{lem:subadd}
Let $X$ be a locally closed semi-Pfaffian set and $Y \sub X$ be closed
in $X.$ Then, the following inequality holds.
\begin{equation}\label{eq:subadd}
\bbm(X) \leq \bbm(X\bs Y)+\bbm(Y).
\end{equation}
\end{lemma}

\begin{proof}(See also~\cite[\oldS11.7]{bcr}.)
Let $C\sub B \sub A$ be compact definable sets. We can triangulate $A$
so that $B$ and $C$ are subcomplexes of $A.$ This yields an exact
sequence 
\begin{equation}\label{seq1}
\cdots \lra H_{i+1}(A,B) \lra H_i(B,C) \lra H_i(A,C) \lra H_i(A,B)
\lra \cdots
\end{equation}
Now, if $X$ is bounded and locally closed, we have $X=U \cap F$ for
$U$ open and $F$ closed. Thus, we have $\fr X=\fr U \cap F,$ so $\fr
X$ is compact, and since $Y$ is closed in $X,$ the set $\fr X \cup Y$
is compact too. Thus, setting $A=\ol{X},$ $B=\fr X \cup Y$ and $C=\fr
X$ in~\eqref{seq1}, we obtain the exact sequence 
\begin{equation*}%
\cdots \lra H_{i+1}(\ol{X},\fr X \cup Y) \lra H_i(\fr X \cup Y,\fr X)
\lra H_i(\ol{X},\fr X) \lra H_i(\ol{X},\fr X \cup Y) \lra \cdots
\end{equation*}
Let $Y'$ be the interior of $Y$ in $\ol{X}.$ It is easy to check that
$\ol{X \bs Y}=\ol{X} \bs Y'$ and that $\fr (X \bs Y) = \fr X \cup (Y
\bs Y') =(\fr X \cup Y) \bs Y'$ (the last one since $\fr X \cap
Y'=\mpty$).  Thus, by excision, we obtain for all $i$ the following
isomorphism
\begin{equation*}
H_i(\ol{X},\fr X \cup Y) \iso H_i(\ol{X} \bs Y', (\fr X \cup Y) \bs
Y')=\Hbm(X \bs Y).
\end{equation*}
Similarly, since $Y$ is closed in $X,$ we have $\fr Y \sub \fr X,$ so
if $Z$ is the interior of $\fr X$ in $\fr X \cup Y,$ we have $x \in
\fr X \cap \fr Y $ if and only if $ x\not\in Z.$ Hence, one obtains
by excision that
\begin{equation*}
H_i(\fr X \cup Y,\fr X)\iso H_i((\fr X \cup Y) \bs Z ,\fr X \bs Z)
=\Hbm(Y).
\end{equation*}
Thus, we end up with the long exact sequence
\begin{equation*}
\cdots \lra \Hbm_{i+1}(X\bs Y) \lra \Hbm_{i}(Y)\lra \Hbm_{i}(X)
\lra \Hbm_{i}(X \bs Y) \lra \cdots
\end{equation*}
The inequality~\eqref{eq:subadd} then follows easily. 
\end{proof}

This result allows to derive immediately an upper-bound for any basic set.

\begin{proposition}\label{prop:bmbasic}
Let $\P=\{p_1, \ldots, p_s\}$ be a family of Pfaffian functions in a
given chain $\f$ of length $\ell$ and degree $\alpha,$ defined on a
domain $\U$ of bounded complexity. Suppose that the maximum of
$\deg_{\f} p_i$ is bounded by $\beta,$ and let $\s \in\Sig$ be a
strict sign condition on $\P.$ Then, if $V$ is a Pfaffian variety of
dimension $d$ defined by equations of degree bounded by $\beta,$ we
have
\begin{equation}\label{eq:bmbasic}
\bbm(S(V;\s)) \leq s^d \, 2^{\ell(\ell-1)/2} O(n \beta+
\min(n,\ell) \alpha)^{n+\ell};
\end{equation}
where the constant depends only on the domain $\U.$
\end{proposition}

\begin{proof}
Without loss of generality, we can assume that we have
\begin{equation*}
S(V;\s)=\{x \in V \mid p_1(x)=\cdots=p_r(x)=0, p_{r+1}>0, \ldots,
p_s(x)>0\}.
\end{equation*}
Let $q=p_{r+1}\cdots p_s,$ and define the sets
\begin{align*}
X &= \{x \in V \mid p_1(x)=\cdots=p_r(x)=0, p_{r+1}\geq0, \ldots,
p_s(x)\geq0\}; \\
Y &= \{x \in V \mid p_1(x)=\cdots=p_r(x)=0, \, q(x)=0,\, p_{r+1}\geq0, \ldots,
p_s(x)\geq0\}.
\end{align*}
The sets $X$ and $Y$ are closed, with $Y \sub X,$ and we have
$S(V;\s)=X \bs Y.$ Thus, by Lemma~\ref{lem:subadd}, we have
\begin{equation*}
\bbm(S(V;\s)) \leq \bbm(X)+\bbm(Y)=b(X)+b(Y);
\end{equation*}
(since $X$ and $Y$ are compact). The upper-bound follows from
the estimates on the sum of Betti numbers for $\P$-closed formulas
appearing in Corollary~\ref{cor:Pclosed}.
\end{proof}

\begin{theorem}\label{th:bm}
Let $V$ be a Pfaffian variety of dimension $d$ and $X$ be a locally
closed semi-Pfaffian subset of $V$ of format $(n, \ell , \alpha ,
\beta ,\gamma , s).$ The rank of the Borel-Moore homology of $X$ verify
\begin{equation}
\bbm(X) \leq s^{2d} \, 2^{\ell(\ell-1)} O(n \beta+
\min(n,\ell) \alpha)^{2(n+\ell)};
\end{equation}
where the constant depends only on the domain $\U.$ 
\end{theorem}

\begin{proof}
Let $\P=\{p_1, \ldots, p_s\}$ be the set of possible functions
appearing in the atoms of the formula defining $X.$ Using
Lemma~\ref{lem:subadd} twice, we obtain
\begin{equation*}
\bbm(X) \leq \bbm(X \cap \{p_1<0\})+\bbm(X \cap \{p_1=0\})+\bbm(X \cap
\{p_1>0\}). 
\end{equation*}
Repeating this inductively for $p_2, \ldots, p_s,$ we obtain
\begin{equation}\label{eq:thbbm}
\bbm(X) \leq \sum_{\s \in \Sig} \bbm(X \cap S(V;\s)).
\end{equation}
Since $X$ is defined on $V$ by a sign condition on $\P,$ the
intersections $X \cap S(V;\s)$ are either empty or equal to $S(V;\s).$
The bound~\eqref{eq:bmbasic} is known for $\bbm(S(V;\s)),$ and since
the number of terms appearing in the right-hand side
of~\eqref{eq:thbbm} is bounded by the number of cells $\Ccal(V;\P)$ of
$\P$ on $V,$ we can combine the bound from Theorem~\ref{th:cells}
to~\eqref{eq:bmbasic} to obtain the above upper-bound on $\bbm(X).$
\end{proof}

\begin{rem}[Compact case]
As mentioned earlier, if $X$ is a compact semi-Pfaffian set, it is
certainly locally closed and verifies $H_*(X)=\Hbm_*(X),$ and thus,
in this case, Theorem~\ref{th:bm} gives an upper-bound on $b(X).$
\end{rem}

We conclude this chapter by giving, for reference purposes, a very
recent result (Summer 2003) of Gabrielov and Vorobjov. It is the most
general upper-bound known for the sum of the Betti numbers of
semi-Pfaffian set, since it does not have require any hypothesis on
the topology of the set or the shape of the defining formula. As
mentioned earlier, it gives for compact sets a sharper bound than
Theorem~\ref{th:bm}.

\begin{theorem}[Gabrielov-Vorobjov~\cite{gv:qf}]\label{th:qf}
Let $X$ be {\em any} semi-Pfaffian set defined by a quantifier free
formula of format $(n,\ell, \alpha, \beta, s).$ 
The sum of the Betti numbers of $X$ admits a bound of the
form 
\begin{equation}\label{eq:qf}
b(X) \leq 2^{\ell(\ell-1)/2} \, s^{2n} \,%
O(n \beta+ \min(n,\ell) \alpha)^{n+\ell};
\end{equation}
where the constant depends only on the definable domain $\U.$
\end{theorem}

The estimate~\eqref{eq:qf} is obtained by constructing a set $X_1$
defined by a $\P^*$-closed formula, where 
\begin{equation*}
\P^*=\{h_i \mid 1 \leq i \leq s\} \cup
\{h_i^2 - \e_j \mid 1 \leq i,j \leq s\}.
\end{equation*}
For a suitable choice $1 \gg\e_1 \gg \cdots \gg\e_s>0,$ we have
$b(X_1)=b(X),$ and the result then follows from Theorem~\ref{th:Pclosed}.

%% file: aa3.tex
The first section of this chapter is devoted to proving the following
theorem. 

\begin{theorem}\label{th:ss}
Let $f:X \to Y$ be a surjective continuous 
{\em compact covering}\footnote{$f: X \to Y$ is compact covering if and only
if for any compact $L \sub Y,$ there exists a compact $K \sub X$ such
that $f(K)=L.$} 
map.  Then, for all $k \in \N,$ we have
\begin{equation}\label{dim}
b_k(Y) \leq \sum_{p+q=k} b_q(\W_f^p(X));
\end{equation}
where $\W_f^p(X)$ is the $(p+1)$-fold fibered product of $X$ over $f,$
\begin{equation}\label{eq:WpX}
\W_f^p(X)=\{(\xx_0, \ldots, \xx_p) \in X^{p+1} \mid
f(\xx_0)=\cdots=f(\xx_p)\}.
\end{equation}
\end{theorem}

The rest of the chapter contains applications of this result to establish
upper-bounds on the Betti numbers of sub-Pfaffian subsets of the cube.
Applications of this theorem to relative closures will be found in
Chapter~5. 

\bigskip

If $I=[0,1]$ and $X$ is a semi-Pfaffian set in the $N=n_0+\cdots+n_\nu$-dimensional
cube $I^N,$ and $Q_1, \ldots, Q_\nu$ is a sequence of alternating
quantifiers, the set
\begin{equation*}
S=\{\xx_0 \in I^{n_0} \mid Q_\nu \xx_\nu \in I^{n_\nu} \ldots Q_1
\xx_1\in I^{n_1}, (\xx_0, \ldots, \xx_\nu) \in X\};
\end{equation*}
is a sub-Pfaffian set~\cite{ga:sub}. If $X$ is semi-algebraic, then
$S$ is semi-algebraic too. When $\nu$ is small, the bounds established
here for $b(S)$ are better than the previously known bounds coming
from cell decomposition~\cite{collins,gv:cyldec,pv:cyldec} or
quantifier elimination.

\bigskip

The chapter is organized as follows: the first section describes the
construction of a spectral sequence $E^r_{p,q}$ that gives
Theorem~\ref{th:ss}. Section~2 contains various topological lemmas,
and section~3 applies the theorem to the case of a set defined by one
quantifier block (Theorem~\ref{th:exist} and
Corollary~\ref{cor:univ}), initiating the induction. In section~4, we
establish an inductive relation for the bound for $\nu$ quantifiers
(Theorem~\ref{th:Erec}), and use it to deduce general upper-bounds
(Corollary~\ref{cor:EPfaff} for the Pfaffian case and
Corollary~\ref{cor:Ealg} for the algebraic case). For semi-algebraic
sets, a comparison is presented in section~5 between those bounds and
the previously available ones using quantifier elimination.

\section{Spectral sequence of a surjective map}

For a {\em closed}\/ surjection $f,$ Theorem~\ref{th:ss} will be proved
in the following way. We will construct a space $J^f(X)$ which is
homotopically equivalent to $Y,$ and has a natural filtration
$J^f_p(X).$ This filtration gives rise to a spectral sequence $E^r_{p,q}$
converging to the homology of $J^f(X),$ as described in the
Appendix. The first term of the sequence  $E^1_{p,q}$ is isomorphic to
$H_q(\W_f^p(X)),$ which will prove the result.
The convergence for a compact-covering map will be deduced from the
closed case in Theorem~\ref{th:ccss}.

\bigskip

In this section, $X,Y$ and $P_0, \ldots P_p$ are topological spaces,
and $f_i: P_i \to Y$ are continuous surjective maps. We denote by
$\Delta^p$ the standard $p$-simplex	
\begin{equation*}
\Delta^p=\{s=(s_0, \ldots, s_p) \in \R^{p+1} \mid s_0\ge
0,\dots,s_p\ge 0,\;s_0+\dots+s_p=1\}.
\end{equation*}

\begin{definition}[Join]\label{joins}
For a sequence $(P_0,\dots,P_p)$ of topological spaces, their {\em
join} $P_0*\dots*P_p$ can be defined as the quotient
\begin{equation*}
P_0\times\dots\times P_p\times\Delta^p/\sim;
\end{equation*}
where $\sim$ is the {\em join} relation 
\begin{equation}\label{equiv}
(x_0,\dots,x_p,s)\sim (x'_0,\dots,x'_p,s') 
\hbox{ iff }
s=s' \hbox{ and } (s_i \neq 0 ) \Rightarrow (x_i=x'_i).
\end{equation}
\end{definition}

\medskip

Recall that if for all $i,$ $f_i: P_i\to Y$ is a continuous surjective
map, we can define the {\em fibered product} 
\begin{equation*}
P_0\times_Y\dots\times_Y P_p
=\{(x_0 , \ldots, x_p) \in P_0\times\cdots\times P_p\mid  
f_0(x_0)=\cdots=f_p(x_p)\}.
\end{equation*}
Note that there is a natural map
\begin{align*}
F: P_0\times_Y\dots \times_Y P_p & \to  Y \\
   (x_0, \ldots, x_p) & \mapsto  f_i(x_i);
\quad \hbox{(taking any $0 \leq i \leq p.$)}
\end{align*}

\begin{definition}[Fibered join]
For $P_0, \ldots, P_p$ as above, we define the {\em fibered join}
$P_0*_Y\dots *_Y P_p$ as the quotient space of
$P_0\times_Y\dots\times_Y P_p\times\Delta^p$ over the join relation
\eqref{equiv}.
\end{definition}

\medskip

The map $F: P_0\times_Y\dots \times_Y P_p \to Y$ extends naturally to
$P=P_0*_Y\dots *_Y P_p.$ Indeed, for any $(x_0 , \ldots, x_p)$ and
$(x'_0,\dots,x'_p)$ in $P_0\times_Y\dots \times_Y P_p \to Y$ such that
$(x_0,\dots,x_p,s)\sim (x'_0,\dots,x'_p,s)$ for some $s \in \Delta^p,$
we must have $x_i=x'_i$ for some $0 \leq i \leq p,$ and thus $F(x_0 ,
\ldots, x_p)=F(x'_0,\dots,x'_p).$ We still denote this map by $F.$ For
any point $y\in Y$ the fiber $F^{-1}(y)$ coincides with the join
$f_0^{-1}(y) *\dots* f_p^{-1}(y)$ of the fibers of $f_i$.

The other natural map is the projection $\pi: P \to \Delta^p.$ If $s$
is in the interior of $\Delta^p,$ the equivalence relation $\sim$ is
trivial over $s,$ so we must have
\begin{equation*}
\forall s \in \mathrm{int}(\Delta^p), \quad
\pi^{-1}(s)=P_0\times_Y\dots \times_Y P_p.
\end{equation*}

\bigskip

For $0 \leq i \leq p,$ define $Q_i$ to be the fibered join 
\begin{equation*}
Q_i=P_0*_Y\dots*_Y P_{i-1}*_Y P_{i+1}*_Y\dots*_Y P_p.
\end{equation*}
Then, one can define a map $\phi_i: Q_i \to P$  by
\begin{equation*}
\phi_i(y_0, \ldots, y_{i-1}, y_{i+1}, \ldots, y_p,t)=
(y_0, \ldots, y_{i-1}, x_i, y_{i+1}, \ldots, y_p,s)
\end{equation*}
where $x_i \in P_i$ is any point and
\begin{equation*}
s_j=
\begin{cases}
t_j &\hbox{ if }  j<i; \\
0  &\hbox{ if } j=i; \\
t_{j-1}  &\hbox{ if } j>i.
\end{cases}
\end{equation*}

\begin{lemma}\label{lem:join101}
The map $\phi_i$ is an embedding $Q_i \to P,$ and
$\phi_i(Q_i)=\pi^{-1}\{s_i=0\}.$ 
Moreover, the space
\begin{equation*}
P/\left(\bigcup_i \phi_i (Q_i)\right);
\end{equation*}
 is homotopy equivalent to the $p$-th suspension%
of $P_0\times_Y\dots\times_Y P_p$.
\end{lemma}

\begin{proof}
The map $\phi_i$ sends $Q_i$ to points $(x_0, \ldots, x_p,s)$ of $P$
such that $s_i=0$ by construction. This means also that $\phi_i$ does
not depend on the choice of the point $x_i \in P_i.$ 
\end{proof}

\bigskip

We now consider the case of a continuous surjection $f: X \to Y.$ For
all $p \in \N,$ define $J^f_p(X)$ to be the fibered join of
$p+1$ copies of $X$ over $f,$
\begin{equation}\label{jeq}
J^f_p(X)=\underbrace{X*_Y\dots*_Y X}_{p+1\> \>{\rm times}}.
\end{equation}

\begin{definition}[Join space]\label{joindef}
If $f$ is as above, the {\em join space} $J^f(X)$ is the quotient
space 
\begin{equation}
\bigsqcup_p J^f_p(X) / \sim;
\end{equation}
where we identify, for all $p \in \N,$
\begin{equation*}
J^f_{p-1}(X) \sim \phi_i(J^f_{p-1}(X)) \sub J^f_p(X), \quad 
\hbox{for all } 0 \leq i \leq p.
\end{equation*}
When $Y$ is a point, we write $J_p(X)$
instead of $J^f_p(X)$ and $J(X)$ instead of $J^f(X)$.
\end{definition}

\begin{lemma}\label{contract}
Let $\phi:J_p(X)\to J(X)$ be the natural map induced by the maps
$\phi_i$.  Then $\phi(J_{p-1}(X))$ is contractible in $\phi(J_p(X))$.
\end{lemma}

\begin{proof}
Let $x$ be a point in $X$.  For $t\in[0,1]$, the maps
$g_t(x,x_1,\dots,x_p,s)\mapsto (x,x_1,\dots,x_p,ts)$ define a
contraction of $\phi_0(J_{p-1}(X))$ to the point $x\in X$ where $X$ is
identified with its embedding in $J_p(X)$ as $\pi^{-1}(1,0,\dots,0)$.
It is easy to see that the maps $g_t$ are compatible with the
equivalence relations in Definition \ref{joindef} and define a
contraction of $\phi(J_{p-1}(X))$ to a point in $\phi(J_p(X))$.
\end{proof}

\begin{proposition}\label{prop:JX}
Let $f:\ X \to Y$ be a closed surjective continuous map, where $X$ and
$Y$ are definable in an o-minimal structure. If $J^f(X)$ is the join
space introduced in Definition~\ref{joindef}, we have
\begin{equation*}
H_*(J^f(X)) \iso H_*(Y).
\end{equation*}
\end{proposition}

\begin{proof}
Let $F:J^f(X)\to Y$ be the natural map induced by $f$.  Its fiber
$F^{-1}y$ over a point $y\in Y$ coincides with the join space
$J(f^{-1}y)$.  
According to Lemma~\ref{contract}, $\phi(J_{p-1}(f^{-1}y))$ is
contractible in $J(f^{-1}y)$, for each $p,$ so $\HB^*(J(f^{-1}y))=0.$ 

\medskip

Since $J(f^{-1} y)$ is a definable quotient, it is locally
contractible, so by Proposition~\ref{prop:Hbar}, we also have
for the reduced Alexander cohomology 
$\HBT^*(J(f^{-1}y))=0.$

\medskip

Since $F: J^f(X) \to Y$ is a closed continuous surjection with fibers
that are trivial for the Alexander cohomology, we can
apply to $F$ the 
Vietoris-Begle theorem (Theorem~\ref{th:VB}) to obtain
$\HB^*(J^f(X)) \iso \HB^*(Y),$ which implies that 
$H_*(J^f(X)) \iso H_*(Y).$
\end{proof}

\begin{notation}[Fibered products]\label{not:Wp}
Throughout this chapter, for $f: X \to Y$ a continuous surjection,
we will denote by $\W_f^p(X)$ the $(p+1)$-fold fibered product of $X$
over $f,$
\begin{equation*}
\W^p_f(X)= \underbrace{X\times_Y \cdots \times_Y X}_{p+1 \,
\mathrm{times}}=\{(\xx_0, \ldots, \xx_p) \in (X)^{p+1} \mid 
f(\xx_0)=\cdots=f(\xx_p)\}.
\end{equation*}
Moreover, if  $A \sub X,$ we will denote by $\W_f^p(A)$ the
corresponding fibered product for the restriction $f|_A.$
\end{notation}

\begin{theorem}[Spectral sequence, closed case]\label{th:spectral-f}
Let $f:\ X \to Y$ be a closed surjective continuous map, where $X$ and
$Y$ are definable in an o-minimal structure.
Then, there exists a spectral sequence $E_{p,q}^r$
converging to $H_*(Y)$ with
\begin{equation}\label{sseq}
E_{p,q}^1 \iso H_q(\W_f^p(X))
\end{equation}
\end{theorem}

\begin{proof}
By Theorem~\ref{th:filtr-ss}, the filtration of $J^f(X)$ by the
spaces $J^f_p(X)$ gives rise to a spectral sequence $E^r_{p,q}$
converging to $H_*(J^f(X)),$ which, by Proposition~\ref{prop:JX}, is
isomorphic to $H_*(Y).$

\medskip

The first term of the sequence is $E^1_{p,q}=H_{p+q}(A_p),$ where
\begin{equation*}
A_p=J^f_p(X)/\left(\bigcup_{q<p} \quad J^f_q(X)\right).
\end{equation*}
From Lemma~\ref{lem:join101}, the space $A_p$
is homotopy equivalent
to the $p$-th suspension of $\W_f^p(X),$ and thus we have $E^1_{p,q}
\iso H_q(\W_f^p(X)),$ proving the theorem.
\end{proof}

\bigskip

One of the features of spectral sequences is that the rank of the
limit of a spectral sequence is controlled by the rank of the initial
terms. Such estimates are discussed in more details in
Corollary~\ref{cor:ss-bound} in the appendix, and applying them to the
present situation yields the estimates of Theorem~\ref{th:ss} for a
closed $f.$

\begin{rem}
The condition of o-minimality is not really important here. 
If $X$ is the difference between a finite CW-complex and
one of its subcomplexes, and $Y$ is of the same type, the spaces $J(f^{-1}y)$
are still locally contractible and the result still holds.
\end{rem}

\medskip

For a locally split map
\footnote{A map $f: X\to Y$ is {\em locally split} if it admits continuous %
sections defined around any point $y\in Y.$ In particular, %
the projection of an open set is always locally split.} %
 $f$, the convergence of the same spectral sequence can be derived from
\cite[Corollary~1.3]{dugger}. 
However, this follows from a more general result: the spectral
sequence converges when $f$ is a {\em compact covering} map. Let us
first recall the definition.

\begin{definition}[Compact covering]\label{df:cc}
A map $f: X \to Y$ is called {\em compact covering} if for all compact
$L \sub Y,$ there exists a compact $K \sub X$ such that $f(K)=L.$
\end{definition}

Note that if $f$ is {\em closed} or {\em locally split,} it is
necessarily compact-covering.

\begin{theorem}[Spectral sequence, compact covering case]\label{th:ccss}
Let $f: X \to Y$ be a definable, compact covering surjection.  
Then, there exists a spectral sequence $E_{p,q}^r$
converging to $H_*(Y)$ with
\begin{equation}\label{bla}
E_{p,q}^1 \iso H_q(\W_f^p(X)).
\end{equation}
\end{theorem}

\begin{proof}
Recall that the singular homology of a space is isomorphic to the
direct limit of its compact subsets \cite[Theorem~4.4.6]{spanier}.
Since $f$ is compact covering, if $K$ and $L$ range over all compact subsets
of $X$ and $Y$ respectively, the following inductive limits verify
\begin{equation}\label{indY}
\limind H_*(f(K))=\limind H_*(L)\iso H_*(Y).
\end{equation}
Let $p$ be fixed and $L_p$ be a compact subset of the fibered product
$\W^p_f(X).$ If for all $0 \leq i \leq p,$ $\pi_i$ denotes the
canonical projection $(\xx_0, \ldots, \xx_p) \mapsto \xx_i,$ we let 
$K_p=\pi_0(L_p)\cup \cdots \cup \pi_p(L_p).$ Observe then that the set
$\W^p_f(K_p)$ 
is a compact subset of $\W^p_f(X)$ containing $L_p.$ Thus, we also
have the following equality
\begin{equation}\label{indW}
\limind H_*(\W^p_f(K_p))=\limind H_*(L_p)\iso H_*(\W^p_f(X)).
\end{equation}
For any compact subset $K$ of $X,$ the restriction $f|_{K}$ is closed,
so by Theorem~\ref{th:spectral-f}, there exists a spectral sequence
$E^r_{p,q}(K)$ that converges to $H_*(f(K))$ and such that
$E^1_{p,q}(K)\iso H_q(\W^p_f(K)).$ By~\eqref{indY} and~\eqref{indW},
the direct limit of $E^r_{p,q}(K)$ when $K$ ranges over all compact
subsets of $X$ is a spectral sequence converging to $H_*(Y)$ 
and verifying~\eqref{bla}.
\end{proof}

\begin{example}
Note that without an additional assumption on $X$ and $Y,$ the
spectral sequence may not converge to $H_*(Y).$ For instance, consider
for $X$ any open segment in $\R^3.$ Let $\{a,b\}=\fr X,$ assume $a\neq
b$ and let $f$ be any projection such that $f$ is 1-to-1 on $X$ and
there exists $c\in X$ such that $f(a)=f(b)=f(c).$ Then, if $Y=f(X),$
we have $b_1(Y)=2,$ but since $X$ is contractible, $b_1(X)=0,$ and
since $f$ is 1-to-1 on $X,$ we have $b_0(\W^1_f(X))=1,$ so $b_1(Y)>
b_1(X)+b_0(\W^1_f(X)).$ The inequality of Theorem~\ref{th:ss} does not
hold in this case.
\end{example}

\begin{rem}
For a map $f$ with 0-dimensional fibers, a similar spectral sequence,
called ``image computing spectral sequence'', was applied to problems
in theory of singularities and topology by Vassiliev \cite{vassiliev},
Goryunov-Mond \cite{gor-mond}, Goryunov \cite{goryunov}, Houston
\cite{houston}, and others. In sheaf cohomology, the corresponding
spectral sequence is known as {\em cohomological
descent}~\cite{deligne}.
\end{rem}

\section{Topological lemmas}

Throughout the rest chapter, $I$ denotes the closed interval $[0,1]$
and {\em open} and {\em closed} are meant in a
cube $I^m$ (for some $m$).

\begin{lemma}\label{lem:open-closed}
Let $X \sub I^{n+p}$ be closed (resp. open). Then, the sets 
\begin{equation*}
Y=\{ \y \mid \exists \xx\in I^p, \, (\xx,\y) \in X\},
\hbox{ and }\quad
Z=\{ \y \mid \forall \xx\in I^p, \, (\xx,\y) \in X\};
\end{equation*}
are both closed (resp. open).
\end{lemma}

\begin{proof}
Let $\pi$ be the canonical projection $\R^{n+p} \to \R^n.$ The sets
$Y$ and $Z$ can be defined by $Y=\pi(X)$ and $I^n \bs Z = \pi(I^{n+p}
\bs X).$ Since $\pi$ is continuous, it sends closed sets to closed
sets, since any closed subset of a cube is compact. Moreover, $\pi$
also sends open sets to open sets. The result then follows easily.
\end{proof}

\begin{lemma}[Alexander duality in the cube~1]\label{lem:Alex1}
Let $X \sub I^n$ be a definable open set. For any $0 \leq q \leq n-1,$
we have
\begin{equation}\label{eq:duality1}
\tilde{H}_q(X \cup J^n) \iso H^{n-q-1}(I^n\bs X);
\end{equation}
where $J^n=(-\e, 1+\e)^n \bs I^n$ for some  $\e>0.$
\end{lemma}

\begin{proof}
Let $S^n=\R^n \cup \{\infty\}$ be the one-point compactification of
$\R^n,$ and let $K$ be the complement of $X \cup J^n$ in $S^n.$ Since
$K$ is closed and not empty,we have by Alexander duality in
$S^n$~\cite[Corollary~VI.8.6]{bredon};
\begin{equation}\label{eq:dual1}
\tilde{H}_q(X \cup J^n)=\tilde{H}_q(S^n \bs K) \iso
\tilde{\ck{H}}^{n-q-1}(K);
\end{equation}
and since $K$ is triangulable, the right-hand side of this equation is
isomorphic to $\tilde{H}^{n-q-1}(K).$ If $C_0, \ldots, C_N$ are the
connected components of $K,$ where $\infty \in C_0,$ then $C_0=S^n \bs
(-\e, 1+\e)^n$ is contractible and $C_1 \cup \cdots \cup C_N=I^n \bs
X.$ This implies that
\begin{equation*}
H^*(K) \iso H^*(I^n \bs X) \oplus H^*(C_0);
\end{equation*}
and as $C_0$ is contractible, we obtain that $\tilde{H}^*(K) \iso
H^*(I^n\bs X).$ Substituting this result in~\eqref{eq:dual1} gives the
lemma.
\end{proof}

To prove a similar result in the case where $X$ is closed, we will
need the following lemma. 

\begin{lemma}\label{lem:complement}
Let $I_0$ be the open interval $(0,1),$ and let $X \sub I^n$ be
closed. Then, we have $H_*(I_0^n \bs X)\iso H_*(I^n\bs X).$
\end{lemma}

\begin{proof}
Consider for $\d>0$ the set $X_\d$ defined by
\begin{equation*}
X_\d=\{x \in I^n \mid \dist(x,X)<\d\}.
\end{equation*}
Since $I^n\bs X_\d$ is a compact subset of $I^n\bs X,$ we have
$H_*(I^n\bs X) \iso H_*(I^n\bs X_\d)$ for $\d\ll 1.$ (The proof uses
the same arguments as the proof of Lemma~\ref{lem:deform}.) But
clearly $I^n \bs X_\d$ is homotopy equivalent to $I_0^n \bs X_\d,$ and
we also have $H_*(I_0^n \bs X_\d)\iso H_*(I_0^n \bs X)$ for $\d\ll1.$
\end{proof}

\begin{lemma}[Alexander duality in the cube~2]\label{lem:Alex2}
Let $X \sub I^n$ be a definable closed subset.  For any $0 \leq q \leq
n-1,$ we have
\begin{equation}\label{eq:duality}
H_q(I^n \bs X)\iso \tilde{H}^{n-q-1}(X \cup \fr I^n).
\end{equation}
\end{lemma}

\begin{proof}
Let us consider the compact set $K=X \cup \fr I^n.$ As in the proof of
Lemma~\ref{lem:Alex1}, Alexander duality in $S^n$ gives
\begin{equation}\label{eq:dual2}
\tilde{H}_q(S^n \bs K) 
\iso \tilde{\ck{H}}^{n-q-1}(K).
\end{equation}
Again, the right-hand side above is isomorphic to
$\tilde{H}^{n-q-1}(K),$ since $K$ is triangulable.  Let $C_0, \ldots,
C_N$ be the connected components of $S^n \bs K,$ with $\infty \in
C_0.$ The component $C_0$ is simply $S^n \bs I^n,$ and thus is
contractible. As before, we can derive from this that $\tilde{H}_*(S^n
\bs K) \iso H_*(I^n \bs K).$ If $I_0=(0,1),$ we have $I^n \bs K=I_0^n
\bs X,$ and we can conclude using Lemma~\ref{lem:complement}.
\end{proof}

\begin{lemma}[Generalized Mayer-Vietoris inequalities]\label{lem:MVG}
Let $X_1, \ldots , X_m \sub I^n$ be all open or all closed
in $I^n$.
Then
\begin{equation}\label{eq:MVG1}
b_i \left( \bigcup_{1 \le j \le m} X_j \right) \le
\sum_{J \sub \{ 1, \ldots , m \}} b_{i-|J|+1} \left(
\bigcap_{j \in J} X_j \right);
\end{equation}
and
\begin{equation}\label{eq:MVG2}
b_i \left( \bigcap_{1 \le j \le m} X_j \right) \le
\sum_{J \sub \{ 1, \ldots , m \}} b_{i+|J|-1} \left(
\bigcup_{j \in J} X_j \right).
\end{equation}
\end{lemma}

\begin{proof}
The proof is by induction by $m.$ For $m=2,$ the result follows from
the exactness of the Mayer-Vietoris sequence of $(X_1,X_2).$ If the
result is true up to $m-1,$ then define $Y=X_2 \cup \cdots \cup X_m$
and $Z=X_2 \cap \cdots \cap X_m.$ Then~\eqref{eq:MVG1}
and~\eqref{eq:MVG2} follow from the Mayer-Vietoris sequences of
$(X_1,Y)$ and $(X_1,Z)$ respectively.
\end{proof}

\section{The one quantifier case}

In this section, we apply the spectral sequence discussed earlier to
sub-Pfaffian sets defined using a single quantifier. 

\bigskip

Let $\f$ be a Pfaffian chain on a domain $\U\sub \R^{n_0+n_1}$ of
bounded complexity $\gamma$ for $\f,$ such that $I^{n_0+n_1} \sub \U.$
Let $\P$ be a set of Pfaffian functions in the chain $\f$ and
$\Phi(\xx_0,\xx_1)$ be a $\P$-closed formula. 
We denote by $\pi_0$ the canonical projection $\R^{n_0+n_1} \to \R^{n_0}.$

\begin{theorem}[Existential bound]\label{th:exist}
Let $\Phi$ be as above, and let $X=\{(\xx_0,\xx_1) \in I^{n_0+n_1} \mid
\Phi(\xx_0,\xx_1)\}$ and $Y=\pi_0(X).$ Then, if the format of $\Phi$ is
$(n_0+n_1,\alpha,\beta,s),$ we have for all $k\in\N,$
\begin{equation}\label{eq:exist-bd}
b_k(Y) \leq (ks+n_0+kn_1)^N \, 2^{L(L-1)/2}
\,O(\,N \beta + \min(N,L) \,\alpha)^{N+L};
\end{equation}
where $N=n_0+(k+1)n_1$ and $L=(k+1)\ell.$

Moreover, if $X'=\{(\xx_0,\xx_1) \in I^{n_0+n_1} \mid \neg
\Phi(\xx_0,\xx_1)\}$ and $Y'=\pi_0(X'),$ the same bound holds for
$b_k(Y').$
\end{theorem}

\begin{proof}
Theorem~\ref{th:ss} is applicable to the map $\pi_0$ restricted to
$X,$ giving 
\begin{equation*}
b_k(Y) \leq \sum_{p+q=k} b_q(\W^p(X));
\end{equation*}
where $\W^p(X)$ is the $(p+1)$-fold fibered product of $X$ over $Y.$
We can build from the Pfaffian chain $\f$ a chain $F=(\f(\xx_0,\y_0),
\ldots, \f(\xx_0,\y_p))$ of length $(p+1)\ell$ and degree $\alpha$ 
by substituting successively each $\y_j$ for $\xx_1.$ The set $\W^p(X)$
is defined by the following quantifier-free formula in that chain.
\begin{equation*}
(\xx_0,\y_0,\ldots,\y_p) \in I^{n_0}\times I^{(p+1)n_1} \wedge
\Phi(\xx_0,\y_0) \wedge \cdots \wedge \Phi(\xx_0,\y_p).
\end{equation*}
This formula is $\P'$-closed for some $\P',$ and its format is 
\begin{equation*}
(n_0+(p+1)n_1, (p+1)\ell, \alpha, \beta, (p+1)s+2[n_0+(p+1)n_1]);
\end{equation*}
and by Corollary~\ref{cor:Pclosed},
we have 
\begin{equation*}
\begin{split}
b(\W^p(X))\leq (ps+n_0+pn_1)^{n_0+(p+1)n_1} \,
2^{(p+1)\ell((p+1)\ell-1)/2}\, \\  O[(n_0+pn_1) \beta +
\min(n_0+pn_1,p\ell)
\alpha]^{n_0+(p+1)(n_1+\ell)}.
\end{split}
\end{equation*}
Summing the above for $0 \leq p\leq k$ gives~\eqref{eq:exist-bd}.

\bigskip

When considering the case of $Y',$ Theorem~\ref{th:ss} is again
applicable, so we still have $b_k(Y') \leq
\sum_{p+q=k} b_q(\W^p(X')),$ but here $\W^p(X')$ is defined by the formula
\begin{equation*}
(\xx_0,\y_0,\ldots,\y_p) \in I^{n_0}\times I^{(p+1)n_1} \wedge
\neg \Phi(\xx_0,\y_0) \wedge \cdots \wedge \neg \Phi(\xx_0,\y_p);
\end{equation*}
which is neither a $\P$-closed formula nor the negation of one, but we
can reduce to that case in the following way.
Let $I_\e$ be the interval $[\,\e, 1-\e],$ and let $\W^p_\e(X')$ be the
set defined by 
\begin{equation*}
\W^p_\e(X')=\{(\xx_0,\y_0,\ldots,\y_p) \in \mathrm{int}(I^{n_0+(p+1)n_1})\mid 
\neg \Phi(\xx_0,\y_0) \wedge \cdots \wedge \neg \Phi(\xx_0,\y_p)\}.
\end{equation*}
For $\e \ll 1,$ we have $b_q(\W^p(X'))=b_q(\W^p_\e(X'))$ for all $q.$
Since $\W^p_\e(X')$ is given by the negation of a $\P'$-closed formula
of format 
\begin{equation*}
(n_0+(p+1)n_1, (p+1)\ell, \alpha, \beta, (p+1)s+4[n_0+(p+1)n_1]).
\end{equation*}
Corollary~\ref{cor:ch2-Alex} is applicable, and yields the same
asymptotic bound as the bound for $b(\W^p(X)).$
\end{proof}

\begin{corollary}[Universal bound]\label{cor:univ}
Let $\Phi$ be as above and 
let $Z=\{\xx_0 \in I^{n_0} \mid \forall \xx_1 \in
I^{n_1}, \neg \Phi(\xx_0,\xx_1)\}.$ Then, if the format of $\Phi$ is
$(n_0+n_1,\alpha,\beta,s),$ we have for all $k\in\N,$
\begin{equation}\label{eq:univ-bd}
b_k(Z)\leq (n_0+(n_0-k)(s+n_1))^{N^*} \, 2^{L^*(L^*-1)/2}
\,O(\,N^* \beta + \min(N^*,L^*) \,\alpha)^{N^*+L^*};
\end{equation}
where $N^*=n_0+(n_0-k)n_1,$ and $L^*=(n_0-k)\ell.$ %

Moreover, if 
$Z'=\{\xx_0 \in I^{n_0} \mid \forall \xx_1 \in
I^{n_1}, \Phi(\xx_0,\xx_1)\},$ the same bound holds for
$b_k(Z').$
\end{corollary}

\begin{proof}
Let $F$ be the closed set $F=I^{n_0} \bs Z.$ 
We have $b_k(Z)=b_k(I^{n_0} \bs Z),$ so by Lemma~\ref{lem:Alex2}, it
is enough to estimate $b_{n-k-1}(F \cup \fr I^{n_0})$ to estimate
$b_k(Z).$ Let 
\begin{equation*}
X=\{(\xx_0,\xx_1) \in I^{n_0+n_1} \mid \Phi(\xx_0,\xx_1)\} \cup \fr
I^{n_0} \times I^{n_1}.
\end{equation*}
Note that $X$ can be given by a quantifier-free $\P$-closed
formula. Moreover, we have $\pi_0(X)=F \cup \fr I^{n_0}.$ Thus, we can
apply Theorem~\ref{th:exist} to estimate the Betti numbers of $F \cup
\fr I^{n_0},$ and thus of $Z.$ The case of $Z'$ is identical.
\end{proof}

\begin{rem}
The $\P$-closed formula hypothesis is not really necessary here.
Note that similar estimates can be established for a compact set $X$
defined by a formula $\Phi$ that is not $\P$-closed, 
replacing the estimates of Corollary~\ref{cor:Pclosed} by the
Borel-Moore estimates of Theorem~\ref{th:bm}.
\end{rem}

\section{The case of two and more quantifiers}

We will now generalize the results of the previous section to the case
of an arbitrary number of quantifiers. Complementation (with the use
of Alexander duality) allows to apply repeatedly the spectral sequence
argument, by forcing the outer quantifier to be existential, and thus
we can deduce estimates by induction.

\bigskip

Let us fix $n_0 \in \N$ and $\n=(n_1, n_2,\ldots)$ a sequence of positive
integers. For any $\nu\geq 0,$ we let $N(\nu)=n_0+\cdots+n_\nu.$

\bigskip

\begin{definition}\label{df:Ecal} 
For $\n$ as above, we let $\Ecal(\n,n_0,\nu,\ell,\alpha,\beta,s)$ be
the maximum of $b(S),$ where $S \sub I^{n_0}$ is a sub-Pfaffian set
defined as follows: $\Phi$ should be a quantifier-free formula with
format $(N(\nu),\ell, \alpha, \beta, s)$ (see Definition~\ref{df:qf}
and Definition~\ref{df:qf_format}), for some Pfaffian chain $\f$
defined on a domain $\U \supseteq I^{N(\nu)}.$ We assume furthermore
that the semi-Pfaffian set $\{\xx \in I^{N(\nu)} \mid \Phi(\xx)\}$ is
either open or closed in $I^{N(\nu)}.$ Then, if $Q_1, \ldots, Q_\nu$
is a sequence of alternating quantifiers, the set $S$ is defined by
\begin{equation}\label{eq:aux}
S=\{ \xx_0 \in I^{n_0} \mid Q_\nu \xx_\nu \in I^{n_\nu} 
\ldots Q_1 \xx_1 \in I^{n_1}, \Phi(\xx_0,\xx_1,\xx_2, \ldots, \xx_\nu)\}.
\end{equation}
\end{definition}

\medskip

If $\Phi_1, \ldots, \Phi_M$ are quantifier-free formulas defined in
the same Pfaffian chain, defining only open or only closed sets in the
cube $I^{N(\nu)},$ and with the same format $(N(\nu),\ell, \alpha,
\beta, s),$ and if $Q_1, \ldots, Q_\nu$ is a fixed sequence of
alternating quantifiers, we can define for all $1 \leq m \leq M$
\begin{equation*}
S_m=\{ \xx_0 \in I^{n_0} \mid Q_\nu \xx_\nu \in I^{n_\nu} 
\ldots Q_1 \xx_1 \in I^{n_1}, \Phi_m(\xx_0,\xx_1,\xx_2, \ldots, \xx_\nu)\}.
\end{equation*}
Then, we also define $\Ecal^M(\n,n_0,\nu,\ell,\alpha,\beta,s)$
to be the maximum of $b(S),$ where 
\begin{itemize}
\item $S=S_1 \cup \cdots \cup S_M,$ if $Q_\nu=\exists;$ or
\item $S=S_1 \cap \cdots \cap S_M,$ if $Q_\nu=\forall.$
\end{itemize}

\begin{theorem}\label{th:Erec}
For any $\nu \geq 1$ and any values of the other parameters,
the quantity $\Ecal^M(\n,n_0,\nu,\ell,\alpha,\beta,s)$ is bounded by
\begin{equation}\label{eq:Erec}
(4Mn_0)^{n_0} \,  %
\Ecal^{n_0}(\n,(n_\nu+1)n_0,\nu-1,n_0\ell,\alpha,\beta,n_0s).
\end{equation}
\end{theorem}

\begin{proof}
Let $\Phi_1, \ldots, \Phi_M$ be quantifier-free formulas as in
Definition~\ref{df:Ecal}, and consider the sets
\begin{equation*}
X_m=\{ (\xx_0,\xx_1) \in I^{n_0+n_1} \mid %
\forall \xx_{\nu-1} \in I^{n_{\nu-1}}
\ldots Q_1 \xx_1 \in I^{n_1}, \Phi_m(\xx_0,\xx_1, \ldots, \xx_\nu)\}.
\end{equation*}
Let $X=X_1 \cup \cdots \cup X_M,$ and let $S_m=\pi_1(X_m),$ where
$\pi_1$ is the canonical projection $\R^{n_0+n_1} \to \R^{n_0}.$ We
will bound the Betti numbers of $S=S_1 \cup \cdots \cup S_M.$ We can
always reduce to this case by taking complement, in the same way as in
the proof of Corollary~\ref{cor:univ}.

\bigskip

\noindent\textbf{Step one:} {\em spectral sequence argument.}
Note that union and projection commute, so we have $S=\pi_1(X).$ Since
$X$ is open or closed (by Lemma~\ref{lem:open-closed}) we can apply
Theorem~\ref{th:ss} to $\pi_1$ to obtain 
\begin{equation*}
b_k(S) \leq \sum_{p+q=k} b_q(\W^p(X));
\end{equation*}
where $\W^p(X)$ is the corresponding fibered product.
If for all $p,$ $\y^p=(\y_0, \ldots, \y_p)$ denotes a bloc of $(p+1)$
times $n_{\nu}$ variables, and if $m_p=n_0+(p+1)n_\nu$ denotes the
total number of variables in $(\xx_0,\y^p),$ we have 
\begin{equation}\label{eq:wedge}
\W^p(X)=\{(x_0,\y^p) \in I^{m_p} \mid \bigwedge_{j=0}^p \bigvee_{m=1}^M 
(\xx_0,\y_j) \in X_m\}.
\end{equation}

\bigskip

\noindent\textbf{Step two:} {\em Mayer-Vietoris and duality.}
If we define for $1 \leq m \leq M$ and $0 \leq j \leq p,$
\begin{equation}\label{eq:Yj}
Y_j^m=\{(\xx_0,\y^p) \in I^{m_p} \mid (\xx_0,\y_j) \in X_m\};
\end{equation}
we then have from~\eqref{eq:wedge}
\begin{equation}
\W^p(X)= \bigcap_{j=0}^p \bigcup_{m=1}^M Y^m_j.
\end{equation}
We can use the generalized Mayer-Vietoris inequality~\eqref{eq:MVG2}
to transform the intersection above in a union; we obtain
\begin{equation}\label{eq:Wp-ineq}
b_q(\W^p(X)) \leq \sum_{J \sub \{0, \ldots, p\}}
b_{q+|J|-1} \left( \bigcup_{j \in J} \bigcup_{m=1}^M Y^m_j\right)
\end{equation}

\bigskip
Define $\NQ_i$ as the opposite quantifier to $Q_i,$ {\em i.e.}
$\NQ_i=\exists$ if $Q_i=\forall$ and vice-versa, and for all $j$ and
$m,$ let $Z_j^m$ be the subset defined by
\begin{align}
Z_j^m &= \{(\xx_0,\y^p) \in I^{m_p}\mid (\xx_0,\y_j) \not \in X_m\} 
\label{eq:Z1}\\
      &= \{(\xx_0,\y^p) \in I^{m_p}\mid \exists  \xx_{\nu-1} \in
      I^{n_{\nu-1}} %
\cdots  \NQ_1 \xx_1 \in I^{n_1} \,
\neg \Phi_m(\xx_0,\xx_1, \ldots, \xx_\nu)\}. \label{eq:Z2}
\end{align}
By comparing~\eqref{eq:Yj} and~\eqref{eq:Z1}, we notice that
$Z_j^m=I^{m_p} \bs Y_j^m.$ 
For all $J,$ we have 
\begin{equation*}
\bigcap_{j\in J}\bigcap_{m=1}^M Z_j^m =
\left( I^{m_p} \bs \bigcup_{j \in J} \bigcup_{m=1}^M Y_j^m\right);
\end{equation*}
and so by Alexander duality (Lemma~\ref{lem:Alex1} or
Lemma~\ref{lem:Alex2}), we obtain (up to the boundary terms which 
we will neglect for the sake of simplifying the notations)
\begin{equation}\label{eq:YZ}
b_{q+|J|-1} \left( \bigcup_{j \in J} \bigcup_{m=1}^M Y_j^m \right)=
b_{m_p-q-|J|+1} \left( \bigcap_{j \in J} \bigcap_{m=1}^M Z_j^m \right).
\end{equation}
Now, using the other Mayer-Vietoris inequality~\eqref{eq:MVG1}, we have
\begin{equation}\label{eq:ZZ}
b_{m_p-q-|J|+1} \left( \bigcap_{j \in J} \bigcap_{m=1}^M Z_j^m \right)
\leq \sum_{K \sub J\times \{1, \ldots, M\}} 
b_{m_p-q-|J|+|K|} \left( \bigcup_{(j,m) \in K} Z_j^m \right).
\end{equation}

\medskip

Note that the sets $Z_j^m$ are subsets of $I^{m_p},$ and therefore we
have $b_i(Z_j^m)=0$ for $i \geq m_p.$ Thus, we can restrict the sum above
to subsets $K$ such that $m_p-q-|J|+|K| \leq m_p-1,$ which gives
\begin{equation}\label{eq:K}
|K| \leq q+|J|-1 \leq q+p.
\end{equation}
Combining this fact with~\eqref{eq:YZ} and~\eqref{eq:ZZ}, we obtain 
\begin{equation}\label{eq:combined}
b_{q+|J|-1} \left( \bigcup_{j \in J} \bigcup_{m=1}^M Y_j^m \right)
\leq \sum_{K \sub J\times \{1, \ldots, M\}, |K|\leq p+q} 
b_{m_p-q-|J|+|K|} \left( \bigcup_{(j,m) \in K} Z_j^m \right).
\end{equation}

\bigskip

\noindent\textbf{Step three:} {\em combinatorial estimates for
$b_q(\W^p(X))$.}  Let $p$ and $q$ be fixed. We will estimate
$b_q(\W^p(X))$ in terms of $\Ecal.$ For any $J \sub \{0, \ldots, p\}$
and $K \sub \{0, \ldots, p\}\times \{1, \ldots, M\},$ we let
\begin{equation*}
Y_J=\bigcup_{j\in J} \bigcup_{m=1}^M 
Y_j^m, \hbox{ and } Z_K=\bigcup_{(j,m) \in K} Z_j^m.
\end{equation*}
We let $j_0=|J|$ and $k_0=|K|,$ and we denote by $\J(j_0)$ and
$\K(k_0)$ respectively the set of subsets of $\{0, \ldots, p\}$ of
cardinality $j_0$ and the set of subsets of 
$\{0, \ldots, p\}\times \{1, \ldots, M\}$ with cardinality $k_0.$
With these notations, the inequality~\eqref{eq:combined} can be
written as
\begin{equation}\label{eq:sumYJ}
b_{q+j_0-1}(Y_J)
\leq \sum_{k_0=1}^{p+q} \sum_{K \in \K(k_0)}  b_{m_p-q-j_0+k_0} (Z_K).
\end{equation}

\bigskip

Let $K \sub \{0, \ldots, p\}\times \{1, \ldots, M\}$ be fixed and
consider the set
\begin{equation*}
\Sigma(j_0,K)=\{J\in \J(j_0) \mid K \sub J\times\{1, \ldots, M\}\}.
\end{equation*}
Then, for every $J \in \Sigma(j_0,K),$ the term $b_{m_p-q-j_0+k_0}
(Z_K)$ appears, -- when bounding $b_{q+j_0-1}(Y_J),$ -- on the
right-hand side of~\eqref{eq:sumYJ}. Thus, if $\s(j_0,K)$ denotes the
cardinality of $\Sigma(j_0,K),$ we obtain, when
summing~\eqref{eq:sumYJ} over all $J \in\J(j_0),$
\begin{equation*}
\sum_{J\in\J(j_0)} b_{q+j_0-1}(Y_J) \leq \sum_{k_0=1}^{p+q} \sum_{K
\in \K(k_0)} \s(j_0,K) \, b_{m_p-q-j_0+k_0} (Z_K).
\end{equation*}
Since $\Sigma(j_0,K) \sub \J(j_0),$ we have the trivial bound
$\s(j_0,K) \leq 2^{j_0+1} \leq 2^{p+1}.$ Using this in the above
inequality, we get
\begin{equation}\label{eq:2eme}
\sum_{J\in\J(j_0)} b_{q+j_0-1}(Y_J) \leq 2^{p+1} \, \sum_{k_0=1}^{p+q}
\sum_{K \in \K(k_0)} b_{m_p-q-j_0+k_0} (Z_K).
\end{equation}

\bigskip

Recall that from~\eqref{eq:Wp-ineq}, we have 
\begin{align*}
b_q(\W^p(X)) \leq \sum_{j_0=1}^{p+1} \sum_{J\in\J(j_0)}
b_{q+j_0-1}(Y_J).
\end{align*}
Thus, summing~\eqref{eq:2eme} for $1 \leq j_0 \leq p+1,$ we obtain
\begin{equation*}
b_q(\W^p(X)) \leq 2^{p+1} \sum_{j_0=1}^{p+1} \sum_{k_0=1}^{p+q}
 \sum_{K \in \K(k_0)} b_{m_p-q-j_0+k_0} (Z_K).
\end{equation*}
We can change the ordering of the sums in the right hand side to
obtain 
\begin{equation*}
b_q(\W^p(X))  \leq 2^{p+1}  \sum_{k_0=1}^{p+q}
 \sum_{K \in \K(k_0)} \sum_{j_0=1}^{p+1} b_{m_p-q-j_0+k_0} (Z_K)
\end{equation*}
and since obviously $\sum_{j_0=1}^{p+1} b_{m_p-q-j_0+k_0} (Z_K) \leq
b(Z_K),$ we get 
\begin{equation}\label{eq:WbZ}
b_q(\W^p(X)) \leq 2^{p+1} \sum_{k_0=1}^{p+q} \sum_{K \in \K(k_0)}
 b(Z_K).
\end{equation}

\medskip

Now, observe that every set $Z_K$ is given by a union of $|K|$ sets
$Z_j^m,$ which by the formula~\eqref{eq:Z2} are sub-Pfaffian subsets
of $I^{m_p}$ given by an alternation of $\nu-1$ quantifier, starting
with $\exists.$ The formula defining $Z_j^m$ is $\Phi_m(\xx_0, \ldots,
\xx_{\nu-1},\y_j).$ Thus, defining $Z_K$ may require up to $p+q$ such
formulas (since $|K|\leq p+q$) for $1\leq m \leq M$ and $0 \leq j \leq
p.$ Since each such formula involves $s$ Pfaffian functions of degree
at most $\beta$ defined in a chain of length $\ell,$ the set $Z_K$ can
be defined with $(p+q)s$ functions which are defined in a Pfaffian
chain of length $(p+1)\ell.$ By definition of $\Ecal^M,$ it follows
that
\begin{equation*}
b(Z_K)  \leq \Ecal^{|K|}(\n,m_p,\nu-1,(p+1)\ell,\alpha,\beta,(p+q)s)
\end{equation*}
Since $|K|\leq p+q,$ using this estimate in~\eqref{eq:WbZ} gives that
$b_q(\W^p(X))$ is bounded by
\begin{equation}\label{eq:WE}
2^{p+1} [M(p+1)]^{p+q}\,
\Ecal^{p+q}(\n,m_p,\nu-1,(p+1)\ell,\alpha,\beta,(p+q)s);
\end{equation}
since we have (as in Remark~\ref{rem:combi})
\begin{equation*}
\sum_{k_0=1}^{p+q} \sum_{K \in \K(k_0)} 1=
\sum_{k_0=1}^{p+q} \binom{M(p+1)}{k_0} \leq [M(p+1)]^{p+q}.
\end{equation*}

\bigskip

\noindent\textbf{Step four:} {\em summing up.} To bound $b(S),$ all we
need to do now is to sum up~\eqref{eq:WE} for $0 \leq k \leq n_0-1$
and $p+q=k.$ For all the terms in the sums, we have $p+q \leq n_0,$
$p+1\leq n_0$ and $m_p\leq (n_\nu+1)n_0,$ so all the terms
$\Ecal^{p+q}$ from~\eqref{eq:WE} will be bounded by 
\begin{equation*}
\Ecal^{n_0}(\n,(n_\nu+1)n_0,\nu-1,n_0\ell,\alpha,\beta,n_0s).
\end{equation*}
All that remains to be estimated is a term of the form
\begin{equation*}
\sum_{k=0}^{n_0-1} \sum_{p+q=k} 2^{p+1} [M(p+1)]^{p+q}
\end{equation*}
which is clearly bounded by $(2M)^{n_0}\sum_{k=0}^{n_0-1} kn_0^k\leq
(4Mn_0)^{n_0},$ and thus the bound~\eqref{eq:Erec} follows.
\end{proof}

\begin{corollary}\label{cor:EPfaff}
Let $u_\nu=2^\nu n_0n_\nu\cdots n_1$ and $v_\nu=2^{2\nu} n_0^2
n_\nu^2\cdots n_3^2n_2.$ Then, we have
\begin{equation*}
\Ecal(\n,n_0,\nu,\ell,\alpha,\beta,s)\leq
2^{O(\nu u_\nu+\ell^2 v_\nu^2)} s^{O(u_\nu)}
[u_\nu(\alpha+\beta)]^{O(u_\nu+\ell v_\nu)}.
\end{equation*}
\end{corollary}

\begin{proof}
Using the Gabrielov-Vorobjov estimate for arbitrary semi-Pfaffian sets
(Theorem~\ref{th:qf}), %
we can generalize Theorem~\ref{th:exist} and Corollary~\ref{cor:univ}
to obtain, for $\nu=1,$ that
\begin{equation}\label{eq:nu1}
\Ecal^M(\n,n_0,1,\ell,\alpha,\beta,s)\leq
2^{n_0\ell(n_0\ell-1)/2}\, (sM)^{2n_0(n_1+1)}\,
O(n_0n_1(\alpha+\beta))^{n_0(n_1+1+\ell)}.
\end{equation}
(Using the fact that union and existential quantifiers commute, as do
intersections and universal quantifiers.)

\bigskip

Let us now apply Theorem~\ref{th:Erec} inductively. After $i$
iterations, we will denote by $N_i$ the number of free variables,
$s_i$ the number of Pfaffian functions, $\ell_i$ the length of the
Pfaffian chain, $M_i$ the number of sets, and a number $F_i$ so that
\begin{equation*}
\Ecal^M(\n,n_0,\nu,\ell,\alpha,\beta,s) \leq
F_i \,\Ecal^{M_i}(\n,N_i,\nu-i,\ell_i,\alpha,\beta,s_i).
\end{equation*}
We let $N_0=n_0, s_0=s, M_0=M, F_0=1, \ell_0=\ell.$ From
Theorem~\ref{th:Erec}, we know that
\begin{equation*}
N_1=(n_\nu+1)n_0, \, s_1=n_0s, \, M_1=n_0, \, F_1=(4Mn_0)^{n_0}\,
\hbox{ and } \ell_1=n_0\ell.
\end{equation*}
Thus, we obtain for these parameters the following inductions
\begin{align*}
N_{i+1}&=(n_{\nu-i}+1)N_i \\
s_{i+1}&=N_is_i=s N_0\cdots N_i \\
M_{i+1}&=N_i\\
F_{i+1}&=F_i(4M_iN_i)^{N_i}\\
\ell_{i+1}&=N_i\ell+\ell_i=(N_0+\cdots +N_i)\ell.
\end{align*}
From the induction on $N_i,$ we obtain that $N_{i+1}\leq
2N_in_{\nu-i+1},$ and thus for all $i,$ 
\begin{equation*}
N_i \leq 2^i n_0 n_\nu\cdots n_{\nu-i+1};
\end{equation*}
After $\nu-1$ iteration, and applying~\eqref{eq:nu1}, one gets
\begin{align*}
\Ecal^M&(\n,n_0,\nu,\ell,\alpha,\beta,s) \leq F_{\nu-1}
\,\Ecal^{M_{\nu-1}}(\n,N_{\nu-1},1,\ell_{\nu-1},\alpha,\beta,s_{\nu-1})\\
&\leq F_{\nu-1}\, 2^{O(N_{\nu-1}^2\ell_{\nu-1}^2)}\,
(s_{\nu-1}M_{\nu-1})^{2N_{\nu-1}(n_1+1)}\,
O(N_{\nu-1}n_1(\alpha+\beta))^{N_{\nu-1}(n_1+1+\ell_{\nu-1})}.
\end{align*}
Using the notations introduced in the statement of the corollary, 
we obtain
\begin{equation*}
N_{\nu-1}n_1=O(u_\nu), \hbox{ and }
N_{\nu-1}\ell_{\nu-1} =N_{\nu-1} \ell\sum_{j=0}^{\nu-2} N_j
\leq N_{\nu-1}(\nu-1)N_{\nu-2}=O(v_\nu).
\end{equation*}
We also have $2N_{\nu-1}(n_1+1)=O(u_\nu),$ so we can bound
the term $(s_{\nu-1}M_{\nu-1})^{2N_{\nu-1}(n_1+1)}$ by
\begin{equation*}
(sN_0\cdots N_{\nu-3}N_{\nu-2}^2)^{O(u_\nu)}
\leq (sN_{\nu-2}^\nu)^{O(u_\nu)}\leq (su_\nu)^{O(u_\nu)};
\end{equation*}
and we also have
\begin{equation*}
F_{\nu-1}= \prod_{i=0}^{\nu-1} 
(4M_iN_i)^{N_i}\leq
M^{N_0}\,(4N_{\nu-2}N_{\nu-1})^{N_0+\cdots+N_{\nu-1}} \leq M^{n_0}\,
2^{O(\nu u_\nu)}.
\end{equation*}
Using the above estimates, one derives easily an upper-bound for
$\Ecal^M(\n,n_0,\nu,\ell,\alpha,\beta,s)$ in terms of the parameters
$(M,n_0,\ldots, n_\nu,\ell,\alpha,\beta,s),$ and the stated result
follows from the case where $M=M_0=1.$
\end{proof}

\begin{corollary}[Semi-algebraic case]\label{cor:Ealg}
Let $S \sub I^{n_0}$ be as in Definition~\ref{df:Ecal}, for a formula
$\Phi$ having as atoms $s$ polynomials of degree bounded by $d.$ Then,
we have
\begin{equation*}
b(S) \leq [ 2^{\nu^2} d\,s\,n_0 n_1\cdots n_\nu 
]^{O(2^\nu \, n_0n_1\cdots n_\nu)}.
\end{equation*}
\end{corollary}

\begin{proof}
The result follows from the proof of Corollary~\ref{cor:EPfaff},
replacing $\beta$ by $d$ and setting $\alpha=\ell=0.$
\end{proof}

\section{Comparison with quantifier elimination}

We will now compare the results obtained in Corollary~\ref{cor:Ealg}
with the bounds that can be established using quantifier
elimination. Similar comparisons can be made in the Pfaffian case
between Corollary~\ref{cor:EPfaff} and effective cylindrical
decomposition as appear in~\cite{gv:cyldec}
and~\cite{pv:cyldec}. However, we will restrict our attention to the
algebraic case for simplicity.

\bigskip

Let $X$ be a semi-algebraic subspace of $n_0+\cdots+n_\nu$-space
defined by $s$ polynomials of degree bounded by $d,$ and let 
\begin{equation*}
S=\{ \xx_0 \in \R^{n_0} \mid Q_\nu \xx_\nu \in \R^{n_\nu} 
\ldots Q_1 \xx_1 \in \R^{n_1}, (\xx_0,\xx_1,\xx_2, \ldots, \xx_\nu) \in X\}.
\end{equation*}
The set $S$ can be effectively described by a quantifier-free formula
$\Psi(\xx_0)$. The best complexity results for $\Psi$ appear
in~\cite{bpr:qe} (see also~\cite{ba:qe}). We have there
\begin{equation*}
\Psi(\xx_0)=\bigvee_{1\leq i \leq I} \bigwedge_{1 \leq j \leq J_i}
\sign(P_{i,j}(\xx_0))=\e_{i,j};
\end{equation*}
for some family of polynomials $P_{i,j}$ and some sign
conditions $\e_{i,j}\in \{=0,>0,<0\},$ with
\begin{align*}
I & \leq s^{\prod_{i \geq 0} (n_i+1)} d^{(n_0+1)\prod_{i \geq 1}
O(n_i)};\\
J_i & \leq s^{\prod_{i \geq 1} (n_i+1)} d^{\prod_{i \geq 1}
O(n_i)};\\
\deg P_{i,j} & \leq d^{\prod_{i \geq 1} O(n_i)}.
\end{align*}

From this, we can derive the following estimate.

\begin{proposition}\label{prop:bqe}
Let $S$ be as above. Then the sum of the Betti numbers of $S$ verifies
\begin{equation*}
b(S) \leq s^{4n_0(n_0+1)\prod_{i \geq 0} (n_i+1)}
\, d^{O(n_0^2n_1\cdots n_\nu)}.
\end{equation*}
\end{proposition}

\begin{proof}
The set $S$ is defined by a quantifier-free formula involving at most
$\s=J_1\cdots J_I$ polynomials of degree $\d \leq d^{\prod_{i \geq 1}
O(n_i)}.$ By~\cite{gv:qf}, the sum of the Betti numbers is bounded by
$O(\s^2\d)^{n_0}.$ Bounding $\s$ and $\d$ gives the proposition.
\end{proof}

Thus, the estimate in Proposition~\ref{prop:bqe} is better than
Corollary~\ref{cor:Ealg} asymptotically in $\nu,$ but
Corollary~\ref{cor:Ealg} is better for small values of $\nu.$ (For a
fixed value of $\nu,$ Corollary~\ref{cor:Ealg} is better when $n_0$
goes to infinity.)

\bigskip

\begin{rem}
Note that Proposition~\ref{prop:bqe} is a much more general result
than Corollary~\ref{cor:Ealg}, since it is not necessary to make
assumptions on the topology of $X$ or $S,$ nor is it necessary to
restrict ourselves to cubes.
\end{rem}

%% file: aa4.tex
In this chapter, we will give effective estimates for the number of
connected components of the relative closure $(X,Y)_0$ of a
semi-Pfaffian couple $(X,Y).$ This estimate is established first in
the smooth case, by estimating the number of local extrema of the 
distance function $\dist(\cdot, Y_\lambda)$ on $X_\lambda.$
In the singular case, deformation techniques are used to reduce to the
case of smooth hypersurfaces.

\bigskip

Note that the case where $Y=\mpty$ is trivial.  Indeed, this implies
that $(\fr X)_+ =\mpty,$ and since $X$ is assumed to be relatively
compact (see Remark~\ref{rem:relcompact}), $X_\lambda$ is compact for
all $\lambda$ and the number of connected components of $X_0$ is
bounded by the number of connected components of a generic fiber
$X_\lambda$ for $\lambda \ll 1.$ Since $X_\lambda$ is semi-Pfaffian,
Theorem~\ref{th:cells} provides an estimate in that case.

\bigskip

Thus, we'll assume throughout the present chapter that $Y \neq \mpty.$
In the first section, we establish a property that proves the
finiteness of $b_0((X,Y)_0).$ This is used in the second part to
provide the quantitative estimates, first in the smooth case
(Theorem~\ref{th:smooth_bound}) and then in the singular case
(Theorem~\ref{th:sing_bound}). These results are then used in the
third section to give upper-bounds in the fewnomial case.

\section{Finiteness of the number of connected components}

We show here how to reduce the problem of counting the number of
connected components of a limit set to a problem in the semi-Pfaffian
setting.

Let $\Phi$ be the (squared) distance function on 
$\R^n \times \R^n:$
\begin{equation}
\label{eq:map}
\begin{split}
\Phi : \R^n &\times \R^n  \lra \R \\
(x,&y) \mapsto |x-y|^2
\end{split}
\end{equation}

For any $\lambda>0,$ we can define the distance to $Y_\lambda$,
$\Psi_\lambda$ on $X_\lambda$ by:
\begin{equation}
\label{eq:Psi_l}
\Psi_\lambda(x)=\min_{y\in Y_\lambda} \Phi(x,y).
\end{equation}

Define similarly for $x \in \ck{X}:$
\begin{equation}
\label{eq:Psi}
\Psi(x)=\min_{y\in \ck{Y}} \Phi(x,y).
\end{equation}

\begin{theorem}
\label{th:finite}
Let $(X,Y)$ be a semi-Pfaffian couple. Then, there exists $\lambda \ll 1$
such that for every connected component $C$ of $(X,Y)_0,$ we can find
a connected component $D_\lambda$ of the set of local maxima of
$\Psi_\lambda$ such that $D_\lambda$ is arbitrarily close to $C.$
\end{theorem}

\begin{proof}
Let $C$ be a connected component of $(X,Y)_0.$ Note that by definition
of the relative closure, if $x$ is in $C,$ it cannot be in $\ck{Y}.$ So
we must have $\Phi(x,y)>0$ for all $y \in \ck{Y},$ and since $\ck{Y}$ is
compact, we must have $\Psi(x)>0.$ Also, any point in $\fr C$ must be
in $\ck{X},$ but not in $(X,Y)_0.$ So we must have $\fr C \sub \ck{Y},$
hence $\Psi|_{\fr C} \equiv 0.$ This means that the restriction of
$\Psi$ to $C$ takes its maximum inside of $C.$

Choose $x_0\in C$, and let $c=\Psi(x_0)>0$.  For a small $\lambda$,
there is a point $x_\lambda\in X_\lambda$ close to $x_0$ such that
$c_\lambda=\Psi_\lambda(x_\lambda)$ is close to $c$, and is greater
than the maximum of the values of $\Psi_\lambda$ over points of
$X_\lambda$ close to $\fr C$.  Hence the set $\{x\in X_\lambda \mid
\Psi_\lambda(x)\ge c_\lambda\}$ is nonempty, and the connected
component $A_\lambda$ of this set that contains $x_\lambda$ is close
to $C$.  There exists a local maximum $x^*_\lambda\in A_\lambda$ of
$\Psi_\lambda$.  If $D_\lambda$ is the connected component in the set
of local maxima of $\Psi_\lambda,$ it is contained in $Z_\lambda$ and
is close to $C$.
\end{proof}

From the above theorem, we can not only deduce that $(X,Y)_0$ has
finitely many connected components, but also derive effective
estimates.

\section{Bounding the number of connected components}

\subsection{Finding local maxima of the distance function}

We will now show how the number of connected components of the set of
local maxima of $\Psi_\lambda$ that appear in Theorem~\ref{th:finite}
can be estimated when the sets $X_\lambda$ and $Y_\lambda$ are smooth.

\medskip

Define for all $p,$ 
\begin{equation}\label{zeq}
Z^p_\lambda=\{(x,y_0,\dots,y_p)\in W^p_\lambda \mid \Phi(x,y_0)=\dots=\Phi(x,y_p)\},
\end{equation}
where
\begin{equation}\label{weq}
W^p_\lambda=\{(x,y_0,\dots,y_p)\in X_\lambda\times(Y_\lambda)^{p+1} \mid y_i\ne y_j, 0\leq
i<j\leq p\}.
\end{equation}

\begin{lemma}\label{zw}
Assume $(X,Y)$ is a Pfaffian couple such that $X_\lambda$ and $Y_\lambda$ are
smooth for all $\lambda>0.$ For a given $\lambda>0,$ let $x^*$ be a local
maximum of $\Psi_\lambda(x)$.  Then, there exists $0\leq p\le\dim(X_\lambda)$
and a point $z^*=(x^*,y^*_0,\dots,y^*_p)\in Z^p_\lambda$ such that $Z^p_\lambda$ is
smooth at $z^*$, and $z^*$ is a critical point of $\Phi(x,y_0)$ on
$Z^p_\lambda$.
\end{lemma}

\begin{proof}
Since $x^*$ is a local maximum of $\Psi_\lambda(x)$, there exists a point
$y^*_0\in Y_\lambda$ such that $\Phi(x^*,y^*_0)=\min_{y\in
Y_\lambda}\Phi(x^*,y)=\Psi_\lambda(x^*)$.  In particular, $d_y\Phi(x,y)=0$ at
$(x,y)=(x^*,y^*_0)$.  If $(x^*,y^*_0)$ is a critical point of
$\Phi(x,y)$ (this is always the case when $\dim(X_\lambda)=0$) the
statement holds for $p=0$.  Otherwise $d_x\Phi(x,y^*_0)\ne 0$ at
$x=x^*$.  Let $\xi$ be a tangent vector to $X$ at $x^*$ such that
$d_x\Phi(x^*,y^*_0)(\xi)>0$.

Assume that for all $y \in Y_\lambda$ such that $\Phi(x^*,y)=\Psi_\lambda(x^*),$ we
have $d_x\Phi(x,y)(\xi)>0$ when $x=x^*.$ Let $\gamma(t)$ be a curve on
$X_\lambda$ such that $\gamma(0)=x^*$ and $\dot{\gamma}(0)=\xi.$ For all $y
\in Y_\lambda,$ there exists $T_y$ such that for all $0<t<T_y,$ the inequality
$\Phi(\gamma(t),y)>\Phi(x^*,y)$ holds. By compactness of $Y_\lambda,$ this
means we can find some $t$ such that that inequality holds for all $y
\in Y_\lambda.$ Hence, $\Psi_\lambda(\gamma(t))> \Psi_\lambda(x^*),$ which contradicts
the hypothesis that $\Psi_\lambda$ has a local maximum at $x^*.$

Since $x^*$ is a local maximum of $\Psi_\lambda(x)$, there exists a point
$y^*_1\in Y_\lambda$ such that $d_x\Phi(x,y^*_1)(\xi)\le 0$ at $x=x^*$ and
$\Phi(x^*,y^*_1)=\Psi_\lambda(x^*)$.
In particular, $y^*_1\ne y^*_0,\quad d_y\Phi(x^*,y)=0$ at $y=y^*_1$,
and $d_x\Phi(x,y^*_1)\ne d_x\Phi(x,y^*_0)$ at $x=x^*$.
This implies that $(x^*,y^*_0,y^*_1)\in Z^1_\lambda$,
and $Z^1_\lambda$ is smooth at $(x^*,y^*_0,y^*_1)$.
If $(x^*,y^*_0,y^*_1)$ is a critical point of $\Phi(x,y_0)$ on $Z^1_\lambda$
(this is always the case when $\dim(X_\lambda)=1$) the statement holds for $p=1$.
Otherwise $d_x\Phi(x,y^*_0)$ and  $d_x\Phi(x,y^*_1)$
are linearly independent at $x=x^*$.
Since $\dim(X_\lambda)\ge 2$, there exists a tangent vector $\xi$ to $X_\lambda$
at $x^*$ such that $d_x\Phi(x^*,y^*_0)(\xi)>0$ and $d_x\Phi(x^*,y^*_0)(\xi)>0$.
Since $x^*$ is a local maximum of $\Psi_\lambda(x)$,  there exists a point
$y^*_2\in Y_\lambda$ such that $d_x\Phi(x,y^*_2)(\xi)\le 0$ at $x=x^*$ and
$\Phi(x^*,y^*_2)=\Psi_\lambda(x^*)$.
This implies that $(x^*,y^*_0,y^*_1,y^*_2)\in Z^2_\lambda$,
and $Z^2_\lambda$ is smooth at $(x^*,y^*_0,y^*_1,y^*_2)$.
The above arguments can be repeated now for $Z^2_\lambda,\ Z^3_\lambda$, etc.,
to prove the statement for all $p\le\dim(X_\lambda)$.
\end{proof}

\bigskip

Assume now that $X_\lambda$ and $Y_\lambda$ are {\em effectively} non-singular,
{\em i.e.} they are of the following form:
\begin{equation}\label{eq:xy}
\begin{split}
X_\lambda=\{x \in \R^n &\mid p_1(x,\lambda)=\cdots=p_{n-d}(x,\lambda)=0\};\\
Y_\lambda=\{y \in \R^n &\mid q_1(y,\lambda)=\cdots=q_{n-k}(y,\lambda)=0\};
\end{split}
\end{equation}
where, for all $\lambda>0,$ we assume that $d_x p_1 \wedge \cdots \wedge
d_x p_{n-d} \neq 0$ on $X_\lambda$ and that $d_y q_1 \wedge \cdots \wedge
d_y q_{n-k} \neq 0$ on $Y_\lambda.$
In particular, we have $\dim(X_\lambda)=d$ and $\dim(Y_\lambda)=k.$

\bigskip

\begin{rem}
Note that we assume that no inequalities appear in~\eqref{eq:xy}. We
can clearly make that assumption for $Y_\lambda,$ since that set has
to be closed for all $\lambda>0.$ For $X_\lambda,$ we observe the
following: if $C$ is a connected component of $C^p_\lambda,$ the
critical set of $\Phi|_{Z^p_\lambda},$ the function $\Phi$ is constant
on $C.$ If $C$ contains a local maximum for $\Psi_\lambda,$ it cannot
meet $\fr X_\lambda$ because $\fr X_\lambda \sub Y_\lambda.$ Hence, we
do not need to take into account the inequalities appearing in the
definition of $X_\lambda.$
\end{rem}

Let us now define for all $p,$
\begin{equation}\label{eq:theta}
\theta_p: (y_0, \ldots, y_p) \in (Y_\lambda)^{p+1} \mapsto \sum_{0 \leq i < j
\leq p} |y_i-y_j|^2.
\end{equation}
Then, for $X_\lambda$ and $Y_\lambda$ as in~\eqref{eq:xy}, the sets $Z^p_\lambda$ are defined
for all $p$ by the following conditions.
\begin{equation}\label{eq:zp}
\begin{cases}
p_1(x,\lambda)=\cdots=p_{n-d}(x,\lambda)=0; \\
q_1(y_i,\lambda)=\cdots=q_{n-k}(y_i,\lambda)=0, \quad 0\leq i \leq p;\\
\Phi(x,y_i)-\Phi(x,y_j)=0, \quad 0\leq i < j\leq p;
\end{cases}
\end{equation}
and the inequality
\begin{equation}\label{eq:zp2}
\theta_p(y_0, \ldots, y_p)>0.
\end{equation}

Under these hypotheses, we obtain the following bound.

\begin{theorem}[Basic smooth case]
\label{th:smooth_bound}
Let $(X,Y)$ be a semi-Pfaffian couple defined in a domain of bounded
complexity $\gamma,$ such that for all small $\lambda>0,$ $X_\lambda$
and $Y_\lambda$ are effectively non-singular basic sets of dimension
respectively $d$ and $k.$ If the fiber-wise format of $(X,Y)$ is
$(n,\ell,\alpha,\beta,s),$
the number of connected components of $(X,Y)_0$ is bounded by
\begin{equation}\label{eq:smooth_cc}
2 \; \sum_{p=0}^d \V((p+2)n,(p+2)\ell,\alpha,\beta_p,\gamma);
\end{equation}
where $\beta_p=\max\{1+(n-k)(\alpha+\beta-1),
1+(n-d+p)(\alpha+\beta-1)\},$ and $\V$ is defined in~\eqref{eq:defV}.
\end{theorem}

\begin{proof}
According to Theorem~\ref{th:finite}, we can chose $\lambda>0$ such that
for any connected component $C$ of $(X,Y)_0,$ we can find a connected
component $D_\lambda$ of the set of local maxima of $\Psi_\lambda$ such that $D_\lambda$ is
close to $C.$ We see that for $\lambda$ small enough, two connected
components $C$ and $C'$ of $(X,Y)_0$ cannot share the same connected
component $D_\lambda,$ since $D_\lambda$ cannot meet $\ck{Y}$ for $\lambda$ small enough.
Indeed, the distance from $D_\lambda$ to $\ck{Y}$ is bounded from below by
the distance from $D_\lambda$ to $Y_\lambda$, -- which is at least $c_\lambda$, --
minus the distance between $Y_\lambda$ and $\ck{Y}.$ But the latter
distance goes to zero, whereas the former goes to a positive constant
$c$ when $\lambda$ goes to zero.

\medskip

Once that $\lambda$ is fixed, all we need to do is estimate the number of
connected components of the set of local maxima of $\Psi_\lambda.$
According to Lemma~\ref{zw}, we can reduce to estimating the number of
connected components of the critical sets $C^p_\lambda$ of the restriction
$\Phi|_{Z^p_\lambda}$ for $0 \leq \lambda \leq d.$

\medskip
 
For the sake of concision, we will drop $\lambda$ from the notations in
this proof, writing $Z^p$ for $Z^p_\lambda,$ $p_i(x)$ for
$p_i(x,\lambda),$ {\em etc\ldots}

\medskip

A point $z=(x,y_0, \ldots, y_p) \in Z^p$ is in $C^p$ if and only if
the following conditions are satisfied: 
\begin{equation}\label{eq:cp}
\begin{cases}
d_y \Phi(x,y_j)=0, \quad 0 \leq j \leq p; \\
\rk(d_x \Phi(x,y_0), \ldots ,d_x \Phi(x,y_p)) \leq p.
\end{cases}
\end{equation}
 For $X_\lambda$ and $Y_\lambda$ as in~\eqref{eq:xy}, those conditions
become:
\begin{equation}\label{eq:cp2}
\begin{cases}
\rk\{\gr_y q_1(y_i), \ldots, \gr_y q_{n-k}(y_i), \gr_y \Phi(x,y_i)\} \leq n-k,
 \quad 0 \leq i \leq p;\\
\rk\{\gr_x p_1(x), \ldots, \gr_x p_{n-d}(x), \gr_x \Phi(x,y_0), \ldots,
\gr_x \Phi(x,y_0)\} \leq n-d+p.
\end{cases}
\end{equation}
Those conditions translate into all the maximal minors of the
corresponding matrices vanishing. These minors are Pfaffian functions
in the chain used to define $X$ and $Y.$ Their degrees are
respectively $1+(n-k)(\alpha+\beta-1)$ and $1+(n-d+p)(\alpha+\beta-1).$

\bigskip

The number of connected components of $C^p$ is bounded by the number
of connected components of the set $D^p$ defined by the conditions
in~\eqref{eq:zp} and~\eqref{eq:zp2}, and the vanishing of the maximal
minors corresponding to the conditions in~\eqref{eq:cp2}.

\medskip

Let $E^p$ be the set defined by the equations~\eqref{eq:zp}
and~\eqref{eq:cp2}, so that $D^p=E^p \cap \{\theta_p>0\}.$ Then, the
number of connected components of $D^p$ is bounded by the number of
connected components of $E^p$ plus the number of connected components
of $E^p \cap \{\theta_p=\e\}$ for a choice of $\e>0$ small enough.

\medskip

Hence, we're reduced to the problem of estimating the number of
connected components of two 
varieties in $\R^{(p+2)n}$ defined in a Pfaffian chain of degree
$\alpha$ and length $(p+2)\ell.$ Using
Theorem~\ref{th:var}, we obtain the bound~\eqref{eq:smooth_cc}.
\end{proof}

\subsection{Bounds for the singular case}
Let's consider now the case where $X_\lambda$ and $Y_\lambda$ may be
singular. We can use deformation techniques to reduce to the smooth
case.  First, the following lemma shows we can reduce to the case
where $X_\lambda$ is a basic set.

\begin{lemma}
\label{union}
Let $X_1,X_2$ and $Y$ be semi-Pfaffian sets such that $(X_1,Y)$ and
$(X_2,Y)$ are Pfaffian families. Then,  $(X_1 \cup X_2,Y)_0=
(X_1,Y)_0 \cup (X_2,Y)_0.$
\end{lemma}

The proof follows from the definition of the relative closure. 

\medskip

\begin{theorem}[Singular case]
\label{th:sing_bound}
Let $(X,Y)$ be a semi-Pfaffian couple defined in a domain of bounded
complexity $\gamma.$ Assume $X_\lambda$ and $Y_\lambda$ are unions of
basic sets of format $(n,\ell, \alpha, \beta,s).$ If the number of basic
sets in $X_\lambda$ is $M$ and the number of basic sets in $Y_\lambda$
is $N,$ then the number of connected components of $(X,Y)_0$ is
bounded by
\begin{equation}
\label{eq:sing_cc}
2\; MN \; \sum_{p=0}^{n-1} \V((p+2)n,(p+2)\ell,\alpha, \beta^*_p,\gamma);
\end{equation}
where $\beta^*_p=1+(p+1)(\alpha+2\beta-1)$ and $\V$ is defined
in~\eqref{eq:defV}. 
\end{theorem}

\begin{proof}
Again, we want to estimate the number of local maxima of the function
$\Psi_\lambda$ defined in~\eqref{eq:Psi_l}.

By Lemma~\ref{union}, we can restrict ourselves to the case where
$X$ is basic. Let $Y=Y_1 \cup \cdots \cup Y_N,$ where all the sets $Y_i$
are basic. For each basic set, we take the sum of squares of the
equations defining it: the corresponding positive functions, which we
denote by $p$ and  $q_1, \ldots, q_N,$ have degree $2\beta$ in the chain.
Fix $\e_i>0,$ for $0 \leq i \leq N,$ and $\lambda>0,$ and let $\X=\{
p(x,\lambda) =\e_0\}$ and for all $1 \leq i \leq N,$ let
$\Y_i=\{q_i(x,\lambda)=\e_i\}.$ 

Since $Y_\lambda$ is compact, if $x^*$ is a point in $X_\lambda$ such that
$\Psi_\lambda$ has a local maximum at $x=x^*,$ there is a point $y^*$ in
some $(Y_i)_\lambda$ such that $\Phi(x,y)=\Psi_\lambda(x).$ Then, we can find a
couple $(x',y') \in \X_\lambda \times (\Y_i)_\lambda$ close to $(x^*,y^*)$ such
that $\Phi(x',y')$ is a local maximum of
the distance (measured by $\Phi$) from $\X_\lambda$ to $(\Y_i)_\lambda.$

Since for small enough $\e_0, \ldots, \e_N,$ the sets $\X_\lambda$ and
$(\Y_i)_\lambda$ are effectively non-singular hypersurfaces, the number of
local maxima of the distance of $\X_\lambda$ to $(\Y_i)_\lambda$ can be bounded
by~\eqref{eq:smooth_cc}, for appropriate values of the parameters. The
estimate~\eqref{eq:sing_cc} follows.
\end{proof}

\begin{corollary}
Let $\U$ be a fixed domain of bounded complexity, and let $(X,Y)$ be a
semi-Pfaffian couple defined in $\U.$ If $X_\lambda$ and $Y_\lambda$
are unions of basic sets of format $(n,\ell, \alpha, \beta,s),$ where 
$X_\lambda$ is the union of $M$ basic sets and $Y_\lambda$ is the
union of $N$ basic sets, the number of connected components of
$(X,Y)_0$ is bounded by 
\begin{equation*}
MN\, 2^{(n\ell)^2}\, O(n^2(\alpha+\beta))^{(n+1)\ell};
\end{equation*}
for a constant that depends on $\U.$
\end{corollary}

\section{Application to fewnomials}
In this section, we will apply our previous results to the case where
the Pfaffian functions we consider are fewnomials.

Recall from Definition~\ref{df:fewnomials} that we can consider the
restriction of any polynomial $q$ to $\U=\{x_1\cdots x_n \neq 0\}$ as
a Pfaffian function whose complexity depends only on the number of non
zero monomials in $q.$ Fix $\K=\{m_1, \ldots, m_r\} \in \N^n$ a set of
exponents, and let $\ell=n+r,$ and $\f=(f_1, \ldots, f_\ell)$ be the
functions defined by:
\begin{equation}
\label{eq:few_chain2}
f_i(x)=
\begin{cases}
&x_i^{-1} \hbox{ if $1 \leq i \leq n,$}\\
& x^{m_{i-n}} \hbox{ if $i>n.$}\\
\end{cases}
\end{equation}
Then, if $q$ is a polynomial whose non-zero coefficients are in $\K,$
it can be written as a Pfaffian function in $\f$ with degree $\beta=1.$

\medskip

Let now $S \sub \U$ be a bounded semi-algebraic set. 
We can define from $S$ a
semi-Pfaffian family $X \sub \R^n \times \R$ by:
\begin{equation}
\label{eq:X}
X=\{(x,\lambda)\in \U\times \R_+\mid x \in S, \:x_1>\lambda, 
\ldots, x_n>\lambda \}.
\end{equation}

If $S$ is defined by $\K$-fewnomials, we can apply the results from
Theorems~\ref{th:smooth_bound} and~\ref{th:sing_bound} to $X,$ to
obtain a bound on the number of connected components of $\ol{S} \cap
\fr \U.$ Note that from Example~\ref{ex:sub-few}, one can build a
$\K$-fewnomial set $S$ such that $\ol{S}$ is not a $\K$-fewnomial set
(see~\cite{ga:cex}).

\begin{theorem}
\label{th:few}
Let $(X,Y)$ be a semi-Pfaffian couple defined by degree~1 functions in
the chain~\eqref{eq:few_chain2}. 
If $X$ and $Y$ are the union of respectively $M$ and $N$ basic sets,
and letting $q=p+2,$ the number of connected components of $(X,Y)_0$
is bounded by
\begin{equation}
\label{eq:sing_few}
MN \sum_{p=0}^{n-1} 2^{q^2(n+r)^2/2} (6n+6)^{q(3n+2r)}q^{q(n+r)}.
\end{equation}
\end{theorem}

\begin{proof}
This bound is obtained using Theorem~\ref{th:sing_bound} 
for $\alpha=2,$ $\beta=1$ and $\ell=n+r.$ 
\end{proof}

\bigskip

Let $X$ be a semi-Pfaffian family such that for all $\lambda>0,$ the set
$X_\lambda$ is defined by $\K$-fewnomials. By definition of a family, 
$\fr X_\lambda$ is restricted for all $\lambda>0.$
By the results contained in ~\cite{ga:fr-and-cl}, this set is 
semi-Pfaffian set in the same chain,
and the format of $\fr X_\lambda$ can be estimated from the format of
$X_\lambda.$ Applying those results together with those of
Theorem~\ref{th:sing_bound}, we can give estimates for the number of
connected components of $(X, \fr X)_0.$

\begin{theorem}\label{th:dX}
Let $\K$ be fixed and $X$ be a semi-Pfaffian family in the
chain~\eqref{eq:few_chain2}. If $X$ is the union of $N$ basic sets of
format $(n,\ell=n+r, \alpha=2, \beta=1,s),$ the number of connected
components of $X_0=(X,\fr X)_0$ is bounded by
\begin{equation}\label{eq:dX_bound}
N^2 s^{N+rO(n^2)} n^{(n+r)^{n^{O(n^2+nr)}}}.
\end{equation}
\end{theorem}

\begin{proof}
Following~\cite{ga:fr-and-cl}, the set $\fr X_\lambda$ can be defined using
the same Pfaffian chain as $X_\lambda,$ using $N'$ basic sets and functions
of degree at most $\beta',$ where, under the hypotheses above, the
following bounds hold.
\begin{equation*}
\beta' \leq n^{(n+r)^{O(n)}}, \quad \hbox{ and } \quad
N' \leq N s^{N+rO(n^2)} N^{(n+r)^{rO(n)}}.
\end{equation*}
The bound on the number of connected components follows readily.
\end{proof}

\emptyevenpage

%% file: aa5.tex
In this chapter, we consider the following: $\U$ is a domain of
bounded complexity in $\R^n \times \R_+,$ for a Pfaffian chain $\f$ of
length $\ell$ and degree $\alpha,$ and $V$ is a Pfaffian variety in
$\U$ of dimension $d+1,$ with $V_\lambda$ compact for $\lambda>0.$

\medskip

Let $X \sub V$ be a semi-Pfaffian family that is defined by a
$\P$-closed formula $\Phi$ on $V.$ If $\fr X_\lambda=\mpty$ for
$\lambda>0,$ we can consider relative closure of $X,$
$X_0=(X,\mpty)_0,$ which is the Hausdorff limit of the family of
compact sets $X_\lambda.$ We will give in Theorem~\ref{th:rcbnd} an
explicit upper-bound on $b_k(X_0)$ for any $k>0.$ This allows in turn to
establish an upper-bound  for the Betti numbers of any relative
closure $(X,Y)_0,$ even when $Y_\lambda \neq \mpty$
(Theorem~\ref{th:rc-betti}).
In both cases, the bounds depend only on the format of generic fibers,
and are not affected by the dependence in the parameter $\lambda.$

\bigskip

The proof relies on the spectral sequence for closed surjections
developed in Chapter~3. Using triangulation, we construct a
surjection $f^\lambda: X_\lambda \to X_0$ for $\lambda \ll 1.$ Then,
we approximate the corresponding fibered products by semi-Pfaffian
sets. 

\begin{rem}
In the more general setting of o-minimal structures, the result of the
present chapter allow to estimate the Betti numbers of any Hausdorff
limit in a definable family in terms of {\em simple} definable sets
(deformations of diagonals in Cartesian products). This is the point
of view adopted in~\cite{z03}. One can reduce to the one parameter
case using the main result of~\cite{ls:Hausdorff}.
\end{rem}

\section{Constructions with simplicial complexes}

We describe here some constructions that involve simplicial complexes
and PL-maps on them. Using the fact that continuous definable
functions in an o-minimal structure can be triangulated, we will be
able to use these construction in the next section.

\medskip

It is important to note that the constructions done in this section
{\em are not explicit.} They can be achieved using general arguments
from o-minimality, but we do not claim to be able to give an effective
procedure to construct triangulations in the Pfaffian case.
Consequently, these constructions in themselves will not suffice to
establish any Betti number bounds. 

\bigskip

First, let us recall some of the notations used in the discussion of
definable triangulations in section~\ref{sec:geom}. See
Definition~\ref{df:simplex} and Definition~\ref{df:scomplex} for more
details. 
If $a_0, \ldots, a_d$ are affine-independent points in $\R^n,$ we
denote by $\s=(a_0,\ldots, a_d)$ the {\em open} simplex and $\sib=[a_0,
\ldots, a_d]$ the {\em closed} simplex defined by those points.
We say that $K=\{\sib_1, \ldots, \sib_k\}$ is a {\em simplicial
complex} if it is closed under taking faces and for all $i,j$ $\sib_i
\cap \sib_j$ is a common face of $\sib_i$ and $\sib_j.$ We denote by
$|K|$ the {\em geometric realization} of $K.$

\subsection{Retraction on a subcomplex}

Let $K$ be a simplicial complex in $\R^n,$ and let $L \sub K$ be a
subcomplex, {\em i.e.} $L$ is also a simplicial complex.

\medskip

Let $S=\mathrm{st}_K(L)$ be the {\em star} of $L$ in $K,$ {\em i.e.}
the union of all open $\s$ such that $\sib \in K$ and has at least one
vertex in $L.$
We will define a continuous retraction $F$ from $S$ to $L.$ 

\bigskip

If $a_0, \ldots, a_d$ are vertices of $K$ ordered such that $a_0,
\ldots, a_k$ are in $L$ and $a_{k+1}, \ldots, a_d$ are not in $L,$
(for some $k$ such that $0<k<d$), the open simplex $\s=(a_0, \ldots,
a_d)$ is contained in $S$ and we will define $F$ on $\s$ by
\begin{equation}\label{eq:F}
F\left(\sum_{i=0}^d w_i\,a_i\right)=
\frac{1}{\sum_{i=0}^k w_i} \quad \sum_{i=0}^k w_i\,a_i.
\end{equation}

\bigskip

\begin{proposition}\label{prop:F}
The formula~\eqref{eq:F} defines a continuous retraction $F$ from $S$ to
$|L|$ %
that maps all points on the segment $[x,F(x)]$ to $F(x).$
\end{proposition}

\begin{proof}
Let $\s$ be an open simplex appearing in $S.$ Then $\s$ is of the form
$(a_0, \ldots, a_d),$ where the vertices $a_i$ are ordered as above
so that the vertices in $L$ are exactly $a_0, \ldots , a_k$ for some
$0<k<d.$ 

\medskip

Fix $x \in \s,$ and let $s=\sum_{i=0}^k w_i.$ Since all the weights
$w_i$ are positive, the inequality $0<k<d$ implies that $0<s<1.$
Thus, the formula~\eqref{eq:F} clearly defines a continuous function
from $\s$ to $|L|.$

\medskip

Let $y=\theta x+ (1-\theta) F(x)$ for $\theta\in (0,1)$ be a point on
the open segment $(x,F(x)).$ We have $y=\sum_{i=0}^k w'_i \, a_i,$
where
\begin{equation*}
w'_i=
\begin{cases}
\theta w_i + (1-\theta)\Frac{w_i}{s}
& \hbox{ if } 0 \leq i \leq k; \\
\theta w_i & \hbox{ if } k+1 \leq i \leq d.
\end{cases}
\end{equation*}
To prove that $F(x)=F(y),$ we must prove that for all $0 \leq i \leq
k,$ %
\begin{equation*}
\frac{w_i}{\sum_{j=0}^k w_j}=\frac{w'_i}{\sum_{j=0}^k w'_j}.
\end{equation*}

Cross-multiplying, we get the following quantities.
\begin{equation*}
w_i\, \sum_{j=0}^k w'_j= 
w_i\, \sum_{j=0}^k \left(\theta w_j+ (1-\theta)\Frac{w_j}{s} \right)
=w_i \left(1-\theta+\theta s\right);
\end{equation*}
and 
\begin{equation*}
w'_i\, \sum_{j=0}^k w_j= \left(\theta w_i +
(1-\theta)\Frac{w_i}{s}\right)s
=(\theta s + (1-\theta))w_i.
\end{equation*}
The two cross-multiplied quantities are indeed equal, so $F(x)=F(y).$
\end{proof}

For any $x \in S \bs |L|,$ we denote by $\tau_x$ the open segment $(x,
F(x)).$ 

\begin{proposition}\label{prop:tau}
Let $x$ and $y$ be points of $S \bs |L|$ such that $\tau_x$ and $\tau_y$
intersect. Then, we have $\tau_x \sub \tau_y$ or $\tau_y \sub \tau_x.$
\end{proposition}

\begin{proof}
Let $z \in \tau_x \cap \tau_y.$ By Proposition~\ref{prop:F}, we have
$F(x)=F(y)=F(z).$ Thus, $\tau_x$ and $\tau_y$ have one endpoint in common,
and at least one point in common. One must be contained in the other.
\end{proof}

\subsection{Level sets of a PL-function}

Assume now that there exists a continuous function $\pi: |K| \to \R$
with the following properties.
\begin{itemize}
\item $\pi$ is affine on each simplex $\sib$ of $K;$
\item $\pi$ is positive on $K;$
\item $|L|=\pi^{-1}(0).$
\end{itemize}
For all $\lambda >0,$ we will denote by $|K|_\lambda$ the level set
$\pi^{-1}(\lambda).$ We define also
\begin{equation}
\lambda_0=\min\{\pi(a) \mid a \hbox{ is a vertex of } K, a \not \in L\}.
\end{equation}

\begin{rem}
Note that for all $0<\lambda < \lambda_0,$ the level set $|K|_\lambda$
is contained in the star $S.$ Indeed, if $\s=(a_0, \ldots ,a_d)$ is a
simplex that is not in $S,$ we must have $\pi(a_i)\geq \lambda_0$ for all
$i$ since none of the $a_i$ are in $L,$ and $\pi$ being affine on $\sib,$ it
follows that $\pi(x) \geq \lambda_0$ for all $x \in \s.$
\end{rem}

\medskip

We want to describe for $0 < \lambda < \lambda_0,$ the restriction of the
retraction $F$ to the level set $|K|_\lambda.$ We will denote by $F^\lambda$
this restriction. A similar construction is outlined 
in~\cite[Exercise~4.11]{coste:pisa}.

\begin{proposition}\label{prop:homeo}
For all $0<\lambda'<\lambda<\lambda_0,$ there exists a homeomorphism $H: |K|_\lambda
\to |K|_{\lambda'}$ such that $F^\lambda \circ H=F^{\lambda'}.$ 
\end{proposition}

\begin{proof}
Let $x \in |K|_\lambda.$ For any $z \in \tau_x,$ if $z=\theta x+
(1-\theta)F(x),$ we have 
\begin{equation*}
\pi(z)=\theta\pi(x)+(1-\theta)\pi(F(x))=\theta \lambda;
\end{equation*}
since $z$ is in any simplex $\sib$ of $K$ that contains $x$ and $F(x)$ and 
since $\pi$ is affine on the simplices of $K.$
Thus, $z\in |K|_{\lambda'}$ if and only if
$\theta=\lambda'/\lambda,$ and so the map $h$ defined by
\begin{equation}
H(x)=\frac{\lambda'}{\lambda}\,x + \left( 1- \frac{\lambda'}{\lambda}\right) F(x);
\end{equation}
maps $|K|_\lambda$ to $|K|_{\lambda'}.$

\medskip

The map $H$ is certainly injective, since by
Proposition~\ref{prop:tau}, two segments $\tau_x$ and $\tau_y$ cannot
intersect if $x$ and $y$ are two distinct points of $|K|_\lambda.$ It is
also surjective, since for $z \in |K|_{\lambda'},$ it is easy to verify
that the point $x$ defined by
\begin{equation*}
x=\frac{\lambda}{\lambda'}z-\left(\frac{\lambda}{\lambda'}-1\right)F(z);
\end{equation*}
is a point in $|K|_\lambda$ such that $H(x)=z.$

\medskip

The continuity of $H$ follows from the continuity of $F.$ 
Since $H(x) \in \tau_x$ by construction, Proposition~\ref{prop:F}
implies that $F(H(x))=F(x).$
\end{proof}

\begin{proposition}\label{prop:normF}
For $F^\lambda$ defined as above, we have
\begin{equation}\label{eq:limFl}
\lim_{\lambda\to 0} \, \max_{x\in |K|_\lambda} |x-F^\lambda(x)| =0.
\end{equation}
\end{proposition}

\begin{proof}
Let $\s=(a_0, \ldots, a_d)$ be an open simplex contained in $S,$ such
that $\s \not \sub |L|.$ As before, we can assume that the vertices of
$\s$ that are in $L$ are $a_0, \ldots, a_k,$ where $0\leq k<d$
Fix $x = \sum_{i=0}^d w_i \, a_i$ is $\s,$ and let $s=\sum_{i=0}^k
w_i.$ We have 
\begin{equation*}
\sum_{i=k+1}^d w_i=\sum_{i=0}^d w_i-\sum_{i=0}^k w_i=1-s;
\end{equation*}
and 
\begin{equation*}
x-F(x)=\sum_{i=0}^d w_i \, a_i - \Frac{1}{s} \sum_{i=0}^k w_i \, a_i
=\left(1-\frac{1}{s}\right)\left(\sum_{i=0}^k w_i\, a_i \right) +
\sum_{i=k+1}^d w_i\, a_i.
\end{equation*}
By the triangle inequality, we obtain
\begin{equation}\label{eq:norm-ineq}
|x-F(x)|\leq \max_{o \leq i \leq d} |a_i| 
\left(\left| 1-\frac{1}{s} \right|\left(\sum_{i=0}^k w_i\right)  + 
\sum_{i=k+1}^d w_i\right)
=2 (1-s)\max_{0 \leq i \leq d} |a_i|.
\end{equation}
If $x \in |K|_\lambda,$ we have $\pi(x)=\sum_{i=k+1}^d w_i \, \pi(a_i)=\lambda.$
Since $\pi(a_i) \geq \lambda_0$ for all $i \geq k+1,$ it follows that 
\begin{equation}
\lambda=\sum_{i=k+1}^d w_i \, \pi(a_i) \geq \lambda_0 \left(\sum_{i=k+1}^d w_i
\right)= \lambda_0 (1-s).
\end{equation}
It follows that $1-s \leq \Frac{\lambda}{\lambda_0}.$ Combining this 
with~\eqref{eq:norm-ineq}, we obtain
\begin{equation*}
|x-F(x)|\leq 2 \,\frac{\lambda}{\lambda_0} \,\max_{0 \leq i \leq d} |a_i| 
\leq 2 \,\frac{\lambda}{\lambda_0} \,\max\{ |a|, \,a \hbox{ vertex of }K\}.
\end{equation*}
Thus, $|x-F(x)|$ is bounded by a quantity independent of $x$ that goes
to zero when $\lambda$ goes to zero, and the result follows.
\end{proof}

\section{Bounds on the Betti numbers of Hausdorff limits}

Fix $\U$ a domain of bounded complexity $\gamma$ for a Pfaffian chain $\f.$
Let $X$ be a semi-Pfaffian family with compact fibers defined on a
variety $V$ such that $\dim(V)=d+1.$ Assume 
that for all $\lambda\in(0,1),$ $X_\lambda$ is compact, so that
$X_0=(X, \mpty)_0.$ 
The main result of this chapter is the following.

\begin{theorem}\label{th:rcbnd}
Let $X$ be a semi-Pfaffian family with compact fibers as above. If the
format of $X$ is $(n,\ell,\alpha,\beta,s),$ we have for all $0 \leq k \leq d,$ 
\begin{equation}\label{eq:H-betti}
b_k(X_0)\leq \sum_{p=0}^k (10s)^{(p+1)d} \, \V((p+1)n,(p+1)\ell,\alpha,
2\beta, \gamma);
\end{equation}
where $\V$ is defined in~\eqref{eq:defV}. In particular, we have
\begin{equation*}
b_k(X_0) \leq s^{d(k+1)} \, 2^{(k\ell)^2}
O(kn\beta+k\min(n,\ell)\alpha)^{(k+1)(n+\ell)};
\end{equation*}
where the constant depends on $\U.$
\end{theorem}

\begin{rem}
If $X$ is not defined by a $\P$-closed formula, the method of proof is
still valid, and one can still establish bounds on $b_k(X_0),$ using
the Borel-Moore estimates from Chapter~2. In that case, the bound
obtained is 
\begin{equation*}
b_k(X_0) \leq s^{2d(k+1)} \, 2^{2(k\ell)^2}
O(kn\beta+k\min(n,\ell)\alpha)^{2(k+1)(n+\ell)};
\end{equation*}
\end{rem}

In the process of proving Theorem~\ref{th:rcbnd}, we will actually
prove a much more general result. Before stating it, we need
the following notation: for any integer $p,$ we let $\rho_p$ be the
function on $(p+1)$-tuples $(\xx_0, \ldots, \xx_p)$ of points in
$\R^n$ defined by
\begin{equation*}
\rho_p(\xx_0, \ldots,\xx_p)=\sum_{0 \leq i < j \leq p} |\xx_i-\xx_j|^2;
\end{equation*}
Then we will prove the following theorem.

\begin{theorem}
\label{th:1param}
Let $X \sub \R^n\times\R_+$ be a bounded set definable in {\em any}
o-minimal structure, such that the fibers $X_\lambda$ are compact for
all $\lambda>0.$ Let $X_0$ be the Hausdorff limit of those fibers when
$\lambda$ goes to zero. Then, there
exists $\lambda>0$ such that for any integer $k,$ we have
\begin{equation*}%
b_k(X_0) \leq \sum_{p+q=k} b_q(D^p_\lambda(\d));
\end{equation*}
for some $\d>0,$ where the set $D^p_\lambda(\d)$ is the expanded diagonal
\begin{equation*}
D^p_\lambda(\d)=\{(\xx_0, \ldots, \xx_p) \in (X_\lambda)^{p+1} \mid \rho_p(\xx_0,
\ldots, \xx_p)\leq \d\}.
\end{equation*}
\end{theorem}

\begin{rem}
Using results on the definability of Hausdorff limits in o-minimal
structures (see for instance~\cite{ls:Hausdorff}), one can even
generalize the above further: we can estimate in this way the Betti
numbers of any Hausdorff limit of a sequence of compact fibers in a
$p$-parameter definable family.  See~\cite{z03} for more details.
\end{rem}

\subsection{Triangulation of the projection on $\lambda$}

Let $A$ be the closure of $X \cap \{0 <\lambda<1\}.$ By
Theorem~\ref{th:triang-map},
there exists a simplicial complex $K$ such that
$|K| \sub \R^{n+1},$ a subcomplex $L \sub K$ and a homeomorphism $\Phi:
|K| \to A$ such that $\Phi(L)=X_0$ and such that $\pi_\lambda \circ \Phi$
is affine on each simplex of $K.$

\medskip

Denote by $F$ the retraction constructed in the previous section. For
all $\lambda < \lambda_0,$ let $f^\lambda=\Phi^{-1}\circ F^\lambda.$

\begin{proposition}\label{prop:fl}
For all $\lambda < \lambda_0,$ the map $f^\lambda$ is a continuous surjection from
$X_\lambda$ to $X_0.$ Moreover, we have
\begin{equation}\label{eq:normf}
\lim_{\lambda \to 0} \, \max_{x \in X_\lambda} \, |x-f^\lambda(x)|=0.
\end{equation}
\end{proposition}

\begin{proof}
Since $\Phi$ is uniformly continuous, Proposition~\ref{prop:normF}
implies~\eqref{eq:normf}. 
\end{proof}

\bigskip

Define for $p \in \N$ and $\lambda \in (0,\lambda_0),$ 
\begin{equation}\label{eq:W}
W_\lambda^p=\{(\xx_0, \ldots, \xx_p) \in (X_\lambda)^{p+1} \mid 
f^\lambda(\xx_0)=\cdots=f^\lambda(\xx_p)\}.
\end{equation}
From Theorem~\ref{th:ss}
we have for any $\lambda \in (0, \lambda_0),$
\begin{equation}\label{eq:ss}
b_k(X_0) \leq \sum_{p+q=k} b_q(W_\lambda^p).
\end{equation}
Thus, the problem is reduced to estimating the Betti numbers of the
sets $W_\lambda^p.$ The first step in that direction is the following.

\begin{proposition}\label{prop:Whomeo}
For all $0<\lambda'<\lambda<\lambda_0,$ the sets $W_\lambda^p$ and
$W_{\lambda'}^p$ are homeomorphic.
\end{proposition}

\begin{proof}
Let $H: |K|_\lambda \to |K|_{\lambda'}$ be the homeomorphism described in
Proposition~\ref{prop:homeo}. Then, the map $h=\Phi\circ H \circ
\Phi^{-1}$ is a homeomorphism between  $X_\lambda$ and $X_{\lambda'},$ and since
$F^\lambda \circ H= F^{\lambda'},$ we also have $f^\lambda \circ h=f^{\lambda'}.$ It is
then easy to check that the map $h^p: (X_\lambda)^{p+1} \to
(X_{\lambda'})^{p+1}$ defined by
\begin{equation}\label{eq:hp}
h^p(\xx_0, \ldots, \xx_p)=(h(\xx_0), \ldots, h(\xx_p));
\end{equation}
maps $W_\lambda^p$ homeomorphically onto $W_{\lambda'}^p.$
\end{proof}

\subsection{Approximating $W^p$}

For $p \in N$ and $\xx_0, \ldots,\xx_p \in \R^n,$ let $\rho_p$ be the
polynomial 
\begin{equation}
\rho_p(\xx_0, \ldots,\xx_p)=\sum_{0 \leq i < j \leq p} |\xx_i-\xx_j|^2.
\end{equation}
For $\lambda \in (0,\lambda_0),$  $\e>0$ and $\d>0,$ we define the following sets.
\begin{align*}
W_\lambda^p(\e)&=\{(\xx_0, \ldots, \xx_p) \in (X_\lambda)^{p+1} \mid 
\rho_p(f^\lambda(\xx_0), \ldots, f^\lambda(\xx_p)) \leq \e\};\\
D_\lambda^p(\d)&=\{(\xx_0, \ldots, \xx_p) \in (X_\lambda)^{p+1} \mid 
\rho_p(\xx_0, \ldots, \xx_p) \leq \d\}.
\end{align*}
We will use these sets to approximate the sets $W^p.$ Namely, we will
show that for any $p\in \N,$ we can find appropriate values of $\lambda$
and $\d$ such that the Betti numbers of $W_\lambda^p$ and $D_\lambda^p(\d)$
coincide. 

\begin{proposition}\label{prop:Wehomeo}
Let $p \in \N$ be fixed. There exists $\e_0>0,$ such that for all $\lambda
\in (0,\lambda_0)$ and all $0<\e'<\e<\e_0,$ the inclusion $W_\lambda^p(\e') \inc
W_\lambda^p(\e)$ is a homotopy equivalence. In particular, this implies
that 
\begin{equation*}
b_q(W_\lambda^p(\e))=b_q(W_\lambda^p);
\end{equation*}
for all $\lambda \in (0,\lambda_0)$ and all
$\e\in(0,\e_0).$
\end{proposition}

\begin{proof}
First, notice that it is enough to prove the result for a fixed $\lambda
\in (0,\lambda_0),$ since if $0<\lambda'<\lambda<\lambda_0$ are fixed, the map $h^p$
introduced in~\eqref{eq:hp} induces a homeomorphism between
$W_{\lambda}^p(\e)$ and $W_{\lambda'}^p(\e)$ for any $\e>0.$

\medskip

Let us fix $\lambda \in (0,\lambda_0).$ By the generic triviality
theorem (Theorem~\ref{th:generic}), there exists $\e_0>0$ such that
the projection
\begin{equation*}
\left\{ (\xx_0, \ldots, \xx_p, \e) \mid \e\in(0,\e_0) \hbox{ and }
(\xx_0, \ldots, \xx_p) \in W_\lambda(\e)\right\} \mapsto \e;
\end{equation*}
is a trivial fibration.
It follows that for all $0<\e'<\e<\e_0,$ the inclusion $W_\lambda^p(\e')
\inc W_\lambda^p(\e)$ is a homotopy equivalence, and thus the homology
groups $H_*(W_\lambda^p(\e))$ are isomorphic for all $\e \in (0,\e_0).$

\medskip

The sets $W_\lambda^p$ and $W_\lambda^p(\e)$ being compact definable
sets, they are homeomorphic to finite simplicial complexes. This means
that their singular and \v{C}ech homologies coincide, and since
$W_\lambda^p=\cap_{\e>0} W_\lambda^p(\e),$ the continuity property of
the \v{C}ech homology implies that $H_*(W_\lambda^p)$ is the
projective limit of $H_*(W_\lambda^p(\e)).$ Since the latter groups
are constant when $\e \in (0,\e_0),$ the result follows.
\end{proof}

\bigskip

\begin{proposition}\label{prop:deltai}
Let $p\in \N$ be fixed. For $\lambda \ll 1,$ there exist definable
functions $\d_0(\lambda)$ and $\d_1(\lambda)$ such that $\lim_{\lambda
\to 0} \d_0(\lambda)=0,$ $\lim_{\lambda \to 0} \d_1(\lambda)\neq0,$
and such that for all $\d_0(\lambda)< \d'<\d<\d_1(\lambda),$ the
inclusion $D_\lambda^p(\d') \inc D_\lambda^p(\d)$ is a homotopy
equivalence.
\end{proposition}

\begin{proof}
Let $\lambda\in(0,\lambda_0)$ be fixed. By the same local triviality
argument as above, there exists $d_0=0< d_1 < \cdots< d_m <
d_{m+1}=\infty$ such that for all $0\leq i\leq m$ and all $d_i < \d' <
\d < d_{i+1},$ the inclusion $D_\lambda^p(\d') \inc D_\lambda^p(\d)$
is a homotopy equivalence.
When $\lambda$ varies, the values $d_i(\lambda)$ can be taken as
definable functions of the variable $\lambda,$ so by
Lemma~\ref{lem:limit} each has a well-defined if possibly infinite
limit when $\lambda$ goes to zero.  We take
$\d_0(\lambda)=d_j(\lambda),$ where $j$ is the largest index such that
$\lim_{\lambda \to 0} d_j(\lambda)=0,$ and take
$\d_1(\lambda)=d_{j+1}(\lambda).$
\end{proof}

\bigskip

We define for $p \in \N,$ 
\begin{equation}\label{eq:etap}
\eta_p(\lambda)=p(p+1)\,\left(4R\,\max_{x\in X_\lambda} |x-f^\lambda(x)|
+2\,(\max_{x\in X_\lambda} |x-f^\lambda(x)|)^2\right);%
\end{equation}
where, as before, $R$  is a constant such that $X_\lambda \sub B(0,R)$ for all
$\lambda >0.$ By Proposition~\ref{prop:fl}, we have 
\begin{equation*}
\lim_{\lambda \to 0} \eta_p(\lambda)=0.
\end{equation*}

\begin{lemma}\label{lem:inclusions}
For all $\lambda \in (0,\lambda_0),$ $\d>0$ and $\e>0,$ the following inclusions
hold.
\begin{equation*}
D_\lambda^p(\d) \sub W_\lambda^p(\d+\eta_p(\lambda)), \hbox{ and } 
W_\lambda^p(\e) \sub D_\lambda^p(\e+\eta_p(\lambda)).
\end{equation*}
\end{lemma}

\begin{proof}
Let $m(\lambda)=\max_{x\in X_\lambda} |x-f^\lambda(x)|.$
For any $\xx_i, \xx_j$ in $X_\lambda,$ the triangle inequality gives
\begin{align*}
|f^\lambda(\xx_i) - f^\lambda(\xx_j)|^2 &\leq [|f^\lambda(\xx_i)-\xx_i| + |\xx_i - \xx_j| +
|\xx_j - f^\lambda(\xx_j)|]^2\\
&\leq [|\xx_i-\xx_j|+2m(\lambda)]^2\\
&\leq |\xx_i-\xx_j|^2 + 8R \, m(\lambda)+
4 m(\lambda)^2.
\end{align*}
Summing this inequality for all $0 \leq i < j \leq p,$ we obtain that
for any $\xx_0, \ldots, \xx_p$ in $X_\lambda,$
\begin{equation*}
\rho_p(f^\lambda(\xx_0), \ldots, f^\lambda(\xx_p)) \leq \rho_p(\xx_0, \ldots,
\xx_p)+\eta_p(\lambda). 
\end{equation*}
The first inclusion follows easily from this inequality. The
second inclusion follows from a similar reasoning.
\hbox{ }
\end{proof}

\begin{proposition}\label{prop:approx}
For any $p\in \N,$ there exists $\lambda \in (0,\lambda_0),$ $\e\in (0,\e_0)$
and $\d>0$  such that
\begin{equation}
H_*(W_\lambda^p(\e)) \cong H_*(D_\lambda^p(\d)).
\end{equation}
\end{proposition}

\begin{proof}
Let $\d_0(\lambda)$ and $\d_1(\lambda)$ be the functions defined in
Proposition~\ref{prop:deltai}. Since the limit when $\lambda$ goes to zero
of $\d_1(\lambda)-\d_0(\lambda)$ is not zero, whereas the limit of $\eta_p(\lambda)$
is zero, we can choose $\lambda >0$ such that
$\d_1(\lambda)-\d_0(\lambda)>2\eta_p(\lambda).$ Then, we can choose $\d'>0$ such that
$\d_0(\lambda)<\d'<\d'+2\eta_p(\lambda)< \d_1(\lambda).$ Taking a smaller $\lambda$
if necessary, we can also assume that $\d'+3\eta_p(\lambda)<\e_0.$

\medskip

Let $\e=\d'+\eta_p(\lambda),$ $\d=\d'+2\eta_p(\lambda)$ and $\e'=\d'+3\eta_p(\lambda).$ 
From Lemma~\ref{lem:inclusions}, we have the following sequence of inclusions;
\begin{equation*}
D_\lambda^p(\d')\stackrel{i}{\inc} W_\lambda^p(\e)
\stackrel{j}{\inc} D_\lambda^p(\d) \stackrel{k}{\inc} W_\lambda^p(\e').
\end{equation*}

By the choice of $\e,\e'$ and $\lambda,\lambda',$ the inclusions $k\circ j$ and
$j\circ i$ are homotopy equivalences. The resulting diagram in
homology is the following;
\begin{equation*}
\xymatrix{
H_*(D_\lambda(\d'))\ar[rr]^{(j\circ i)_*}_{\cong} 
\ar[dr]^{i_*}
& & H_*(D_\lambda(\d))\ar[dr]^{k_*} & \\
& H_*(W_\lambda(\e))\ar[rr]^{(k\circ j)_*}_{\cong} 
\ar[ur]^{j_*}
& & H_*(W_\lambda(\e'))
}
\end{equation*}
Since $(j \circ i)_*=j_*\circ i_*,$ is an isomorphism, $j_*$ must be
surjective, and similarly, the fact that $(k \circ j)_*=k_*\circ j_*$
is an isomorphism implies that $j_*$ is injective.  Hence, $j_*$ is an
isomorphism between $H_*(W_\lambda(\e))$ and $H_*(D_\lambda(\d)),$
as required.
\end{proof}

\subsection{Proof of Theorem~\ref{th:rcbnd}}

Recall that from the spectral sequence inequality~\eqref{eq:ss}, all
we need to do to bound $b_k(X_0)$ is to give explicit estimates on the
Betti numbers of the sets $W_\lambda^p$ for all $0 \leq p \leq k.$

\bigskip

Let $p \in \N$ be fixed, and choose $\e, \lambda$ and $\d$ as in
Proposition~\ref{prop:approx}. Since $\e<\e_0,$ the Betti numbers of
$W_\lambda^p(\e)$ and $W_\lambda^p$ are the same. Thus, we are reduced to
estimating $b_q(D_\lambda^p(\d)).$ 
This set is a semi-Pfaffian subset of $V_\lambda^{p+1}$ defined by a
$\mathcal{Q}$-closed formula, where $\mathcal{Q}$ is a set of
$s(p+1)+1$ Pfaffian functions in $n(p+1)$ variables, of degree bounded
by $\beta$ in a chain of length $(p+1)\ell$ and degree $\alpha.$
More explicitly, the chain under consideration is 
\begin{equation*}
\boldsymbol{f}^p=(f_1(\xx_0,\lambda), \ldots, f_\ell(\xx_0,\lambda), \ldots, f_1(\xx_p,\lambda),
\ldots, f_\ell(\xx_p,\lambda)), 
\end{equation*}
where $\lambda$ is kept constant. 
It follows from Theorem~\ref{th:Pclosed} that
\begin{equation*}
b_q(D_\lambda^p(\d)) \leq b(D_\lambda^p(\d)) \leq
(10s)^{d(p+1)} \, \V((p+1)n,(p+1)\ell,\alpha,
2\beta, \gamma).
\end{equation*}
The bound~\eqref{eq:H-betti} follows.
\hfill $\Box$

\section{Betti numbers of a relative closure}

In addition to the bound of Theorem~\ref{th:rcbnd}, the techniques
developed in the present chapter allow us to give a rough estimate
for the Betti numbers of a relative closure of a semi-Pfaffian
couple. In that case, -- as for the Hausdorff limits, -- the Betti numbers
of the limit depends on the format of the fibers, but not on the
families' dependence on the parameter $\lambda.$

\begin{theorem}\label{th:rc-betti}
Let $(X,Y)$ be a semi-Pfaffian couple. Then, the Betti numbers of its
relative closure $(X,Y)_0$ can be bounded in terms of the format of
the fibers $X_\lambda$ and $Y_\lambda$ for $\lambda \ll 1.$
\end{theorem}

\begin{proof}
Let $\d_0>0$ and define
\begin{equation*}
K_0=\{x \in \ck{X} \mid \dist(x,\ck{Y}) \geq \d_0\}.
\end{equation*}
Recall that $(X,Y)_0=\{x \in \ck{X} \mid \dist(x,\ck{Y})>0\},$ so
$K_0$ is a compact subset of $(X,Y)_0.$ For $\d_0 \ll 1,$ the Betti
numbers of $K_0$ and $(X,Y)_0$ coincide.
Let $\d(\lambda)$ be any definable function such that $\lim_{\lambda \to 0}
\d(\lambda)=\d_0,$ and let 
\begin{equation*}
K=\{(x,\lambda) \in X \mid \dist(x,Y_\lambda) \geq \d(\lambda)\}.
\end{equation*}
The set $K$ is a definable family with compact fibers for $\lambda>0.$ The
Hausdorff limit of this family when $\lambda$ goes to zero is $K_0.$ 
By Theorem~\ref{th:1param}, we have
for all $k\in\N$ and all $\lambda \ll 1,$ 
\begin{equation*}
b_k(K_0)\leq \sum_{p+q=k} b_q(D_\lambda^p(\eta));
\end{equation*}
where $\eta>0$ is fixed and
\begin{equation*}
D_\lambda^p(\eta)=\{(\xx_0, \ldots, \xx_p) \in (K_\lambda)^{p+1} \mid 
\rho_p(\xx_0, \ldots, \xx_p) \leq \eta\}.
\end{equation*}

\medskip

Consider now the Cartesian product
\begin{equation*}
T_\lambda^p=\{(\xx_0, \ldots, \xx_p, \y_0, \ldots, \y_p) \mid \xx_i \in X_\lambda,
\y_j \in Y_\lambda\}.
\end{equation*}
We can consider a cylindrical cell decomposition of $T_\lambda^p$ that
would be compatible with the subsets $\{\rho_p(\xx_0, \ldots, \xx_p)
=\eta\}$ and $\{|\xx_i-\y_i|=\d(\lambda)\}.$ The number of cells in such a
decomposition depends only on $p$ and on the formats of $X_\lambda$ and
$Y_\lambda,$ and the projection of this decomposition on the variables
$(\xx_0, \ldots, \xx_p)$ is compatible with $D_\lambda^p.$ Thus, the total
number of cells in the decomposition of $T_\lambda^p$ bounds $b(D_\lambda^p),$
and an upper-bound on $b_k(K_0)=b_k((X,Y)_0)$ follows. Explicit
bounds, which would be doubly exponential in $kn,$ can be derived
from~\cite{gv:cyldec,pv:cyldec}.
\end{proof}

\begin{rem}
The explicit bound that would be obtained by the above method is very
bad (doubly exponential). In particular, the bounds obtained are worse
than those obtained in Chapter~4 in the case where $k=0.$ 
Better estimates, that coincide with Chapter~4 when $k=0,$ are work in
progress~\cite{z03:rc}. 
\end{rem}

%% file: spectral.tex
\chapter[Spectral sequences]{Spectral sequence associated to a filtered chain complex}

\section{Homology spectral sequence}

We consider here only first quadrant homology spectral sequences.
Assume $\{E^r_{p,q}\}$ are modules over some ring $R,$ which are non
zero only for $p,q,r \geq 0.$

This is a spectral sequence if for all $r,p,q$ there is a differential 
$d^r_{p,q}: E^r_{p,q} \lra E^r_{p-r,q+r-1}$ such that $d^r\circ
d^r=0,$ and such that 
\begin{equation}\label{eq:Er1}
E^{r+1}_{p,q}=\Frac{\ker (d^r_{p,q})}
{d^r_{p+r,q-r+1} \left(E^r_{p+r,q-r+1}\right)}.
\end{equation}

Note that for $r>p$ we have $\ker (d^r_{p,q})=E^r_{p,q}$ since the
image is a term that lies in $p<0.$ Similarly, for $r>q+1,$ the
module $E^r_{p+r,q-r+1}$ is zero, hence $d^r_{p+r,q-r+1}
\left(E^r_{p+r,q-r+1}\right)$ is zero too. It follows
from~\eqref{eq:Er1} that for all $r>\max(p,q+1),$ we must have
$E^{r+1}_{p,q}=E^r_{p,q}.$
We denote by $E^\infty_{p,q}$ the 
term at which $E^r_{p,q}$ stabilizes.

\bigskip

Let $H$ be a chain complex with a an increasing filtration
\begin{equation*}
0 \sub F_0H \sub F_1H \sub \cdots \sub F_pH \sub \cdots \sub H
\end{equation*}
such that $\cup_p F_pH=H.$ We say that $E^r_{p,q}$ converges to $H$ if
for all $p$ and $q$ we have 
\begin{equation}\label{eq:Er2}
E^\infty_{p,q}\cong\Frac{F_pH_{p+q}}{F_{p-1}H_{p+q}}
\end{equation}
It is usually denoted by $E^r_{p,q} \Rightarrow H.$

\section{Sequence associated with a filtered complex}

Let $C$ be a chain complex such that $C_n=0$ for $n<0.$ Denote by
$d_n$ the differential from $C_n$ to $C_{n-1}.$ Assume that there is a
filtration of $C$ by subcomplexes $\{F_p\}_{p\in\N}$ such that $F_p$
is increasing and $\cup_p F_p=C.$

\medskip

For $n \in \N,$ denote by $Z_n=\ker(d_n)$ the cycles in $C_n$ and by
$B_n=d_{n+1}C_{n+1}$ the boundaries. Let $H_n=Z_n/B_n.$ We can use the
filtration $F_p$ to approximate $Z_n$ and $B_n,$ in the following
fashion: define, for all $r\geq0,$
\begin{equation*}
A^r_p=\{c \in F_p \mid dc \in F_{p-r}\}.
\end{equation*}
The elements of $A^r_p$ are cycles 'up to $F_{p-r}.$' For $r>p,$ they
are really cycles.

\bigskip

One then defines approximate cycles $Z^r_p=A^r_p/F_{p-1}$ and
$B^r_p=dA^{r-1}_{p+r-1}/F_{p-1}.$ Note that the indices are chosen so
that both are submodules of $F_p/F_{p-1},$ and we have the following
inclusions: 
\begin{equation*}
0 \sub B^0_p \sub B^1_p \sub \cdots \sub B^\infty_p \sub Z^\infty_p
\sub \cdots \sub Z^1_p \sub Z^0_p=F_p/F_{p-1}.
\end{equation*}
Here, $B^\infty_p$ denotes the increasing union of the modules $B^r_p$
and $Z^\infty_p$ is the decreasing intersection of the modules $Z^r_p.$
We then define
\begin{equation*}
E^r_{p,q}=(Z^r_p)_{p+q}/(B^r_p)_{p+q}.
\end{equation*}

\begin{theorem}\label{th:filtr-ss}
With the above definitions, $E^r_{p,q}$ is a spectral sequence that
converges to $H(C).$
\end{theorem}

See~\cite{user} for a proof.

\begin{corollary}\label{cor:ss-bound}
Let $C$ be a chain complex such that $C_n=0$ for $n<0,$ such that
there is an increasing exhaustive filtration of $C$ and let
$E^r_{p,q}$ be the associated homology spectral sequence.
Then, we have for all $n,$
\begin{equation*}
\rk \,H_n(C) \leq \sum_{p+q=n} \rk (E^1_{p,q}).
\end{equation*}
\end{corollary}

\begin{proof}
From~\eqref{eq:Er1}, it is clear that $\rk E^r_{p,q}$ is a decreasing
sequence when $p$ and $q$ are fixed and $r$ increases. Thus, we have
$\rk E^\infty_{p,q} \leq \rk E^1_{p,q}.$ Since $E^r_{p,q} \Rightarrow
H(C),$ it follows from~\eqref{eq:Er2} that
\begin{equation*}
H_n(C) \iso \bigoplus_{p+q=n} E^\infty_{p,q}.
\end{equation*}
The result follows easily.
\end{proof}

\section{Alexander cohomolgy and applications}

This section contains a short list of auxiliary results that play a
role in the construction of the filtration that gives rise to the spectral
sequence of Theorem~\ref{th:spectral-f}.
We refer the reader to Chapter~6 of~\cite{spanier} for more general
statements, proofs and additional details.


\bigskip

In this section, $\HB^*$ denotes the Alexander cohomology. 


\begin{definition}
A topological space $X$ is said to be {\em homologically connected} if
for all $x\in X$ and all neighbourhood $U$ of $x,$ there exists a
neighbourhood $V \sub U$ such that the map $\HT_q(V) \to \HT_q(U)$
given by inclusion is trivial for all $q.$
\end{definition}

\begin{definition} 
A topological space $X$ is said to be {\em locally contractible} if
for all $x\in X$ and all neighbourhood $U$ of $x,$ there exists a
neighbourhood $V \sub U$ such that $V$ can be deformed to $x$ in $U.$
\end{definition}

If $X$ is locally contractible, it is homologically connected. In
particular, all sets that are definable in some o-minimal structures
are locally contractible.

\begin{proposition}\label{prop:Hbar}
Let $X$ be homologically connected. We have $\HB^*(X) \iso H^*(X),$
where $H^*(X)$ is the singular cohomology of $X.$
\end{proposition}

\begin{theorem}[Vietoris-Begle]\label{th:VB}
Let $F: A \to B$ be a closed, continuous surjection between
paracompact Hausdorff spaces. If for all $q,$ we have
$\HBT^q(F^{-1}y)=0$ for all $y \in B,$ the map $F^*: \HB^*(B) \to
\HB^*(A)$ is an isomorphism.
\end{theorem}

See~\cite[p. 344]{spanier} for a proof. Example~16 on the same page
shows that the theorem does not hold if $F$ is not closed.

\emptyevenpage

%% file: bibliographie.tex
\vfill

\hbox{ }